\documentclass[12pt,a4paper,oneside]{article}
\usepackage{a4wide}
\usepackage[T1]{fontenc}
\usepackage[utf8]{inputenc}
\usepackage{amsmath}
\usepackage{times}
\usepackage{graphicx}
\usepackage[final]{pdfpages}
\usepackage{longtable}
\usepackage{float}
\usepackage{amssymb,bm}
\usepackage{url}
\usepackage{xcolor}
\usepackage{mathtools}
\usepackage{tcolorbox}
\usepackage{soul}
\definecolor{pigpink}{HTML}{FDD7E4}
\definecolor{lcyan}{HTML}{E0FFFF}
\definecolor{mint}{HTML}{98FF98}

\usepackage{marginnote}
\usepackage{soul}

\usepackage{amsthm}
\newtheorem{theorem}{Theorem}[section]

\newtheorem{remark}[theorem]{Remark}

\newtheorem{algorithm}[theorem]{Algorithm}

\newcommand{\su}{\bm{\alpha}}
\newcommand{\suf}{\bm{\alpha},k}
\newcommand{\iloc}{\mathcal{I}_{\theta}^{loc}}
\newcommand{\aloc}{\mathcal{A}_{c}^{i,loc}}
\newcommand{\iiso}{\mathcal{I}_{\theta}}

\usepackage{indentfirst}

\title{Global phase-amplitude description of oscillatory dynamics via the parameterization method}
\author{Alberto P\'erez-Cervera$^{1,2}$, Tere M-Seara$^{1}$ and Gemma Huguet$^{1}$\\
	\parbox{13.5cm}{
		\small
		\begin{itemize}
			\item[$^1$]
			Departament de Matem\`atiques, Universitat Polit\`ecnica de
			Catalunya, Avda. Diagonal 647, 08028 Barcelona. 
			\item[$^2$]
			Department of Nonlinear Dynamics and Complex Systems, Institute of Computer Science of the Czech Academy of Sciences, Pod Vodárenskou věží 271/2
			182 07 Prague.
		\end{itemize}
}}

\begin{document}
	%----------------------------------------------------------------------------
	\date{}
	\maketitle
	\vspace{-0.3cm}
	\noindent \textbf{Corresponding author:} Alberto P\'erez-Cervera,
	\texttt{perez@cs.cas.cz} \\
	
	\noindent \textbf{Keywords:} Phase-Amplitude variables, non-linear oscillators, Phase and Amplitude Response Functions, isochrons, isostables,  parameterization method. \\
	
	\noindent \textbf{MSC2000 codes:} 37D10, 92B25, 65P99, 37N25 
	
	%\noindent %\textbf{Abbreviated title:} To be filled \\

	\section*{Abstract}
	
	 In this paper we use the parameterization method to provide a complete description of the dynamics of an $n$-dimensional oscillator beyond the classical phase reduction. The parameterization method allows, via efficient algorithms, to obtain a parameterization of the attracting invariant manifold of the limit cycle in terms of the phase-amplitude variables. The method has several advantages. It provides analytically a Fourier-Taylor expansion of the parameterization up to any order, as well as a simplification of the dynamics that allows for a numerical globalization of the manifolds. Thus, one can obtain the local and global isochrons and isostables, including the slow attracting manifold, up to high accuracy, which offer a geometrical portrait of the oscillatory dynamics. Furthermore, it provides straightforwardly the \textit{infinitesimal Phase and Amplitude Response Functions}, that is, the extended infinitesimal Phase and Amplitude Response Curves, which monitor the phase and amplitude shifts beyond the asymptotic state. Thus, the methodology presented yields an accurate description of the phase dynamics for perturbations not restricted to the limit cycle but to its attracting invariant manifold. Finally, we explore some strategies to reduce the dimension of the dynamics, including the reduction of the dynamics to the slow stable submanifold. 
	 We illustrate our methods by applying them to different three dimensional single neuron and neural population models in neuroscience.

%\section*{Highlights}
%
%We extend the applications of the parameterization method to compute the full set of phase-amplitude coordinates for high dimensional oscillators. 
%We  use the Floquet normal form and automatic differentiation techniques to drastically reduce the computational cost of the required calculations. Our methods provide an analytical expression for the local isochrons, isostables, and infinitesimal Phase and Amplitude Response Functions in a neighbourhood of the limit cycle, while allow for the globalization of these objects and functions to the full basin of attraction.  We illustrate our methodology by applying it to relevant single neuron and neural population models in neuroscience. Moreover, we perform a perturbation study of a single neuron model and study the scope of validity of different dynamical reductions, namely the slow-manifold reduction and the phase reduction. Our results provide an efficient methodology that allows for a geometrical understanding of the dynamics of high dimensional oscillators.
%
%\newpage	
	
%&&&&&&&&&&&&&&&&&&&&&&&&&&&&&&&&&&&&&&&&&&&&&&&&&&&&&&&&&&&&&&&&&&&&&&&

\section{Introduction}

Oscillations are ubiquitous in a wide range of physical,  chemical  and  biological processes \cite{kuramoto2003chemical, winfree2001geometry, strogatzbook, PIK01}. 
%\cite{rapp1987so}. Some classical examples are the cardiac pacemakers \cite{michaels1987mechanisms}, the flashing fireflies \cite{ermentrout1991adaptive} or the circadian rhythms \cite{liu1997cellular}, and 
Of particular interest for this paper are the neural oscillators ranging from single neuron  \cite{izhikevich2007} to neural populations \cite{buzsaki2006rhythms}. From the mathematical perspective, oscillations correspond to attracting limit cycles in the phase space whose dynamics can be described by a single variable: the phase. The description of the dynamics by means of the phase variable, assuming that trajectories remain close to the limit cycle, is known as the phase reduction. This reduction has been extensively used to study weakly coupled oscillators because it simplifies the dynamics of a high-dimensional oscillator to a simple one-dimensional dynamical system \cite{hoppensteadt2012, ErmentroutTerman2010, ErmentroutKopell91, winfree1967biological}.

The phase variable can be extended out of the limit cycle under generic conditions via the concept of asymptotic phase and isochrons \cite{winfree1974patterns}. When the limit cycle is hyperbolic, it is known that the orbit of each point on the basin of attraction approaches asymptotically the orbit of a point on the limit cycle \cite{hirsch1970stable, guckenheimer1975}. Therefore, we can assign to any point an asymptotic phase corresponding to the phase of the point on the limit cycle whose orbit approaches to. The isochrons are thus defined as the set of points having the same asymptotic phase \cite{winfree2001geometry}. Brief external perturbations displace trajectories from one isochron to another, thus causing a phase shift.
%, the phase variable, besides simplifying dynamics, provides a geometric portrait of the synchronization properties of a given oscillator \cite{oprisan2002influence, winfree2001geometry}. 
A useful tool characterizing the response of a given oscillator to brief perturbations is the phase response curve (PRC) \cite{ErmentroutTerman2010, schultheiss2011phase, brown2004phase}. The PRC quantifies the phase change of an oscillator due to perturbations acting at different phases of the cycle. If the strength of the perturbation is weak, the PRC can be approximated by its first order term, the infinitesimal PRC (iPRC). Nevertheless, the validity of the phase reduction and the iPRC is limited by different factors. For instance, if the time between stimuli is short, or the contraction rate to the limit cycle is weak, or the amplitude of the perturbation is strong, the trajectories may remain away from the limit cycle and the classical phase reduction approach fails.

Due to the widespread presence of oscillatory dynamics, there are many cases in which the limitations of the phase reduction become relevant. Examples can be found in the field of circadian rhythms \cite{WilsonE19}, in control theory \cite{monga2019optimal} or in the study of the role of brain rhythms in cognitive functions \cite{buzsaki2006rhythms, schultheiss2011phase,  canavier2015phase, tiesinga2010, PerezSH20}. These examples, amongst many others, motivate the growing mathematical research trying to overcome the limitations inherent to the phase reduction. Amongst the different approaches taken, we find the extension of the phase variable to the basin of attraction of the limit cycle through the computation of the global isochrons \cite{osinga2010continuation, Detrixheetal16}, the computation of PRCs beyond the linear approximation given by the iPRC \cite{suvak2010quadratic, takeshita2010higher, PerezCervera2019, OprisanPC04}, or generalizations of the phase reduction approach for perturbations having specific frequency features \cite{park2016weakly, pyragas2015phase}.

An alternative strategy, which has gained a lot of interest in the recent years, consists in introducing extra variables which describe the dynamics along the directions transversal to the limit cycle. The addition of these variables allows to describe the effects of the perturbations acting on trajectories away from the limit cycle. Over the last decade, several groups have proposed different ways to define these extra variables. In \cite{wedgwood2013phase}, the authors consider orthogonal directions to the limit cycle. However, the most extended approach consists in defining a coordinate system in which the extra transverse coordinates move along the isochrons with a decay rate given by the Floquet exponents of the limit cycle. In the work \cite{mauroy2016global}, these coordinates were called isostables and computed using the Koopman operator. They are defined as the level sets of the slowest decaying Koopman eigenfunction. This approach has been adapted for stable fixed points \cite{mauroy2013isostables} and limit cycles \cite{moehliswilsonpre2016, shirasaka2017phase}. Furthermore, the Koopman operator led to new approaches to globalize isochrons, isostables and compute high order perturbations of limit cycles \cite{mauroy2018global}. Recent works also adapted the concept of isostable coordinate to extend the phase dynamics in a neighborhood of the limit cycle \cite{wilson2018greater, wilson2020}.

The idea of using extra coordinates decaying at a rate given by the Floquet exponents was introduced in the context of 2-dimensional biological oscillators in \cite{guillamon2009computational, castejon2013phase}. In these works, this variable, dynamically equivalent to the above mentioned isostable variable, was called amplitude variable. The addition of this amplitude variable -- that parameterizes the isochrons -- allowed to define the Phase Response Functions (PRFs), which extend the PRCs to the basin of attraction of the limit cycle \cite{guillamon2009computational}. Furthermore, complementing the isochrons, the authors also defined the level sets of the amplitude variable as A-curves in \cite{castejon2013phase} (equivalent to the isostables \cite{mauroy2013isostables}), corresponding to points on the basin of attraction which share the same asymptotic convergence. Analogously to the PRFs, they defined the Amplitude Response Functions (ARFs) (Isostable Response Functions in \cite{mauroy2018global}) quantifying the change in the amplitude variable due to a perturbation \cite{castejon2013phase}.

The approach in \cite{guillamon2009computational} relies on the parameterization method \cite{cabre2003parameterization, cabre2003parameterization2, cabre2005parameterization}. This method is specially well suited to study the dynamics around a hyperbolic attracting limit cycle in terms of the phase-amplitude variables. In particular, the parameterization method can be seen as a change of coordinates to the phase-amplitude variables providing a systematic way to compute the global (un)stable manifolds around the limit cycle. Besides being successfully applied in \cite{guillamon2009computational, huguet2013computation, castejon2013phase} to compute local and global isochrons, A-curves (corresponding to 2D isostables), iPRF and iARFs for planar systems, it has also been applied to compute invariant curves and PRCs in different problems of neuroscience \cite{PerezCervera2019, PerezCerveraHS17, castejon2020phase}.

In this paper we aim to extend the applications of the parameterization method in \cite{guillamon2009computational, castejon2013phase, huguet2013computation} to provide efficient algorithms to compute the full set of phase-amplitude coordinates, as well as local and global isochrons, isostables, iPRF and iARFs for $d$-dimensional oscillators, $d \geq 2$. Since our approach relies on the parameterization method, it benefits from all the previous solid theoretical and numerical framework in this area (see \cite{haro2016} for a review) providing -- thanks to the dynamical equivalence between the amplitude and the isostable variable -- a complete framework complementing most of the results of the existing approaches (we refer the reader to the Discussion section for more details). Although all the methods and algorithms presented herein are completely analogous in any dimension, for clarity of exposition we present them in dimension 3. The efficiency of the algorithms relies mainly on the use of Floquet theory to solve periodic linear high-dimensional systems \cite{castelli2015parameterization} and the use of automatic differentiation techniques \cite{griewank2008evaluating, haro2016} to compute the compositions of power series with elementary functions, thus, avoiding the computation of high-order derivatives of the vector field.  We apply these algorithms and compute an approximate parameterization in Fourier-Taylor series up to any degree, which provides analytically the local isochrons, isostables, iPRF and iARFs in a neighbourhood of the limit cycle. We stress that, in addition to these objects, the parameterization naturally provides the slow submanifold, which corresponds to the manifold associated to the smallest in modulus Floquet exponent \cite{cabre2003parameterization}. Furthermore, we use a numerical strategy based on backwards integration to globalize the slow manifold and then use it as a skeleton to globalize the rest of the objects in the full basin of attraction.  We illustrate the techniques by applying them to a selection of different 3-dimensional models for single neuron and neural populations in neuroscience. We stress that the objects and functions computed are obtained in the full basin of attraction of the limit cycle with a numerical accuracy that we previously established. Finally, we consider a periodically perturbed 3-dimensional single neuron model, and we use the computed iPRF and iARFs to explore different approximations of the stroboscopic map to describe the dynamics. Namely, we consider the full 3D phase-amplitude map, the 2D map that considers the reduction to the slow submanifold, thus considering the phase and the slow amplitude variable, and the classical 1D phase map. We emphasize that the parameterization $K$ and the isochrons provide a geometrical interpretation of the effects of the perturbation and the validity of the different dimensional reductions.

The structure of the paper is the following: in Section \ref{sec:section2} we provide the theoretical background to tackle the rest of the paper. In Section \ref{sec:section3} we introduce the theoretical and computational methodology to obtain the parameterization of the invariant manifolds of the limit cycle and we use it to globalize the isochrons, isostables, iPRF and iARFs. In Section \ref{sec:section4} we illustrate the methodology by computing these objects for different models in neuroscience. In Section \ref{sec:section5}, we apply the tools developed in the previous Sections to perform a study of the dynamics of a perturbed 3-dimensional single neuron model, while we explore the scope of validity of different dimensional reductions. We end the paper in Section \ref{sec:section6}, in which we present the conclusions of our work and its relation with other results in the field.

\section{Background and statement of the problem}\label{sec:section2}

%\underline{In this Section} we go through the background on the main tools that will be related later in sections \ref{sec:section3} and \ref{sec:section4}. Although these tools will be defined next for vector fields in $\mathbb{R}^d$, for the purposes of this paper we will restrict our attention to $d=3$ from Section \ref{sec:section3} on.

\subsection{Phase variable and Isochrons}\label{sec:phaseSec}

Consider an autonomous system of ODEs 
\begin{align}\label{eq:mathDef_1}
\dot{x} = X(x), \quad x \in \mathbb{R}^{d}, \quad d \geq 2 ,
\end{align}
whose flow is denoted by $\phi_t(x)$. Assume that $X$ is an analytic vector field and that system \eqref{eq:mathDef_1} has a $T$-periodic hyperbolic attracting limit cycle $\Gamma$, parameterized by the phase variable $\theta = t/T$ as
\begin{equation}
\begin{aligned}\label{eq:mathDef_2}
\gamma:\mathbb{T}:= \mathbb{R}/\mathbb{Z} &\to \mathbb{R}^{d}\\
\theta &\mapsto \gamma(\theta),
\end{aligned}
\end{equation}
so that it has period 1, that is, $\gamma(\theta) = \gamma(\theta +1)$ and $x(t) = \gamma(\frac{t}{T})$ satisfies \eqref{eq:mathDef_1}. Thus, the dynamics of \eqref{eq:mathDef_1} on $\Gamma$ can be reduced to a single equation
\begin{equation}
\dot{\theta} = \frac{1}{T}, \quad \quad \theta \in \mathbb{T}.
\end{equation}

\begin{remark}\label{rm:varRemakr}
We recall that a periodic orbit is hyperbolic attracting if all its characteristic exponents have negative real part except the trivial one which is 0 (or equivalently, the Floquet multipliers are inside the unit circle except the trivial one which is 1). The characteristic exponents of $\Gamma$ can be obtained by solving the variational equations of system \eqref{eq:mathDef_1} along the solution $\gamma(t/T)$. More precisely, the variational equations are given by the linearisation of the vector field $X$ around $\gamma$, that is,
\begin{equation}\label{eq:mjVarEqs}
\dot{\Phi} = DX(\gamma(t/T)) \Phi, \quad \quad \quad \text{with} \quad \Phi(0) = Id.
\end{equation}
The solution of the system above $\Phi(t)$ evaluated at $t = T$, i. e. $\Phi(T)$, is known as the monodromy matrix. The eigenvalues of $\Phi(T)$, namely $\mu_i = e^{\lambda_i T}$, $i = 0, ..., d-1$, are known as the Floquet multipliers of the limit cycle $\Gamma$ and the values $\lambda_i$ as the characteristic exponents. The index $i = 0$ will be assigned from now on to the trivial multiplier $\mu_0 = 1$, so $\lambda_0 = 0$.
\end{remark}

As we consider a hyperbolic attracting limit cycle $\Gamma$, the orbit of any point $p$ in the basin of attraction $\Omega$ of the limit cycle $\Gamma$ will approach asymptotically the orbit of a point $q$ in $\Gamma$ \cite{HirschPS77,guckenheimer1975}, that is,
\begin{align}\label{eq:synapseCurrent}
\lim_{t \to \infty} |\phi_t(q) - \phi_t(p)| = 0.
\end{align}
We will say that the two points $p$ and $q \in \Omega$ have the same asymptotic phase. We define the isochron $\mathcal{I}_{\theta}$ as the set of points having the same asymptotic phase $\theta$, that is,
\begin{equation}\label{eq:isochronsDef}
\mathcal{I}_{\theta} = \{x \in \Omega \mid \ |\phi_t(x) - \phi_t(\gamma(\theta))| = |\phi_t(x) - \gamma\left(\theta + \frac{t}{T}\right)| \rightarrow 0 \quad \text{as} \quad t \rightarrow \infty \}.
\end{equation}

Thus, the phase can be
extended under generic conditions to a neighbourhood of the limit cycle via the concepts of asymptotic phase
and isochrons \cite{winfree1967biological, guckenheimer1975}. When $\Gamma$ is a hyperbolic attracting periodic orbit, the isochrons correspond to the leaves of the stable manifold $\mathcal{M}$ of $\Gamma$, which coincides with its basin of attraction $\Omega$. The sets of points where the asymptotic phase is not defined are called \emph{phaseless sets} \cite{guckenheimer1975}.

%\begin{remark}\label{rm:phaselessSets}
%	Following former work \cite{guckenheimer1975}, we will define the phaseless sets as the sets of points where the asymptotic phase is not defined. Clearly, for an attracting Normally Hyperbolic Invariant Manifold the phaseless sets ar contained in $\mathbb{R}^n \backslash \mathcal{M}$  
%\end{remark}	

\subsection{Phase-Amplitude variables and the parameterization method} 

In this Section we present the so-called parameterization method adapted to our problem, which allows us to provide a description of the dynamics of the oscillator in terms of phase and amplitude variables. 
%which will be important in the results to follow

The parameterization method, introduced in \cite{cabre2003parameterization, cabre2003parameterization2, cabre2005parameterization}, is a general functional analytic scheme to study $n$-dimensional invariant manifolds of a dynamical system. The method consists in finding a parameterization of the invariant manifold by means of solving a functional equation. This equation characterizes the invariance of the manifold and expresses the dynamics on this manifold in the coordinates induced by the parameterization. The coordinates will be chosen in such a way that the dynamics expresses in the simplest way possible (in some cases linear).
%, which relates to the theory of normal forms. 

Mathematically, one looks for an embedding $K: \mathcal{U} \subset \mathbb{R}^n \rightarrow \mathbb{R}^d$ and a vector field $\mathcal{X}: \mathcal{U} \subset \mathbb{R}^n \rightarrow \mathbb{R}^n$, $n \leq d$ and $\mathcal{U}$ open, such that
\begin{equation}\label{eq:invariantgen}
DK \cdot \mathcal{X}=X \circ K, 
\end{equation}
so that \[\mathcal{M}:= Range(K)=\{ K(u) \in \mathbb{R}^d \quad |\quad u \in \mathcal{U} \subset \mathbb{R}^n \}\] 
is invariant under the flow of $X$ and the vector field $\dot{u}=\mathcal{X}(u)$ describes the dynamics 
on the invariant manifold $\mathcal{M}$.
Thus, the method provides information not only on the embedding but also on the dynamics on the manifold. 
Moreover, the method leads to efficient numerical algorithms to compute high order Taylor expansions of $K$ and $\mathcal{X}$. The book \cite{haro2016} contains a detailed description of the analytical and numerical aspects of the method in several contexts as well as a complete list of references.

\begin{remark} %\marginpar{pensar si es una conjugacio}
In our case, we look for a parameterization $K$ of the $d$-dimensional stable manifold $\mathcal{M}$ of the periodic orbit, which coincides with the basin of attraction $\Omega$ of $\Gamma$. This is equivalent to look for a change of variables $K$ that conjugates the vector field $X$ in $\Omega$ to a vector field $\mathcal{X}$ with a simpler expression of the dynamics.%, i.e. $\mathcal{X}_{|\Omega}=K^{-1} \circ X \circ K$.
\end{remark}

In our case, we look for a local analytic diffeomorphism
\begin{equation}\label{eq:kThetaSigma}
\begin{aligned}
K : \mathbb{T} \times \mathcal{B} \subset \mathbb{T} \times \mathbb{R}^{d-1} &\rightarrow \mathbb{R}^d  \\
(\theta, \sigma) &\rightarrow K(\theta, \sigma), 
\end{aligned}
\end{equation}
such that it satisfies the following invariance equation
\begin{equation}\label{eq:mjInvEq}
\frac{1}{T}\frac{\partial}{\partial\theta}K(\theta, \sigma) + \sum_{i=1}^{d-1} \lambda_i \sigma_i \frac{\partial}{\partial \sigma_i}K(\theta, \sigma) = X(K(\theta, \sigma)),
\end{equation}
where $T$ is the period of the limit cycle $\Gamma$ and $\lambda_1$, ..., $\lambda_{d-1} \in  \mathbb{R}$ its non-trivial characteristic exponents. Equation \eqref{eq:mjInvEq} will be the centrepiece of our approach.  

We can think of \eqref{eq:mjInvEq} as a change of coordinates. The new coordinates will be the phase  $\theta \in \mathbb{T}$ introduced in \eqref{eq:mathDef_2} and the amplitude coordinates $\sigma_1$, ..., $\sigma_{d-1} \in \mathbb{R}$, corresponding to transverse directions to the limit cycle. Thus, the dynamics of the vector field $X$ in \eqref{eq:mathDef_1} expressed in these new variables $(\theta, \sigma) \in \mathbb{T} \times (\mathcal{B} \subset \mathbb{R}^{d-1})$ is given by
\begin{equation}\label{eq:flipoTete}
\dot{\theta} = \frac{1}{T}, \quad \quad \quad  \dot{\sigma} = \Lambda \cdot \sigma, \quad \quad \text{with} \quad \quad \Lambda = 
\begin{pmatrix}
\lambda_1 & & \\
& \ddots & \\
& & \lambda_{d-1}  
\end{pmatrix}.
\end{equation} 
That is, the variable $\theta$ rotates at a constant speed $1/T$, while the variables $\sigma_i$ contract at a rate $\lambda_i$. Note that the vector field in \eqref{eq:flipoTete} is the vector field that we named $\mathcal{X}$ in \eqref{eq:invariantgen}. With this choice of $\mathcal{X}$, Eq.~\eqref{eq:invariantgen} in our general presentation becomes Eq.~\eqref{eq:mjInvEq}.

In Section \ref{sec:section3} we will show that, provided that the characteristic exponents satisfy certain non-resonance conditions, one can indeed solve the functional equation above and find a map $K(\theta, \sigma)$, at least formally. We foresee here that, for the purposes of this paper, we are going to assume that the characteristic exponents are real and distinct, that is, $\lambda_{1}<\cdots<\lambda_{d-1}<0$. Moreover, we are going to assume that they are non-resonant, that is
\begin{equation}\label{eq:resonance}
\sum_{i=1}^{d-1} m_i \lambda_i \neq \lambda_k, \quad \textrm{for any}\enskip k=1,\ldots,d-1,\enskip \textrm{for all} \enskip m_i \in \mathbb{Z} .
\end{equation}
The case of complex exponents or degeneracies, which is beyond the scope of this paper, is discussed in \cite{castelli2015parameterization}.

\begin{remark}
Notice that the characteristic exponents are negative because we assume that the limit cycle is hyperbolic attracting. Indeed, this is the interesting case in the neuroscience context. The case of a hyperbolic repelling limit cycle, though, is completely analogous, just reversing time. 
\end{remark}

\begin{remark}\label{rem:nres0}
The non-resonance condition \eqref{eq:resonance} will be necessary to solve the invariance equation \eqref{eq:mjInvEq} as we will see in Section~\ref{sec:section3} (see Remark~\ref{rem:nonres}). The existence of a resonance at some order does not prevent the existence of the 
stable invariant manifold and therefore the isochrons (see for instance, \cite{guckenheimer2013nonlinear}). However, it prevents the existence of an analytic conjugacy that conjugates the vector field $X$ to a linear one of the form \eqref{eq:flipoTete}. In that case, the parameterization method still works, but it is necessary to conjugate to a polynomial vector field instead of the linear field as done here. For a thorough discussion of resonances and the parameterization method for fixed points and equilibria, see \cite{cabre2003parameterization}. 
\end{remark}

\begin{remark}
Notice that Eq.~\eqref{eq:mjInvEq} does not have a unique solution. Indeed, if $K$ is a solution then $\tilde K(\theta,\sigma)=K(\theta + \omega, b \sigma)$ is also a solution, for any $\omega \in \mathbb{T}$ and $b \in \mathbb{R}^{d-1}$. The meaning of $\omega$ is the choice of the origin of time, and $b$ corresponds to the choice of units in $\sigma$. In Section~\ref{sec:section32} we discuss how to choose these constants to guarantee the numerical stability of the method.
\end{remark}
%

%The existence of the map $K(\theta, \sigma)$ --whose existence is ensured in the hyperbolic case by the stable manifold theorem (see \cite{hirsch1970stable})-- provides a parameterization of the stable manifold of $\Gamma$, and 
%therefore has important geometrical implications. 

The evolution of the flow $\phi_t$ in the coordinates $(\theta,\sigma)$ given in \eqref{eq:kThetaSigma} becomes
\begin{equation}\label{eq:aboveEq}
\phi_t(K(\theta, \sigma)) = K\Big(\theta + \frac{t}{T}, e^{\Lambda t}\sigma\Big).
\end{equation}

The map $K$ in \eqref{eq:kThetaSigma} allows to define a scalar function $\Theta$ that assigns the asymptotic phase to any  point $x$ in a neighbourhood $\Omega$ of the limit cycle $\Gamma$. Indeed, 
\begin{equation}
\begin{aligned}
\Theta:\Omega \subset \mathbb{R}^{d} &\to \mathbb{T},\\
x &\mapsto \Theta(x) = \theta \quad \quad \text{where} \enskip x = K(\theta, \sigma),\enskip \textrm{for some} \enskip \sigma \in \mathbb{R}^{d-1}.
\end{aligned}\label{eq:mathDef_3}
\end{equation}
Notice that $\Theta(\phi_t(x)) = \Theta(x) + \frac{t}{T}$. The level curves of $\Theta$ correspond to the isochrons $\mathcal{I}_{\theta}$ introduced in \eqref{eq:isochronsDef}, that is,
\begin{equation}
\mathcal{I}_{\theta} = \{x \in \Omega \quad | \quad \Theta(x) = \theta \}.
\end{equation}

Analogously, the map $K$ in \eqref{eq:kThetaSigma} also allows us to define the scalar functions $\Sigma_i$, for $i=1,\ldots, d-1$ that assign the amplitude variable $\sigma_i$ to any point $x \in$ $\Omega$: 
\begin{equation}\label{eq:sigma3Dcase}
\begin{aligned}
\Sigma_i : \Omega \subset \mathbb{R}^{d} &\to \mathbb{R},\\
x &\mapsto \Sigma_i(x) = \sigma_i, \quad  \text{where} \enskip x = K(\theta, \sigma), \enskip \textrm{for some} \enskip \theta \in \mathbb{T}.
\end{aligned}
\end{equation}
Notice that $\Sigma_i(\phi_t(x)) = \Sigma_i(x)e^{\lambda_i t}$, for $i=1,\ldots,d-1$. The level curves of $\Sigma_i$ are called  \textit{isostables}  (see \cite{mauroy2018global, mauroy2013isostables}) or \emph{A-curves} for the 2-dimensional case (see \cite{castejon2013phase}), and correspond to the sets of points
\begin{equation}\label{eq:aCurvesDef}
\mathcal{A}_{\sigma_i}^i = \{x \in \Omega \quad | \quad \Sigma_i(x) = \sigma_i \}.
\end{equation}
We will denote by $\Sigma$ the vector-valued function $\Sigma(x):=(\Sigma_1(x),\ldots,\Sigma_{d-1}(x))$.

\subsection{Phase and Amplitude Response Functions}\label{sec:phaseAmpFun}

Let us consider that an instantaneous pulse of amplitude $A$ is applied to the trajectory at time $t=t_s$ in the direction $\vec v \in \mathbb{R}^{d}$, that is
\begin{equation}
\dot{x} = X(x) + A \vec{v} \delta(t-t_s),
\end{equation}
where $\delta(t)$ is the Dirac delta function. We can assume without loss of generality that $t_s=0$. Such instantaneous perturbation, when acting over a trajectory on a point $x = K(\theta, \sigma) \in \Omega$, will displace the trajectory from the point $x$ to a new point $x + \Delta x=K(\theta_{new},\sigma_{new}) \in \Omega$ having a new phase $\theta _{new}$ and new amplitude $\sigma_{new}$ given by:
\begin{equation}\label{eq:defparf}
\begin{aligned}
\theta_{new} = \Theta(x + \Delta x) = \Theta(x) + PRF(A, \theta,\sigma), \\
\sigma_{new} = \Sigma(x + \Delta x) = \Sigma(x) + ARF(A, \theta, \sigma), 
\end{aligned}
\end{equation}
where the PRF and ARF are the Phase Response Function the Amplitude Response Function, respectively (see \cite{castejon2013phase, guillamon2009computational}). These functions quantify the shift in phase and amplitude due to the effect of the perturbation at a given point $x = K(\theta, \sigma) \in \Omega$, respectively.  

For the case of an instantaneous pulse of weak amplitude $|A| \ll 1$, then $\Delta x \ll 1$, and therefore
\begin{equation}
\begin{aligned}\label{eq:arcs_2}
\theta_{new} &= \Theta(x + \Delta x) = \Theta(x) + \nabla\Theta(x) \cdot \Delta x + O(|\Delta x|^2),\\
\sigma_{new} &= \Sigma(x + \Delta x) = \Sigma(x) + \nabla\Sigma(x) \cdot \Delta x + O(|\Delta x|^2), \quad \quad \quad \text{for} \quad x \in \Omega.
\end{aligned}
\end{equation}
Notice that the first order approximation of the PRF and the ARFs is given by $\nabla\Theta(x)$ and $\nabla\Sigma(x)$, which correspond to the infinitesimal PRF (iPRF) and the infinitesimal ARFs (iARFs), respectively.

\begin{remark}\label{rm:prcsRM}
	The PRF and the ARFs are the natural extension of the Phase Response Curve (PRC) and Amplitude Response Curves (ARCs). In particular, the PRC and ARCs correspond to the restriction of the PRF and the ARFs to the limit cycle ($\sigma = 0$), respectively. That is,
	\begin{equation}
	\begin{aligned}
	PRC(A, \theta) = PRF(A, \theta, 0),\\ 
	ARC(A, \theta) = ARF(A, \theta, 0). 
	\end{aligned}
	\end{equation}
	Moreover, the infinitesimal PRC (iPRC) and infinitesimal ARCs (iARC) correspond to
	\begin{equation}
	\begin{aligned}
	iPRC(\theta) = \nabla\Theta(\gamma(\theta)),\\ 
	iARC(\theta) = \nabla\Sigma(\gamma(\theta)). 
	\end{aligned}
	\end{equation}
\end{remark}	 

\begin{remark}\label{rem:sobreprf}
Notice that the PRF (resp. ARFs) are defined in \eqref{eq:defparf} as scalar-valued functions, while the iPRF $\nabla \Theta$ (resp. iARFs $\nabla \Sigma$) are defined in \eqref{eq:arcs_2} as vector-valued functions. Indeed, its range is a vector with $d$-components. However, in the neuroscience literature, since the perturbation $\Delta x$ occurs usually in the direction of the voltage (typically corresponding to the first component $x_1=V$), sometimes the iPRF (resp. iARFs) refer only to the real-valued function corresponding to $\partial \Theta (x)/\partial V$ (resp. $\partial \Sigma (x)/\partial V$). In this paper, abusing  language, we will refer to $\nabla \Theta$ (resp. $\nabla \Sigma$) and its first component indistinctly as iPRF (resp. iARFs).
\end{remark}

The iPRF $\nabla\Theta(x)$ and the iARFs $\nabla\Sigma(x)$ can be computed by means of the parameterization $K(\theta, \sigma)$ in \eqref{eq:kThetaSigma}. Indeed, taking derivatives on both sides at the expression $K(\Theta(x), \Sigma_1(x), ..., \Sigma_{d-1}(x)) = x$, we have for $x = (x_1, ..., x_d) \in \Omega$:
\begin{equation}\label{eq:prodidentity}
Id_{d \times d}
=
\begin{bmatrix} DK(\theta, \sigma) \end{bmatrix} \begin{bmatrix} \nabla\Theta(x) \\ \nabla\Sigma_1(x) \\
\vdots\\
\nabla\Sigma_{d-1}(x) \end{bmatrix}, 
\end{equation}
and therefore
\begin{equation}\label{eq:usefulSpace} \begin{bmatrix} \nabla\Theta(x) \\ \nabla\Sigma_1(x) \\
\vdots\\
\nabla\Sigma_{d-1}(x) \end{bmatrix} 
=
{\begin{bmatrix} DK(\theta, \sigma) \end{bmatrix}}^{-1}.
\end{equation}%\marginpar{invertim la K}
%
%In Section~\ref{sec:section3}, we will show how to obtain $K$ as a Taylor series in $\sigma$. As a consequence, we obtain expressions of $\nabla \Theta$ and $\nabla \Sigma_i$ at any order in $\sigma$ through Eq. \eqref{eq:usefulSpace}. Section~\ref{sec:section3} also illustrates how the numerical computation of $K(\theta, \sigma)$ is restricted to a domain $\Omega_{loc} \subseteq \Omega$. However, as for points $x = K(\theta, \sigma) \in \Omega$, $\nabla \Theta$ and $\nabla \Sigma_i$ satisfy the adjoint equations \cite{guillamon2009computational, castejon2013phase} 
%%
%\begin{equation}
%\begin{aligned}\label{eq:prcs_4b}
%\frac{d \nabla\Theta(\phi_t(x))}{dt} &= -DX^T(\phi_t(x))\nabla\Theta(\phi_t(x)),\\
%\frac{d\nabla\Sigma_i(\phi_t(x))}{dt} &= \left( \lambda_i - DX^T(\phi_t(x)) \right)\nabla\Sigma_i(\phi_t(x)),
%\end{aligned}
%\end{equation}
%%
%we can use values of $\nabla \Theta$ and $\nabla \Sigma_i$ in \eqref{eq:usefulSpace} as initial seeds for Eqs. in \eqref{eq:prcs_4b} and use these equations to extend $\nabla \Theta$ and $\nabla \Sigma_i$ beyond $\Omega_{loc}$ (see Section \ref{sec:isoComp3D} for more details). 

Moreover, as it is shown in \cite{guillamon2009computational, castejon2013phase}, for points $x = K(\theta, \sigma) \in \Omega$, $\nabla \Theta$ and $\nabla \Sigma_i$ satisfy the following adjoint equations
\begin{equation}
\begin{aligned}\label{eq:prcs_4b}
\frac{d \nabla\Theta(\phi_t(x))}{dt} &= -DX^T(\phi_t(x))\nabla\Theta(\phi_t(x)),\\
\frac{d\nabla\Sigma_i(\phi_t(x))}{dt} &= \left( \lambda_i - DX^T(\phi_t(x)) \right)\nabla\Sigma_i(\phi_t(x)),
\end{aligned}
\end{equation}
which is specially relevant for computational purposes. Indeed, as we discuss in Section \ref{sec:section32}, we use Eq. \eqref{eq:usefulSpace} to compute 
$\nabla \Theta$ and $\nabla \Sigma$ in a local neighbourhood of the limit cycle $\Omega_{loc} \subseteq \Omega$ using the local approximate expression for $K$. Moreover, this computation provides an initial condition for system~\eqref{eq:prcs_4b}, so that one can obtain the values of $\nabla \Theta$ and $\nabla \Sigma$ beyond $\Omega_{loc}$ by means of backwards integration.

%We foresee here that we a$\nabla \Theta$ is in first order in $\sigma$ a vector-valued function orthogonal to the space spanned by $\Phi(\theta T)v_i$, where $v_i$ are 
%the eigenvectors of the monodromy matrix $\Phi(\theta T)$ in \eqref{eq:mjVarEqs}. 
%Equivalently, $\nabla \Sigma_i$ is in first order orthogonal to the space spanned by the vector field $X(\gamma(\theta))$ Theta and $\Phi(\theta T)v_j$, where $v_j$ are 
%the eigenvectors of the monodromy matrix associated to $\lambda_j$ for $j=1, \dots, d-1$, $j \neq i$.
%\end{remark}

\begin{remark}
The classical adjoint method \cite{ErmentroutKopell91, brown2004phase, ErmentroutTerman2010} considers $x \in \Gamma$ in Eq.~\eqref{eq:prcs_4b} top, that is, with $\phi_t(x)=\gamma(t/T)$ and looks for a $T$-periodic solution to obtain the iPRC. However, if we want to extend the adjoint method to a neighbourhood of the limit cycle, we cannot impose periodicity conditions.  Thus, the problem lacks conditions to be solved uniquely, unless one knows a way to obtain initial conditions. We use that Eq.~\eqref{eq:usefulSpace} provides an expression for  $\nabla\Theta$ and $\nabla\Sigma$  for points not restricted to the limit cycle, and therefore it provides an initial condition to solve the adjoint equations \eqref{eq:prcs_4b} uniquely. So, the computation of $K(\theta, \sigma)$ allows us to extend the use of the adjoint equations beyond the limit cycle.
\end{remark}

\section{Methodology and Numerical Algorithm}\label{sec:section3}

In this Section, we describe the numerical methodology to solve the invariance equation \eqref{eq:mjInvEq} for the unknowns $K, T$ and $\Lambda$. In particular, in Section \ref{sec:section31}, we compute a formal expression for the parameterization $K(\theta, \sigma)$ in Fourier-Taylor series. Next, in Section \ref{sec:sectionNum}, we discuss the numerical implementation of the computation of $K(\theta, \sigma)$. As the methodology used to numerically obtain $K$ involves the truncation of the power series up to a given order $L$, in Section \ref{sec:section32} we discuss the domain of accuracy of the truncated $K$. Section~\ref{sec:local} is devoted to use the local approximation of $K$ to compute local isochrons, isostables, iPRF and iARFs. Finally, in Section \ref{sec:isoComp3D} we introduce a method to globalize the local approximation of the map $K(\theta, \sigma)$ to a larger domain which automatically provides the globalization of  the isochrons, isostables, iPRF and iARFs. Although the methodology can be applied to any system in $\mathbb{R}^d$, the details are given for the case $d=3$. %See \cite{castelli2015parameterization} for the general case.

%Without loss of generality we will assume that $\Gamma$ is stable. We are going to divide the method in several steps. IN order to perform each of them, we will use some standard methods/algorithms that we sketch now and that we explain in detail in this section. Tal i como dijimos lo haremos para el caso $d=3$

%\begin{algorithm}\label{alg:algorithm_new0}
%	\textbf{Computation of the PRC}.
%	Given $K_0(\theta)$ a parameterization of the limit cycle, and $g_A(\theta)$ an approximate solution of equation 
%	\eqref{eq:invariance_mod0} perform the following operations:
%	%
%	\begin{enumerate}
%		\item Compute $E(\theta) = F_A(K_0(\theta)) - K_0(g_A(\theta))$.
%		\item Compute $DK_0(g_A(\theta))$.
%		\item Compute $\Delta g_A = \frac{<DK_0 (g_A(\theta)),E(\theta)>}{<DK_0 (g_A(\theta)), DK_0(g_A(\theta))>}$.
%		\item Set $g_A(\theta) \leftarrow g_A(\theta)+\Delta g_A(\theta)$.
%		\item Repeat steps 1-4 until the error $E$ is smaller than the established tolerance. Then $PRC (\theta) = g_A(\theta) - \left(\theta + T'/T \right) $.
%	\end{enumerate}
%\end{algorithm}

\subsection{A formal solution for the invariance equation}\label{sec:section31}

The dynamics of system \eqref{eq:mathDef_1} for $d=3$, in terms of the phase-amplitude variables is given by
\begin{equation}\label{eq:thetaSigma3var}
\dot{\theta} = \frac{1}{T}, \quad \quad \quad  \dot{\sigma} = \Lambda \cdot \sigma, \quad \quad \text{with} \quad \quad \Lambda = 
\begin{pmatrix}
\lambda_1 & 0 \\
0 & \lambda_2  
\end{pmatrix},
\end{equation}
where $\sigma = (\sigma_1, \sigma_2) \in \mathbb{R}^2$, and  $\lambda_1$ and $\lambda_2$ are the characteristic exponents of $\Gamma$. We recall that we assume that $\lambda_1, \lambda_2$ are real and distinct (see Eq.~\eqref{eq:resonance}).

The invariance equation \eqref{eq:mjInvEq} in the case $d=3$ writes as 
\begin{equation}\label{eq:mjInvEq2}
\frac{1}{T}\frac{\partial}{\partial\theta}K(\theta, \sigma) + \sum_{i=1}^{2} \lambda_i \sigma_i \frac{\partial}{\partial \sigma_i}K(\theta, \sigma) = X(K(\theta, \sigma)).
\end{equation}

In order to solve the above invariance equation we assume a formal series solution for Eq.~\eqref{eq:mjInvEq2} of the form: 
\begin{equation}\label{eq:mjFourierTaylor}
K(\theta, \sigma) = \sum_{m=0}^{\infty} \sum_{\alpha=0}^{m} K_{\alpha, m-\alpha}(\theta) \sigma^{\alpha}_1 \sigma^{m-\alpha}_2,
\end{equation}
where the functions $K_{\alpha,m-\alpha}: \mathbb{T} \rightarrow \mathbb{R}^d$, for $\alpha = 0, ..., m$ and $m \in \mathbb{N}$.
Therefore, we substitute $K(\theta, \sigma)$ in \eqref{eq:mjFourierTaylor} in Eq.~\eqref{eq:mjInvEq2} and expand the vector field $X(K(\theta, \sigma))$ in Taylor series with respect to the variable $\sigma$ about $\sigma=0$. To obtain the expression for $K(\theta, \sigma)$ one just has to collect terms with the same power of $\sigma$ and solve the resulting equations. Next, we explain how to solve Eq.~\eqref{eq:mjInvEq2} for each degree.

For $m = 0$ the term $K_0(\theta) := K_{00}(\theta)$ satisfies the equation 
\begin{equation}\label{eq:eqMzero}
\frac{1}{T}\frac{d}{d\theta}K_0(\theta) = X(K_0(\theta)).
\end{equation}
Clearly, the solution of \eqref{eq:eqMzero} is the limit cycle itself, that is $K_0(\theta) = \gamma(\theta)$.
\begin{remark}\label{rm:remark31}
	Notice that if $K_0(\theta)$ is a solution, then $K_0(\theta + \omega)$ is also a solution for any $\omega \in [0,1)$. This means that the phase of a given oscillation can be fixed arbitrarily. As we will see in Section \ref{sec:section4}, for the examples in this paper, we follow the standard criterion in neuroscience which sets the zero phase at the maximum value of the voltage coordinate. 
\end{remark}

For $m = 1$, the equations for $K_{10}(\theta)$ and $K_{01}(\theta)$ are 
\begin{equation}\label{eq:eqMone}
\begin{aligned}
\frac{1}{T}\frac{d}{d\theta}K_{10}(\theta) + \lambda_1 K_{10}(\theta) &= DX(K_0(\theta))K_{10}(\theta), \\
\frac{1}{T}\frac{d}{d\theta}K_{01}(\theta) + \lambda_2 K_{01}(\theta) &= DX(K_0(\theta))K_{01}(\theta),
\end{aligned}
\end{equation}
respectively. The solutions for these equations are given by 
\begin{equation}\label{eq:eqMzeroSols}
K_{10}(\theta) = \Phi(\theta T)e^{-\lambda_1 \theta T }v_1, \quad \quad \quad K_{01}(\theta) = \Phi(\theta T )e^{-\lambda_2 \theta T}v_2,
\end{equation}
where $\Phi(t)$ is the solution of the variational equations \eqref{eq:mjVarEqs} and $v_i$ is the eigenvector of $\Phi(T)$ (the monodromy matrix) associated to the non-trivial $i$-th Floquet multiplier $\mu_i$, for $i = 1,2$.

\begin{remark}\label{rm:remark32}
Notice that the solutions $K_{10}(\theta)$ and $K_{01}(\theta)$ of \eqref{eq:eqMone} are non unique. If $v_1$ and $v_2$ are eigenvectors of the monodromy matrix so are $b_1 v_1$ and $b_2 v_2$ for any $b_1, b_2 \in \mathbb{R}$, giving rise to new solutions $b_1 K_{10}$ and $b_2 K_{01}$. Even though all the choices of $K_{10}$ and $K_{01}$ are mathematically equivalent, the choice affects the numerical properties of the algorithm. See Section \ref{sec:section32} for a more detailed discussion.
\end{remark}

Finally, for $m \geq 2$, the terms $K_{\alpha, m-\alpha}(\theta)$,  $\alpha = 0, ..., m$, satisfy the so-called \textit{homological equations}: 
\begin{equation}\label{eq:homologicalEqs}
\frac{1}{T}\frac{d K_{\alpha, m-\alpha}(\theta)}{d\theta} + (\alpha \lambda_1 + (m-\alpha) \lambda_{2}) K_{\alpha, m-\alpha}(\theta) = DX(K_0(\theta)) K_{\alpha, m-\alpha}(\theta) + B_{\alpha, m-\alpha}(\theta),
\end{equation}
where $B_{\alpha, m-\alpha}(\theta)$ is the coefficient of the term $\sigma_1^\alpha \sigma_2^{m-\alpha}$ in the Taylor expansion of
\begin{equation}\label{eq:Bexpand}
\emph{X}\left(  \sum_{n=0}^{m-1} \sum_{\alpha=0}^{n} K_{\alpha, n-\alpha}(\theta) \sigma^{\alpha}_1 \sigma^{n-\alpha}_2 \right).
\end{equation}
Notice that $B_{\alpha, m-\alpha}(\theta)$ is an explicit polynomial depending only on the terms of order lower than $m$, that is, the functions $K_{\alpha, n-\alpha}(\theta)$ for $n<m$ and whose coefficients are the derivatives of $X$ evaluated at $K_0$. They can be numerically computed using automatic differentiation techniques \cite{haro2016} (see Appendix \ref{sec:autoDifAp}).

The Eq.~\eqref{eq:homologicalEqs} can be solved assuming that $K_{\alpha, m-\alpha}(\theta)$ can be written in Fourier series so its coefficients are the unknowns. The resulting system of equations for the Fourier coefficients is linear, but it involves a large dimensional matrix which has a high computational cost. To avoid this numerical drawback, in the next Section, we review the method proposed in \cite{castelli2015parameterization} to solve the homological equation \eqref{eq:homologicalEqs} in an efficient way using the Floquet normal form.

\subsubsection{Reducibility of the homological equations via Floquet normal form}\label{sec:floqSchm}

In this Section we use the Floquet normal form to solve the homological equations \eqref{eq:homologicalEqs}. This allows us to transform the homological equations to a linear system with diagonal constant coefficient matrix in Fourier space, following \cite{castelli2015parameterization}. A similar idea has been also applied in \cite{huguet2013computation}. 

To avoid stodgy notation, from now on we will use $\su:=(\alpha,m-\alpha)$.

First of all recall that, by Floquet theory \cite{floquet1883equations}, the fundamental matrix $\Phi(t)$ of system \eqref{eq:mjVarEqs} can be written as
\begin{equation}\label{eq:floqNf}
\Phi(t) = \mathcal{Q}(t)e^{t R}, 
\end{equation}
where $\mathcal{Q}(t)$ is a $T$-periodic $3 \times 3$ matrix and $R$ is a real-valued $3 \times 3$ matrix. As $\Phi(t)$ is a solution of \eqref{eq:mjVarEqs}, using \eqref{eq:floqNf}, we have
\begin{equation}\label{eq:floquetDifInv}
\frac{1}{T} \left(\frac{d Q(\theta)}{d \theta} + Q(\theta) TR \right) = DX(\gamma(\theta)) Q(\theta),
\end{equation}
where $Q(\theta) := \mathcal{Q}(\theta T)$.

Then, we introduce the function $w: \mathbb{T} \rightarrow \mathbb{R}^3$ and write $K_{\su}(\theta)$ as 
\begin{equation}
K_{\su}(\theta) =  Q(\theta)w(\theta),
\end{equation} 
and substituting it in \eqref{eq:homologicalEqs}, we have
\begin{equation}\label{eq:mjIneed}
\frac{1}{T} \left(\frac{dQ(\theta)}{d\theta}w(\theta) + Q(\theta)\frac{dw(\theta)}{d\theta} \right) + \Upsilon Q(\theta)w(\theta) = DX(K_0(\theta))Q(\theta)w(\theta) + B_{\su}(\theta),
\end{equation}
where we have introduced the constant matrix $\Upsilon := (\alpha \lambda_1 + (m-\alpha) \lambda_{2})\cdot Id_{3 \times 3}$.

Then, using \eqref{eq:floquetDifInv} in equation \eqref{eq:mjIneed}, we obtain
\begin{equation}
\frac{1}{T} \Big(-Q(\theta)TRw(\theta) +  Q(\theta)\frac{d w(\theta)}{d \theta} \Big) + \Upsilon Q(\theta)w(\theta) = B_{\su}(\theta),
\end{equation}
and multiplying both sides by $Q^{-1}(\theta)$ we have
\begin{equation}
\label{eq:batiburillo}
\frac{1}{T}\frac{d w(\theta)}{d \theta} = (-\Upsilon + R)w(\theta) + Q^{-1}(\theta)B_{\su}(\theta).
\end{equation}

Finally, we assume that the matrix $R$ in \eqref{eq:floqNf} can be diagonalized, that is, there exists a matrix $C$ such that
\begin{equation}\label{eq:tereJ}
J = C^{-1}RC = \begin{pmatrix}
\lambda_0 & 0 & 0 \\
0 & \lambda_1 & 0  \\
0 & 0 & \lambda_2 
\end{pmatrix}.
\end{equation}
We make a final coordinate transformation and define $u(\theta)$ as $w(\theta) = Cu(\theta)$, and multiplying both sides by $C^{-1}$, expression \eqref{eq:batiburillo} reads as
\begin{equation}\label{eq:finalEqU}
\frac{1}{T} \frac{d u(\theta)}{d \theta} = (-\Upsilon + J)u(\theta) + A_{\su}(\theta),
\end{equation}
where 
\begin{equation}\label{eq:Afunction}
A_{\su}(\theta) = C^{-1}Q^{-1}(\theta)B_{\su}(\theta).
\end{equation}

Finally, we write $u(\theta)$ and $A_{\su}(\theta)$ in Fourier series, that is,
\begin{equation}\label{eq:mjIneedAgain}
u(\theta) = \sum_{k = -\infty}^{\infty} u_{k}e^{2\pi ik\theta}, \quad \quad A_{\su}(\theta) = \sum_{k = -\infty}^{\infty} A_{k}e^{2\pi ik\theta}, \quad \quad \quad A_{k}, u_{k} \in \mathbb{C}^3,
\end{equation}
and substitute expressions \eqref{eq:mjIneedAgain} in Eq.~\eqref{eq:finalEqU}. We obtain a linear system for the Fourier coefficients $u_k=(u_k^{(1)},u_k^{(2)},u_k^{(3)}) \in \mathbb{C}^3$ which is diagonal and can be solved componentwise, thus obtaining the following expression for the Fourier coefficients:
\begin{equation}\label{eq:mjFourierCoefs}
u^{(j)}_{k} = \frac{1}{\frac{2\pi ik}{T} + \alpha \lambda_1 + (m-\alpha) \lambda_{2} - \lambda_{j-1}}A_{k}^{(j)},
\end{equation}
for $j = 1, 2, 3$. Notice that the superindex $(j)$ refers to each component of the vectors $u_k$ and $A_k$. Finally, the solution $K_{\alpha, m-\alpha}(\theta)$ of Eq.~\eqref{eq:homologicalEqs} is given by 
\begin{equation}\label{eq:Kfunction}
K_{\alpha, m-\alpha}(\theta) = K_{\su}(\theta) = Q(\theta)Cu(\theta).
\end{equation}
\begin{remark}\label{rem:nonres}
	The Fourier coefficients $u^{(j)}_{k}$ in \eqref{eq:mjFourierCoefs} are formally well defined to all orders provided that, for any $k \in \mathbb{N}$ and $m\geq 2$, $\alpha=0,\ldots,m$, we have
	\[\frac{2\pi ik}{T} + \alpha \lambda_1 + (m-\alpha) \lambda_{2} - \lambda_{j-1} \neq 0, \, \quad j =1,2,3.\]
	Notice that this condition is always satisfied since we assumed that the characteristic exponents $\lambda_i$ are real negative and distinct, together with the the non-resonant condition \eqref{eq:resonance} (recall that $\lambda_0=0$).
\end{remark}

\begin{remark}\label{rem:general}
Notice that the Floquet reduction is computed only once and it is then used to find the solution of the homological equations at any degree. Indeed, from expressions \eqref{eq:Afunction}-\eqref{eq:mjIneedAgain}-\eqref{eq:Kfunction}, it is clear that to obtain the terms $K_{\su}$ for different $\su$'s, only the term $B_{\su}$ needs to be recomputed, while the matrices $C$ and $Q$ are always the same. 
\end{remark}

\subsection{Numerical computation of $K$}\label{sec:sectionNum}

In this Section we explain how to numerically solve Eqs.~\eqref{eq:eqMzero}, \eqref{eq:eqMone} and \eqref{eq:homologicalEqs} using the methodology described in Section \ref{sec:section31} and, thus obtain a local approximation of $K(\theta, \sigma)$.

\begin{itemize}
	
	\item For $m=0$ (see Eq.~\eqref{eq:eqMzero}), we need to compute the periodic solution $\Gamma$. To do so, we construct a Poincar\'e section and use a Newton method to find a fixed point of the corresponding Poincar\'e map. By doing this, we obtain a point $x_0 \in \Gamma$ and the period $T$. 
	
	We integrate system \eqref{eq:mathDef_1} with initial condition $x(0) = x_0$ and the variational equations \eqref{eq:mjVarEqs} altogether for a time $T$ to obtain $x(\theta T) =: K_0(\theta)$ and $\Phi(\theta T)$ for $\theta \in [0, 1)$. We store them for equidistant values of $\theta$; that is $\theta_i = i/N$ for $i = 0, ..., N-1$, which is equivalent to store the coefficients of the Fourier series up to order $N$. Indeed, we can switch between real and Fourier space by means of a Fast 
	Fourier Transform (FFT) algorithm \cite{brigham1978fast}.
	
	\item For $m=1$ (see Eq.~\eqref{eq:eqMone}), we consider the monodromy matrix $\Phi(T)$ and obtain its eigenvalues $\mu_i$ (Floquet multipliers), the Floquet exponents $\lambda_i = \frac{1}{T} \ln(\mu_i)$ and their respective eigenvectors $v_i$ for $i = 1,2$. We compute $K_{10}$ and $K_{01}$ according to the formulas \eqref{eq:eqMzeroSols}, and we store them again for the same equidistant values of $\theta$.
	
	\item For $m\geq2$, we use the scheme described in Section \ref{sec:floqSchm} to solve Eq.~\eqref{eq:homologicalEqs}. We solve these equations up to order $m=L$, $L \in \mathbb{N}$.
	
	To obtain the functions $B_{\su}$ in the homological equation~\eqref{eq:homologicalEqs} we need to compute the Taylor expansion of $X(f(\sigma_1, \sigma_2))$, where $f(\sigma_1, \sigma_2)$ is a particular Taylor polynomial in $\sigma_1$, $\sigma_2$ up to order $m-1$ (see Eq.~\eqref{eq:Bexpand}). As the vector field $X$ is analytic and consists of a combination of elementary functions, we can use automatic differentiation  techniques \cite{griewank2008evaluating, haro2016} to compute the coefficients $B_{\su}$ of the Taylor expansion of $X$ up to arbitrary order (see Appendix \ref{sec:autoDifAp} for more details).

	To apply the scheme described in Section \ref{sec:floqSchm}, we first compute the matrices $\mathcal{Q}$ (and therefore $Q$) and $R$ in \eqref{eq:floqNf} and the matrix $C$ satisfying \eqref{eq:tereJ} (see Appendix \ref{sec:floquetComputation}). We use the formula \eqref{eq:Afunction}  to obtain the functions $A_{\su}$ in real space and apply the FFT algorithm to obtain the Fourier coefficients $A_k$ (see~\eqref{eq:mjIneedAgain}). Finally, we use formulae \eqref{eq:mjFourierCoefs} and \eqref{eq:Kfunction} together with the inverse FFT to obtain the function $K_{\su}$ in real space. Again, we save this function for the same equidistant values of $\theta$.

%Concretament, a la pàgina 40. Capítol 2.3.1-2.3.2-2.3.3-2.3.4 

	To check the accuracy of the each of the solutions $K_{\su}$ obtained, we substitute them in their corresponding equation (Eq.~\eqref{eq:eqMzero} for $m=0$, Eq.~\eqref{eq:eqMone} for $m=1$ and Eq.~\eqref{eq:homologicalEqs} for $m\geq2$) for discrete values of $\theta$, that is, $\theta_i = i/N$ for $i = 0, ..., N-1$. For each value $\theta_i$, this substitution provides an error value $E_{\su}(\theta_i)$. For instance, for Eq.~\eqref{eq:eqMzero}, we obtain
	\[E_{0}(\theta_i)=\frac{1}{T} \frac{d}{d \theta} K_0(\theta_i)-X (K_0(\theta_i)),\]
	and, analogously, for Eq~\eqref{eq:eqMone} and \eqref{eq:homologicalEqs}.
	Finally, we compute the discrete $\ell_1$ norm of the error to get the accuracy, that is,
	\begin{equation}\label{eq:errortol}
	\Vert E_{\su} \Vert_{\ell_1} = \frac{1}{N} \sum_{i=0}^{N-1} |E_{\su}(\theta_i)|. 
	\end{equation}
	
%	\begin{remark}\label{rm:errorsK}
%		Although the solutions for equations X, X and X, from which we obtain X for m = X, X and X respectively are exact, they have some error when computing it numerically. As solutions for equations X and X are obtained through integration, the way of avoiding errors for X and X, is to use a numerical integrator as precise as possible. Then, to avoid errors

\end{itemize}

\subsection{Numerical errors and domain of accuracy of the approximate solution}\label{sec:section32}

Given a Taylor truncation at order $L \in \mathbb{N}$ and a Fourier truncation at order $N \in \mathbb{N}$, the above procedure provides an approximate solution of Eq.~\eqref{eq:mjInvEq2} of the form
\begin{equation}\label{eq:truncatedSeries3D}
\bar{K}(\theta, \sigma) = \sum_{m=0}^{L} \sum_{\alpha=0}^{m} \bar{K}_{\alpha, m-\alpha}(\theta) \sigma^{\alpha}_1 \sigma^{m-\alpha}_2, \enskip \text{where} \enskip  \bar{K}_{\su}(\theta) = \sum_{k=-N/2}^{N/2} c_{\suf} e^{2 \pi ik \theta}, \enskip c_{\suf} \in \mathbb{C}.
\end{equation}
Recall that we have defined $\su:=(\alpha,m-\alpha)$.

Next we discuss how to choose $L$ and $N$, as well as how to determine how good is the approximate solution and its domain of accuracy.

\subsubsection{Number of Fourier coefficients $N$}

%\begin{remark}\label{rm:taylorCoefs}
		To decide how many Fourier coefficients $N$ we have to compute (or equivalently, how many points on real space we have to store), we use the same criterion as in \cite{guillamon2009computational}. Thus, we pick a value $N$ such that the norm of the series with the last $10\%$ of Fourier coefficients is smaller than a given tolerance $E_{tail}$ (in the examples considered we set it at $10^{-10}$), that is
		\begin{equation}\label{eq:theTail}
		|K_{\su}^{tail}| = 2 \sum_{k=\lfloor 0.9N/2 \rfloor}^{N/2} |c_{\suf}| < E_{tail},
		\end{equation}
%	\end{remark}	
for all computed $\su$'s. We start with a value $N$ for which this condition is satisfied for $K_0$, $K_{10}$ and $K_{01}$, and we compute the rest of the $K_{\su}$, for $m=2 \ldots L$. At the end of the computation we check this condition for the new $K_{\su}$. Whenever it is not satisfied, we recompute all $K_{\su}$'s, for $|\su| = 0, \ldots, L$ again with $2N$ Fourier coefficients.

\subsubsection{Local approximation and number of Taylor coefficients $L$}

The solution $\bar{K}$ of \eqref{eq:mjInvEq2} is computed as a power series in $\sigma$, yet we do not expect that it is a good approximate
solution for all $\vert \sigma \vert > 0$. Rather we expect that it is good only in a neighbourhood of $\sigma=0$ (the limit cycle). Of course, the domain where 
the approximation is valid will depend on the error tolerance $E_{tol}$ and the order of the approximation $L$. 

Let us define the error function $E$ as
\begin{equation}\label{eq:errorEq}
E(\theta, \sigma) := \frac{1}{T}\frac{\partial}{\partial \theta} \bar{K}(\theta, \sigma) + \sum^2_{i = 1} \lambda_i \sigma_i \frac{\partial}{\partial \sigma_i}\bar{K}(\theta, \sigma) - X(\bar{K}(\theta, \sigma)).
\end{equation}

Therefore, for a given error tolerance $E_{tol} > 0$, we can define a numerical domain of approximation $\Omega_{loc}(E_{tol})$ for the solution $\bar K$ in the following way
\begin{equation}\label{eq:mjOmegaLoc}
\Omega_{loc}(E_{tol}):= \{x \in \Omega \subset \mathbb{R}^3 \mid x=K(\theta, \sigma) \textrm{ for } (\theta,\sigma) \in \mathbb{T} \times \mathbb{R}^2 \textrm{ and } \Vert E(\theta, \sigma) \Vert < E_{tol} \}, 
\end{equation} 
where $\parallel \cdot \parallel$ is the euclidean norm in $\mathbb{R}^3$. From now on we will refer to the domain $\Omega_{loc}$ without writing explicitly the dependence on $E_{tol}$.

In order to compute numerically the domain $\Omega_{loc}$, we perform the following strategy. For a fixed $\theta$, we write $\sigma_1=r \cos \varphi$ and $\sigma_2=r \sin \varphi$, for $r \in \mathbb{R}^+$ and $\varphi \in [0,2 \pi)$. For a fixed value $\varphi$, we look for the maximum $r$ such that
\[x=K(\theta,r \cos \varphi, r \sin \varphi) \in \Omega_{loc}.\]
Thus, we define the functions $R_{\theta}: [0,2\pi) \rightarrow \mathbb{R}^+$, such that
\begin{equation}\label{eq:smax}
x=K(\theta,\sigma_1,\sigma_2) \in \Omega_{loc} \iff \Vert (\sigma_1,\sigma_2) \Vert < R_{\theta} (\varphi), \, \textrm{ where } \tan \varphi=\frac{\sigma_2}{\sigma_1}.  
\end{equation}
Moreover, we define the set $\mathcal{B}_{loc}(\theta)$ as
\begin{equation}\label{eq:bloc}
\mathcal{B}_{loc}(\theta):=\{ (\sigma_1,\sigma_2) \in \mathbb{R}^2 \quad | \quad \Vert (\sigma_1,\sigma_2) \Vert < R_{\theta} (\varphi), \, \textrm{ where } \tan \varphi= \frac{\sigma_2}{\sigma_1}\}.
\end{equation}

As we already mentioned in Remark \ref{rm:remark32}, the choice of $v_1$ and $v_2$ in \eqref{eq:eqMzeroSols} is non unique. Although, theoretically, we can choose any vectors $b_1 v_1$, $b_2 v_2$, for any $b_1, b_2 \in \mathbb{R}$, their choice affects the numerical stability of the method and, more importantly, the size of the local approximation. Indeed, once $b_{1}$ and $b_{2}$ are fixed, the monomial $K_{\alpha, m-\alpha}$ will be multiplied by a factor $b^\alpha_{1}b^{m-\alpha}_{2}$. Therefore, if one chooses small values for $b_{1}$ and $b_{2}$, the terms $K_{\alpha, m-\alpha}$ will become small very fast as $m$ increases, so
increasing the order $L$ does not provide any extra information since we are just adding terms that are smaller than the machine error. By contrast, if one chooses large values for $b_{1}$ and $b_{2}$,  the functions $K_{\alpha, m-\alpha}$ will blow up fast as $m$ increases and,
eventually, those values will be too large to computationally operate with them. Therefore, to obtain a local approximation that extends to a larger neighbourhood, one has to choose appropriate values for $b_{1}$ and $b_{2}$ so that the coefficients $K_{\alpha, m-\alpha}$ can be kept at order 1 and the round-off errors are greatly reduced. We refer the reader to \cite{huguet2013computation} for a more detailed discussion on the role of $b_1$ and $b_2$. 

We determine the value of $b_1$ and $b_2$ and the order of the expansion $L$ by numerical experimentation. We typically have in mind a certain error tolerance for \eqref{eq:errortol} (around $10^{-6}$) and we stop at a certain order $L$ whenever this error is larger than the tolerance and cannot be made smaller by changing $b_1$ and $b_2$. So, we look for solutions computed to the highest order and having the largest accurate domains.

\subsection{Local isochrons, isostables, iPRF and iARFs}\label{sec:local}

In the previous Section, we have discussed the domain of accuracy of $\Omega_{loc}$ for the local approximation $\bar K$. Thus, in this domain, we define the local isochrons $\iloc$ and the local isostables ($A$-surfaces) $\aloc$, $i=1,2$, as
\begin{eqnarray}\label{eq:localIsos}
\iloc&:= \{ x \in \Omega_{loc} \quad \mid \quad \Theta(x) = \theta\}, \quad \theta \in \mathbb{T}\\
\label{eq:localIsostables}
\aloc&:= \{ x \in \Omega_{loc} \quad \mid \quad \Sigma_i(x) = c\}, \quad c \in \mathbb{R}.
\end{eqnarray}
Thus, $\iloc$ can be computed by evaluating the function 
\[
\begin{array}{rccc}
K(\theta,\cdot): &  \mathbb{R}^2 & \rightarrow & \mathbb{R}^3 \\
& \sigma & \rightarrow &  K(\theta,\sigma),
\end{array}
\]
for points $\sigma \in \mathcal{B}_{loc}(\theta)$ defined in \eqref{eq:bloc}. In Section~\ref{sec:isoComp3D} we will discuss how to choose a grid on $\mathcal{B}_{loc}(\theta)$, where we will evaluate the function above.

Moreover, the local iPRF $\nabla \Theta^{loc}$ and the local iARFs $\nabla \Sigma_i^{loc}\enskip i=1,2$, can be computed straightforwardly using formula \eqref{eq:usefulSpace}, for points $x = K(\theta, \sigma) \in \Omega_{loc}$. In Appendix \ref{sec:expPotenciesGradient} we show how to compute the power expansions in $\sigma$ of these functions up to any order.

\subsection{Globalization of $K$, isochrons, isostables, iPRF and iARFs}\label{sec:isoComp3D}

In this Section we explain how to extend to a larger domain $\Omega_c \subset \Omega$, the isochrons $\mathcal{I}_\theta$, the isostables $\mathcal{A}^i_c$, the iPRF $\nabla \Theta$ and the iARFs $\nabla \Sigma$, which are known in 
a local domain $\Omega_{loc}$ (Section~\ref{sec:local}). To do so, we will use that $K(\theta,\sigma)$ is invariant by the flow $\phi_t$ of the vector field $X$ (see Eq.~\eqref{eq:aboveEq}).
%property of the parameterization $K$ such that the flow $\phi_t$ of the vector field $X$ satisfies
%\[ \phi_t(K(\theta,\sigma))=K(\theta+t/T,\sigma e^{\Lambda}t).\]
We refer to this procedure as the globalization process.

%\begin{remark}
%We refer to the computational domain as the domain $\Omega_c \subset \Omega$ where we aim to compute the global
%isochron. For more details and examples see Section~\ref{sec:section4}.
%\end{remark}

Let us start by explaining how to extend the isochron $\iloc$. 
%Since the isochrons satisfy that
%\begin{equation}\label{eq:isoback}
%\phi_{-nT}(\mathcal{I}_{\theta}) \subset \mathcal{I}_{\theta}, \quad \textrm{for} n \in \mathbb{Z},
%\end{equation}
%and the backwards dynamics is expanding, integrating points on the isochron backwards for a multiple of the period we can extend the isochron to a larger domain than %$\Omega_{loc}$. 

Since the flow of the vector field $X$ takes isochrons to isochrons, we can obtain several points on the isochron of phase $\theta$, $\iiso$, by integrating backwards for a time $\Delta t$ points on the isochron of phase $\theta + \Delta t/T$, i.e, $\mathcal{I}_{\theta+ \Delta t/T}$, which is known at least locally, that is,
\begin{equation}\label{eq:isoback}
\phi_{- \Delta t}(\mathcal{I}_{\theta+ \Delta t/T}) \subset \mathcal{I}_{\theta}.
\end{equation}
Since we are integrating backwards and the backward dynamics is expanding, this procedure extends the isochron to a larger domain than $\Omega_{loc}$. In particular, we can use $\Delta t= n T, \enskip n \in \mathbb{N}$, that is, we can use points on the the same isochron to globalize it. We will present the method for the latter case, but we stress that it can be adapted to consider also points on other isochrons (see \cite{guillamon2009computational}).

The globalization of isochrons by means of backwards integration using the property \eqref{eq:isoback} in larger dimensional spaces ($d > 2$) presents a fundamental challenge: it involves different expanding directions, each of them associated to a particular Floquet multiplier of $\Gamma$. Therefore, by taking an homogeneous distribution of points on the local isochron $\iloc$, most of the points will escape in the direction of the largest in modulus Floquet exponent. To avoid this drawback, 
we have taken different actions:
\begin{itemize}
 \item The use of higher order expansions for $\bar K$ (large $L$ in Eq.~\eqref{eq:truncatedSeries3D}), which allows us to start the computations at a relatively large distance from the limit cycle. This fact decreases the integration time to expand the isochron. Thus, our methods represent an advantage in front of other methods that consider only first order expansions. %\cite{OsingaMoehlis} 
 \item To choose the points in $\Omega_{loc}$ in a clever way so that one obtains a homogeneous distribution of points on the isochron $\mathcal{I}_{\theta}$.
Next, we provide a detailed description of the scheme we follow, adapted from \cite{simo1990analytical}.
\end{itemize}

Assume without loss of generality that $|\lambda_1| > |\lambda_2|$. Consider the $2$-dimensional isostable
\begin{equation}\label{eq:slowManifold}
\mathcal{S} := \mathcal{A}^1_0 =\{x \in \Omega \quad | \quad \Sigma_1(x) = 0  \}.
\end{equation}
Notice that $\mathcal{S}$ corresponds to the \emph{slow} stable manifold. Here, following \cite{cabre2003parameterization, cabre2005parameterization}, the concept of slow manifold refers to the invariant submanifold of the full stable manifold associated to the slowest direction (corresponding to the smallest in modulus Floquet exponent $\lambda_2$). These manifolds are important because while the faster directions get rapidly suppressed, the slowest directions dominate the asymptotics of convergence. In many applications one can use the slow manifold to study possible reductions of the dynamics. This will be done for a periodically forced system in Section~\ref{sec:section5}.
%is interested in phenomena that occur at a slow time scale. Thus, the slow manifolds are the geometric objects that organize the dynamics. As an exemple of so, we will

%\begin{remark}\label{rem:slow}\marginpar{\textcolor{red}{potser treuria aquesta remarca}}
%In the literature the concept of slow manifold appears in several contexts. For instance, the centre manifold associated to the zero eigenvalues is sometimes referred as the slow manifold \cite{}. 
%In the study of slow-fast systems, the slow manifold refers to the geometric object where the velocity of the fast variable vanishes \cite{}. 
%\marginpar{es el mateix aqui?}
%Here, we are dealing with non-resonant hyperbolic limit cycles, and slow manifold refers to the direction associated to the smallest multiplier, and hence, controls the asymptotic behaviour. See \cite{cabre2003parameterization, cabre2005parameterization} for more details. 
%\end{remark}

We want to emphasize that the parameterization method provides straightforwardly the  slow manifold. Indeed,
\[\mathcal{S}= \{ x \in \Omega \quad | \quad x=K(\theta, 0, \sigma_2), \quad \theta \in \mathbb{T}, \quad \sigma_2 \in \mathbb{R}\}.\]
%This illustrates the fact that the parameterization method allows computation of invariant submanifolds of the
%full stable manifold associated with slow eigendirections. (Here we use the term “slow stable
%manifold” in the sense discussed in [1, 3]. Briefly, we are considering the manifold defined by
%conjugating the phase space dynamics near the periodic orbit to the linear dynamics associated
%with the slowest stable eigenvalue, in this case the eigenvalue λ 9 = −0.5730.)

%the component along v is the one that decays more slowly and hence, the one which controls the asymptotic behavior. Theinvariant manifold associated to this eigenvalue is a nonlinear analogue and canalso be used to study the asymptotic behavior of the iterates ofF. It is usuallycalled a slow manifold

%In many applications, one is interestedin the phenomena that happen at a slow time-scale. The fast phenomena dis-appear quickly and are not observable. Very often, these slow phenomena aregoverned by geometric objects that are described asslow manifolds. Two sciencesin which such phenomena happen and have been considered are atmospheric sci-ences and chemical kinetics.

It will be convenient to consider the foliation of the slow manifold $\mathcal{S}$ given by $\{S^{\theta}\}_{\theta \in \mathbb{T}}$, where:
\begin{equation}\label{eq:sltheta}
\mathcal{S}^\theta = \{x \in \Omega \quad | \quad \Theta(x) = \theta, \quad \Sigma_1(x) = 0  \}.
\end{equation}
Notice that $\mathcal{S}^{\theta} \subset \iiso$. Thus, our method to compute isochrons in $\Omega_c$ starts by computing points on the leaf of the slow manifold $\mathcal{S}^{\theta}$ in the domain $\Omega_c$, and then use it to obtain the isochron $\iiso$ (see Fig.~\ref{fig:globalMethod3d}). %\marginpar{revisar la Fig1 amb la notacio del text}

We present first an algorithm to extend the foliation of the slow manifold $\mathcal{S}$ provided by the parameterization $\bar{K}$ locally in the domain $\Omega^{loc}$ to a domain $\Omega_c$. Thus, for each $\mathcal{S}^{\theta}$, we are going to generate a sequence of points ${x_0,...,x_m}$ on it, such that they are at a distance (Euclidean norm in $\mathbb{R}^3$) smaller than some tolerance $\Delta_{max}$. The strategy is depicted in Fig.~\ref{fig:globalMethod3d}A.

\begin{algorithm}\label{alg:slowMani}	
	\textbf{Computation of the slow manifold leaf $\mathcal{S}^\theta$}. Given an approximate solution $\bar K$ (see \eqref{eq:truncatedSeries3D}) in a local domain $\Omega_{loc}$, a constant $\Delta_{max}$ that determines the maximum distance between points on the manifold $\mathcal{S}^\theta$ and a constant $\sigma_{max}:=R_{\theta}(\pi/2)$, where $R_{\theta}$ is defined in \eqref{eq:smax}, perform the following operations: 
	\begin{enumerate}
		\item Compute the point $x_0 = K_0(\theta)$. Start the list $\{x_0\}$.   
		\item Set $k = 0$, $n=0$, $\sigma=0$ and $\Delta \sigma=0.8 \sigma_{max}$. Notice that here $\sigma, \Delta \sigma \in \mathbb{R}$. 
		\item Compute $x_{int} = K(\theta, 0, \sigma+\Delta \sigma)$ and $x_{k+1} = \phi_{-nT}(x_{int})$.
		\begin{itemize}
			\item If $ \Vert x_{k+1} - x_{k} \Vert < \Delta_{max}$ then add $x_{k+1}$ to the list $\{x_0, \ldots, x_{k}\}$ and set $\sigma \leftarrow \sigma + \Delta \sigma$ and $k \leftarrow k+1$.
			\item Else, divide $\Delta \sigma $ by 2. 
		\end{itemize}
		 Repeat this step until $\sigma + \Delta \sigma > \sigma_{max}$. 
		\item Then set $\sigma \leftarrow \sigma_{min}:=\sigma_{max} e^{\lambda_2 T}$, $\Delta \sigma=0.8 (\sigma_{max}-\sigma_{min})$, and $n \leftarrow n + 1$, and repeat from step 3 until $x_{k+1}$ is out of the computational domain $\Omega_c$.
	\end{enumerate}
\end{algorithm}

\begin{remark}
Notice that for $n=0$, $x_{k+1}=\phi_0(x_{int})=x_{int}$ in step 3.
\end{remark}

\begin{remark}
To extend the leaf $\mathcal{S}^\theta$ for negative values of $\sigma$ we use Algorithm~\ref{alg:slowMani} with 
$\sigma_{max}=-R_{\theta}(3\pi/2)$. 
\end{remark}

Next, we present an algorithm that extends the local isochron $\iloc$ to $\Omega_c$ from the previously computed manifold $\mathcal{S}^{\theta}$. The strategy is the same as in the previous algorithm and is depicted in Fig.~\ref{fig:globalMethod3d}B.

\begin{algorithm}\label{alg:extendIso}
	\textbf{Computation of the isochron $\iiso$}. Given $\{x_0,..,x_m\}$ points on $\mathcal{S}^\theta$ and pairs\\ $\{(\sigma^ {(0)},N_0),\ldots (\sigma^ {(m)},N_m)\}$, where $(\sigma^{(k)},N_k) \in \mathbb{R} \times (\mathbb{N} \cup \{0\})$ such that $x_k=\phi_{-N_k T}(x_k^{(0)})$, where $x_k^{(0)}=K(\theta,0,\sigma^{(k)})$, $k=0,\ldots,m$, provided by Algorithm \ref{alg:slowMani}. 
	
	Consider also given an approximate solution $\bar K$ in a local domain $\Omega_{loc}$, a constant $\Delta_{max}$ that determines the maximum distance between points on the manifold $\mathcal{I}_{\theta}$ and a constant $\sigma_{max}$ (which is computed using the functions $R_{\theta}$ defined in \eqref{eq:smax}), perform the following operations:
	\begin{enumerate}
		
	\item Set $k=1$ and $q=0$. 
	
	\item Set $x_k^{(0)}=x_k$. Start the list $\{x^{(0)}_k\}$. 
	
	\item Set $\sigma=0$, $\Delta t=N_k T$ and $\Delta \sigma=0.8 \sigma_{max}$. Notice that here $\sigma, \Delta \sigma \in \mathbb{R}$. 
	%R_(\theta}(\varphi)$, where $\varphi=\arctan(\sigma^{(k)}/)$. 
	
		\item Compute $x_{int} = K(\theta, \sigma+\Delta \sigma,\sigma^ {(k)})$ and $x_{k}^{(q+1)} = \phi_{- \Delta t}(x_{int})$.
		\begin{itemize}
			\item If $\Vert x^{(q+1)}_k - x^{(q)}_{k} \Vert < \Delta_{max}$ then add $x_{k}^{(q+1)}$ to the list $\{x_k^{(0)}, \ldots, x_{k}^{(q)}\}$ and set $\sigma \leftarrow \sigma + \Delta \sigma$ and $q \leftarrow q+1$.
			\item Else, divide $\Delta \sigma $ by 2. 
		\end{itemize}
	 Repeat this step until $\sigma + \Delta \sigma > \sigma_{max}$.

	 %\item Pick $ \Delta \theta $ such that $\phi_{- \Delta \theta T} (K(\theta, \sigma, \sigma^{(k)})) \in \Omega_{loc}^{\theta + \Delta \theta}$. Then set $\sigma \leftarrow \sigma e^{\lambda_1 \Delta \theta T}$, $\sigma^{(k)} \leftarrow \sigma^{(k)} e^{\lambda_2 \Delta \theta T}$, $\Delta t \leftarrow (N_k + \Delta \theta) T$, $\Delta \sigma \leftarrow 0.8 S_1(\theta+ \Delta \theta)$ and $\theta \leftarrow \theta + \Delta \theta$. Repeat from step 4 until $x_k^{(q+1)}$ is out of the computational domain $\Omega_c$.
	
	\item Then set $\sigma^{(k)} \leftarrow \sigma^{(k)} e^{\lambda_2 T}$, recompute $\sigma_{max}$ for the new $\sigma^{(k)}$ and set $N_k \leftarrow N_k + 1$, and repeat from step 3 (with $\sigma \leftarrow \sigma_{max} e^{\lambda_1 T}$) until $x_{k+1}$ is out of the computational domain $\Omega_c$. 
	
	\item Set $k \leftarrow k+1$, $q \leftarrow 0$ and repeat from step 2.
	
	\end{enumerate}
\end{algorithm}	

\begin{remark}
	Notice that to completely extend the isochron around $\Sigma_2 = c$ one has to repeat algorithm \ref{alg:extendIso} with $\Delta \sigma \leftarrow - \Delta \sigma$.
\end{remark}

%\begin{remark}
%Notice that a variation of this algorithm and the previous one would be to use also computed points on other isochrons using the property \eqref{eq:isoback} (see for instance \eqref{guillamon2009computational}).
%\end{remark}

\begin{remark}
If the contraction in the direction of $\sigma_1$ is strong enough so that the points  $x_k^{(q)}$ escape fast from the isostable $\Sigma_1=0$, we have that they fall outside of the computational domain $\Omega_c$ without reaching the value $\sigma_{max}$. Thus step 5 in the previous algorithm does not need to be applied (See Fig.~\ref{fig:globalMethod3d}B). We will see that this is the case in all the examples considered in Section~\ref{sec:section4},
\end{remark}

Notice that the way we use to globalize the isochrons is to compute first the isostable $\Sigma_1=0$ and later the isostables $\Sigma_2=c$, for $c \in \mathbb{R}$ (see Fig.~\ref{fig:globalMethod3d}C).

%\begin{remark}
%The map $\phi_{-\Delta t}$ is computed by integrating the vector field $X$ backwards for a time $\Delta t$.
%\end{remark}

%\begin{remark}\label{rem:zeroint}
%Notice that $N_k$ can be zero. In this case, $x^{(q+1)}_k=\phi_0(x_{int})=x_{int}$ in step 4.
%\end{remark}

\subsubsection{Computation of the iPRF and iARFs}

The strategy used in Algorithm~\ref{alg:extendIso} also permits the computation of the iPRF $\nabla\Theta(x)$ and the iARFs $\nabla\Sigma(x)$ using expression \eqref{eq:usefulSpace} and equations \eqref{eq:prcs_4b}. Notice first that each point $x$ on the isochron $\iiso$ has been obtained either by direct evaluation of the parameterization $\bar{K}$, that is $x=\bar{K}(\theta, \sigma)$ for some $(\theta, \sigma) \in \mathbb{T} \times \mathbb{R}^2$ or by integrating backwards for a time $\Delta t$ the flow $\phi_t$ of the vector field $X$ starting at a point $x_{int} \in \Omega_{loc}$, which at its turns is obtained from the parameterization $\bar{K}$, that is, $x_{int} = \bar{K}(\theta, \sigma)$ for some $(\theta, \sigma) \in \mathbb{T} \times \mathbb{R}^2$ (see Algorithm \ref{alg:extendIso}). Thus, for each computed point $x \in \iiso$, we can obtain 
$\nabla \Theta (x)$ and $\nabla \Sigma(x)$ in the following way:

\begin{itemize}
\item If $x \in \Omega_{loc}$ then $x=K(\theta^*, \sigma^*)$ for some $(\theta^*, \sigma^*) \in \mathbb{T} \times \mathbb{R}^2$. Using the local approximation $\bar K$ and the values $(\theta^*,\sigma^*)$, $\nabla \Theta (x)$ and $\nabla \Sigma(x)$ are obtained from expression \eqref{eq:usefulSpace}.

\item If $x \notin \Omega_{loc}$ then $x=\phi_{-\Delta t} (x_{int})$, where $x_{int}=K(\theta^*, \sigma^*)$ for some $(\theta^*, \sigma^*) \in \mathbb{T} \times \mathbb{R}^2$
and $\Delta t$. Using the local approximation $\bar K$ and the values $(\theta^*,\sigma^*)$, we compute $\nabla \Theta (x_{int})$ and $\nabla \Sigma(x_{int})$ using expression \eqref{eq:usefulSpace} and then
\[ \nabla \Theta (x) = \Psi^{\theta}_{-\Delta t} (\nabla \Theta (x_{int})) \quad \textrm{ and } \quad \nabla \Sigma (x) = \Psi^{\sigma}_{-\Delta t} (\nabla \Sigma (x_{int})),\]
where $\Psi^{\theta}_t, \Psi^{\sigma}_t$ is the flow of \eqref{eq:prcs_4b}.
\end{itemize}

Thus, expression \eqref{eq:usefulSpace} provides an accurate approximation of these functions in the domain $\Omega_{loc}$, and  Eqs.~\eqref{eq:prcs_4b} allow to globalize them in $\Omega_c$.

\subsubsection{Computation of isostables $\mathcal{A}^i_{c}$}

Algorithm \ref{alg:slowMani} provides a relevant isostable (the slow manifold $\mathcal{S}:= \mathcal{A}^1_0$), while Algorithm \ref{alg:extendIso} provides the isostables $\mathcal{A}^2_c$, i.e, $\Sigma_2=c$, for $c\in\mathbb{R}$. To compute any other isostable $\mathcal{A}^i_{c}$ for $i=1,2$, defined in \eqref{eq:aCurvesDef}, one can adapt straightforwardly the strategy of the mentioned algorithm. Indeed, the isostable $A^i_{c}$ can be obtained from the locally computed isostable $A^i_{c^*}$ (for $|c|>|c^*|> 0$), just integrating the vector field $X$ backwards for a time
\begin{equation}
t = \frac{1}{\lambda_i} \ln \left(\frac{c}{c^*}\right).
\end{equation}

\begin{figure}[H]
\begin{minipage}[c]{0.65\textwidth}
{\includegraphics[width=98mm]{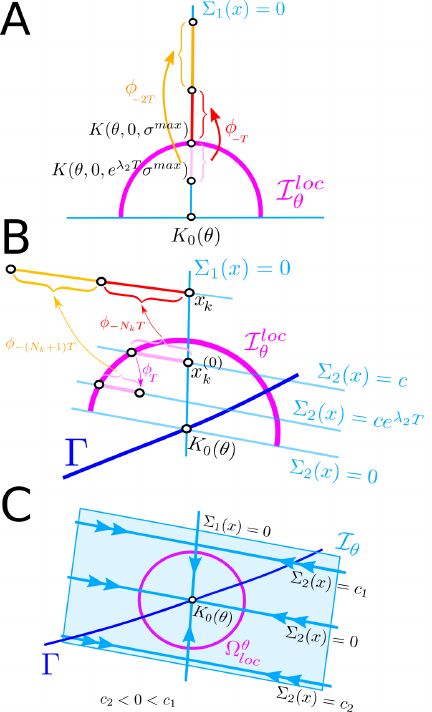}}
\end{minipage}\hfill
\begin{minipage}[c]{0.3\textwidth}
\caption{Schematic representation for the globalization procedure of the isochrons. The purple curve represents the border of the local isochron $\mathcal{I}^{loc}_{\theta}$, the blue curves correspond to the isostables $\Sigma_1=ct$ and $\Sigma_2=ct$, and the dark blue curve is a piece of the limit cycle $\Gamma$.
	(A) Strategy presented in Algorithm~\ref{alg:slowMani} to globalize the leaf $\mathcal{S}^\theta$ ($\Sigma_1(x)=0$) of the slow manifold $\mathcal{S}$. (B) Strategy presented in Algorithm \ref{alg:extendIso} to use the points $x_k$ on $\mathcal{S}^\theta$ to extend the isochron $\mathcal{I}_\theta$ along the curve of constant $\sigma_2$, $\Sigma_2(x)=\Sigma_2(x_k)$ (C) By repeating the strategy in (B) for all the points $x_k$ in the leaf $\mathcal{S}^\theta$ of the slow manifold, one obtains the  global isochron.  See text for more details.} \label{fig:globalMethod3d}
\end{minipage}
\end{figure}

%\textcolor{blue}{Q: Notation does not coincide with the notation in the text and it is confusing. Not enough accurate. Arrows indicating slow and fast directions would help a lot. My guess is that the clue is to find a foliation  of the slow manifold?}

\begin{remark}
Globalization methods require precise integration routines and good numerical precision of the local approximation $\bar{K}$. We have used a 8th-order Runge-Kutta Fehlberg method (rk78) with a tolerance $10^{-14}$. For the models we have considered, the use of stiff and non stiff-methods yields similar results. The numerical precision of the local approximation $\bar{K}$ was found to be even more determinant than the numerical integration method chosen. For instance, in the globalization of the slow manifold, one starts at a point $x_0 = K(\theta, 0, \sigma_2^{(0)})$, but numerically, it corresponds to a value $K(\theta, \varepsilon_1, \sigma_2^{(0)})$, where $\varepsilon_1$ is a small error. Since $|\lambda_1| > |\lambda_2|$, this error may grow as we integrate backwards in time
\begin{equation}
\phi_{\substack{-\Delta t}}(x_0) = K(\theta - \Delta t/T, \varepsilon_1 e^{-\lambda_1 \Delta t},  \sigma^{(0)}_2e^{-\lambda_2 \Delta t}),
\end{equation}
and thus generate isochrons not accurate enough (it is not possible to obtain a homogeneous coverage of the space, specially for points near the manifold $\Sigma_1(x) = 0$).

Thus, a good numerical approximation provides a smaller $\varepsilon_1$ and a high order expansion for $\bar{K}$ provides local points that are far from the invariant object 
$\Gamma$, thus reducing the integration time $\Delta t$ and therefore the error. Moreover, in the algorithms presented, the time $\Delta t$ is always a multiple of the period $T$, since we globalize the isochron by integrating backwards points on the same isochron. However, one can adapt the algorithms to globalize the isochrons using points on a different isochron (see Eq.~\eqref{eq:isoback}) and therefore integrating for a shorter time $\Delta t$ (see \cite{guillamon2009computational}).

%Therefore, the larger the isochron one aims to compute, the more numerical precision one needs to use. Other strategies to improve this precision is to use double-precision digits instead of float and large precision functions and integrators such as Taylor methods (see \cite{}). 
\end{remark}

\section{Numerical Examples}\label{sec:section4}

We have carried out the numerical implementation of the methodology presented in Section \ref{sec:section3} and applied it to some representative single neuron and neuronal population models showing oscillatory dynamics. In this Section we present the results obtained and discuss the relevant aspects. %\textcolor{red}{The examples will show that the geometric objects (isochrons and isostables) provide information on the transient dynamics far from the limit cycle.} 

%and has a limit cycle with two characteristic exponents with different orders of magnitude and one of them very small.

\subsection{Models}

In this Section, we present the four models that we have used to apply the techniques introduced in Section \ref{sec:section3}. 
There are two models for spiking dynamics of single neurons, including a version of the classical Hodgkin-Huxley model, and two models for mean firing rate dynamics of neuronal populations, including a 3D version of the classical Wilson-Cowan model. The examples are chosen either because they are 
relevant to illustrate the properties of the numerical methodology or because they are representative of classical models.
The parameter values and the functions for each model can be found in Appendix \ref{sec:appendix}. 

% In this section, we apply our method to representative examples, rangingfrom the most simple instances of Hopf and SNIC (saddle-node on an invariant curve) bifur-cations and the classical van der Pol oscillator to more sophisticated neuronal models. Apartfrom obtaining isochrons, PRCs, and PRSs, through these examples we want to illustrate dif-ferent facts: (a) what are the clues to explain the transition from “Type I” PRCs to “Type II”PRCs; (b) the numerical problems that arise when dealing with slow-fast systems; and, (c) upto which degree PRSs show disagreement with PRCs in the same phase and how this can affecthigh frequency stimulation. We end the paper with a discussion on these facts in section6.

%
\begin{itemize}
	
	\item A single thalamic neuron model introduced in \cite{rubin2004high}, that we refer to as $RT$, with sodium, potassium and low-threshold calcium currents:
	\begin{equation}
	\begin{aligned}\label{eq:rtEDOs}
	C_m \dot{V} &= - I_L(V) - I_{Na}(V, h) - I_{K}(V, h) - I_{T}(V, r) + I_{app}, \\
	\dot{h} &= \frac{h_{\infty}(V) - h}{\tau_h(V)}, \\
	\dot{r} &= \frac{r_{\infty}(V) - r}{\tau_r(V)},
	\end{aligned}
	\end{equation}	
	where $V$ describes the membrane potential and $h$ and $r$ are the gating variables. 
	
	This model is interesting because it has a slow-fast dynamics (notice that as $\tau_r(V) \gg \tau_h(V) \gg 1$, the variable $r$ is much slower than $V$ and $m$) and it allows us to illustrate how the slow-fast dynamics affects the geometric objects (isochrons and isostables). Moreover, as we will see in Section~\ref{sec:section5}, we can take advantage of these dynamics to explore the amplitude-phase description as an alternative to the phase reduction. This model was also studied in \cite{wilson2018greater}. 
	
	\item A reduced Hodgkin-Huxley-like system that we refer to as $HH$ \cite{izhikevich2007}, with sodium and potassium currents, and two gating variables:
	\begin{equation}
	\begin{aligned}\label{eq:hhEDOs}
	C_m \dot{V} &= - I_L(V) - I_{Na}(V, h) - I_{K}(V, n) + I_{app}, \\
	\dot{n} &= \frac{n_{\infty}(V) - n}{\tau_n(V)}, \\
	\dot{h} &= \frac{h_{\infty}(V) - h}{\tau_h(V)},
	\end{aligned}
	\end{equation}	
	where $V$ describes the membrane potential and $n$ and $h$ are the gating variables.
	
	This model is a 3D reduction of the classical 4D Hodgkin-Huxley model (just setting the fast variable $m$ to its steady-state value and slightly modifying some parameters). It also has a 
	slow-fast nature but less noticed than in the $RT$ model.

	\item An extension of the Wilson-Cowan equations \cite{wilson1972excitatory} including dynamics for the inhibitory synapses that we refer to as $WC_{Syn}$:
	\begin{equation}
	\begin{aligned}\label{eq:wcEDOs3D}
	\tau_e \dot{E} &= -E + \delta_E(c_1 E - c_2 s + P), \\
	\tau_i \dot{I} &= -I + \delta_I(c_3 E - c_4 s + Q), \\
	\tau_d \dot{s} &= -s + \tau_d I,
	%\dot{s} &= -\frac{s}{\tau_d} + \frac{I}{\tau_r},
	\end{aligned}
	\end{equation}
	where $E$ and $I$ are the mean firing rates of excitatory and inhibitory populations, respectively, whereas $s$ describes the inhibitory synaptic dynamics.

	\item A model for the mean field activity of a population of heterogeneous quadratic integrate-and-fire neurons  \cite{devalle2017firing} that we refer to as $QIF$:
	\begin{equation}
	\begin{aligned}\label{eq:qfEDOs}
	\tau_m \dot{V} &= V^2 - (\pi \tau_m R)^2 - J\tau_mS + \Theta, \\
	\tau_m \dot{R} &= -\frac{\Delta}{\pi \tau_m} + 2RV,\\
	\tau_d \dot{S} &= -S + R,
	\end{aligned}
	\end{equation}	
	where $V$ is the mean membrane potential, $R$ is the mean firing rate of the population  and $S$ is the synaptic activation.	
	
	We have chosen this model because it is a representative example of a new generation of neural field models \cite{montbrio2015macroscopic, coombes2019next}.

\end{itemize}

\subsection{Numerical results}\label{sec:section42}

In the four examples considered, there exists a hyperbolic attracting limit cycle $\Gamma$. Moreover, there exists an unstable equilibrium point of saddle-focus type having a two-dimensional unstable manifold and a one-dimensional stable manifold. We set the zero phase value at the maximum of the voltage variable $V$ (in the $WC_{Syn}$ model it will be the maximum of the $E$ variable). For each model, the limit cycle $\Gamma$ has two characteristic exponents, $\lambda_1$ and $\lambda_2$, with $|\lambda_1| > |\lambda_2|$. The numerical integration has been performed using an 8th-order Runge-Kutta Fehlberg method (rk78) with a tolerance of $10^{-14}$. In the neighbourhood of $\Gamma$ we have computed a Taylor expansion up to order $L=10$ of the parameterization $K$ as in \eqref{eq:truncatedSeries3D}  and we have computed $N+1$, with $N=2048$, Fourier coefficients for the periodic functions $K_{\alpha,m-\alpha}$, $\alpha=0,\ldots,m$ and $m=0,\ldots,L$. Notice that, as in all the examples the dimension $d=3$, to obtain expansions up to an order $L$ we use $\sum_{n=1}^{L+1} n=(L+1)(L+2)/2 = 60$ monomials. Recall that the number of Fourier coefficients chosen was such that the residuals are smaller than $E_{tail}=10^{-10}$ (see Eq. \eqref{eq:theTail}). The domain of the local approximation $\Omega_{loc}$ is computed with an error smaller than $E_{loc}=10^{-8}$, except for HH that was $10^{-6}$ (see Eq.~\eqref{eq:mjOmegaLoc}). We emphasize that the numerical computation of $\bar{K}$ for the above mentioned values of $L$ and $N$, only takes a few seconds (around 10s)
	on a regular laptop because of the combination of automatic differentiation techniques and the diagonalization achieved by means of the Floquet normal form. Table \ref{table:modelPrms} contains the values of the parameters described above for each model.\\

\begin{table}[H]
	\renewcommand{\arraystretch}{1.5}
	\begin{center}
		\begin{tabular}{|c|c|c|c|c|}
			\hline
			Model         & $T$ & $\lambda_1$ & $\lambda_2$ & %$L $ & $N $ & $E_{tail}$ & $E_{loc}$ & 
			Equilibrium Point \\
			\hline
			$RT$         & 8.395 & $-0.368$ & $-0.022$ & %10 & 2048 & $10^{-10}$ & $10^{-8}$ & 
			(-39.1, 0.38, 1.3 $10^{-5}$) \\
			\hline
			$HH$         & 7.586 & $-1.73$ & $-0.2$ & %10 & 2048 & $10^{-10}$ & $10^{-6}$ & 
			(-49.1, 0.564,   0.137)\\
			\hline
			$WC_{Syn}$ & 24.43 & $-0.445$ & $-0.246$ & %10 & 2048 & $10^{-10}$ & $10^{-8}$ & 
			(0.272, 0.033, 0.198) \\
			\hline
			$QIF$ & 27.58 & $-0.408$ & $-0.06$ & %10 & 2048 & $10^{-10}$ & $10^{-8}$ & 
			(0.018, -0.267, 0.018) \\
			\hline
		\end{tabular}
	\end{center}
	\renewcommand{\arraystretch}{1}
	\caption{Numerical values for the different models considered; $T$ period of the periodic orbit $\Gamma$; $\lambda_1$, $\lambda_2$ characteristic exponents associated to $\Gamma$; and coordinates of the unstable equilibrium point of saddle-focus type.}
	\label{table:modelPrms}
\end{table}

\begin{figure}
	\centering
	{\includegraphics[width=160mm]{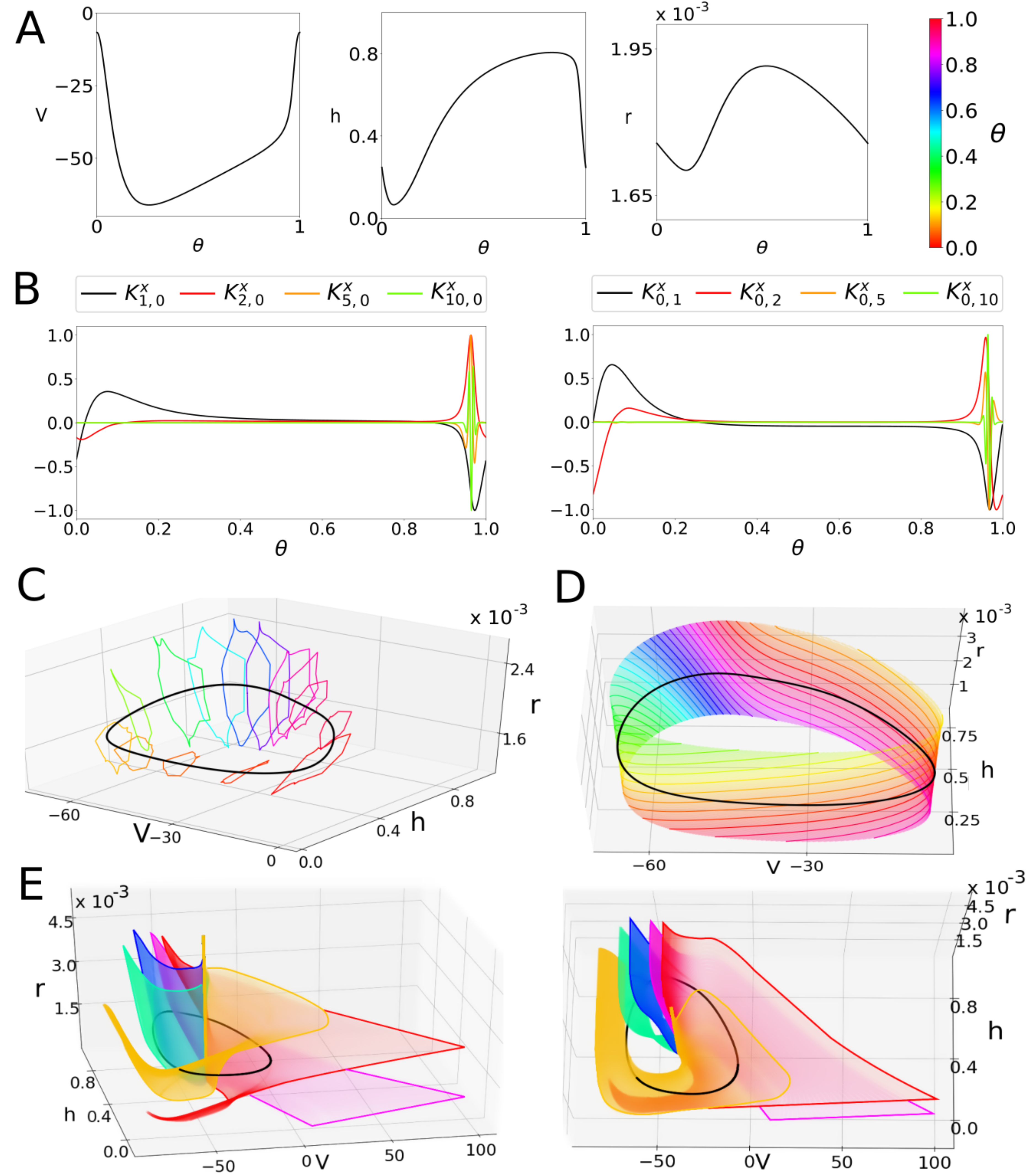}}
	\caption{For the $RT$ model \eqref{eq:rtEDOs} we show: 
	(A) Coordinates $V, h,$ and $r$ of the parameterization $\gamma(\theta)$ in \eqref{eq:mathDef_2} of the periodic orbit. 
	(B) First coordinate of the coefficients $\bar{K}_{i,0}$ and $\bar{K}_{0,i}$ for $i=1,2,5,10$ of the local approximation $\bar K$ normalized by the maximum. 
	(C) Local isochrons $\mathcal{I}^{loc}_{\theta_i}$ for $\theta_i=i/2048$ and $i=0, 64, 128, 192, 256, 512, 768, 1024, 1280, 1536, 1792, 1856, 1920, 1984$. 
	(D) Slow attracting 2D manifold $\mathcal{S}$ of $\Gamma$, corresponding to the isostable $\mathcal{A}^1_0$, and 1D leaves $\mathcal{S}_{\theta_i}$ (solid curves) for $\theta_i=i32/2048$ and $i=0,\ldots,63$. (E) Two perspectives of the globalized isochrons $\mathcal{I}_{\theta_i}$ corresponding to the phases $\theta_i$ with $i=0, 256, 1024, 1280, 1792$. In panels CDE the black curve corresponds to the limit cycle.} \label{fig:rtPanel}
	%= 0.0, 0.125, 0.5, 0.625, 0.875
\end{figure}

\begin{figure}
	\centering
	{\includegraphics[width=160mm]{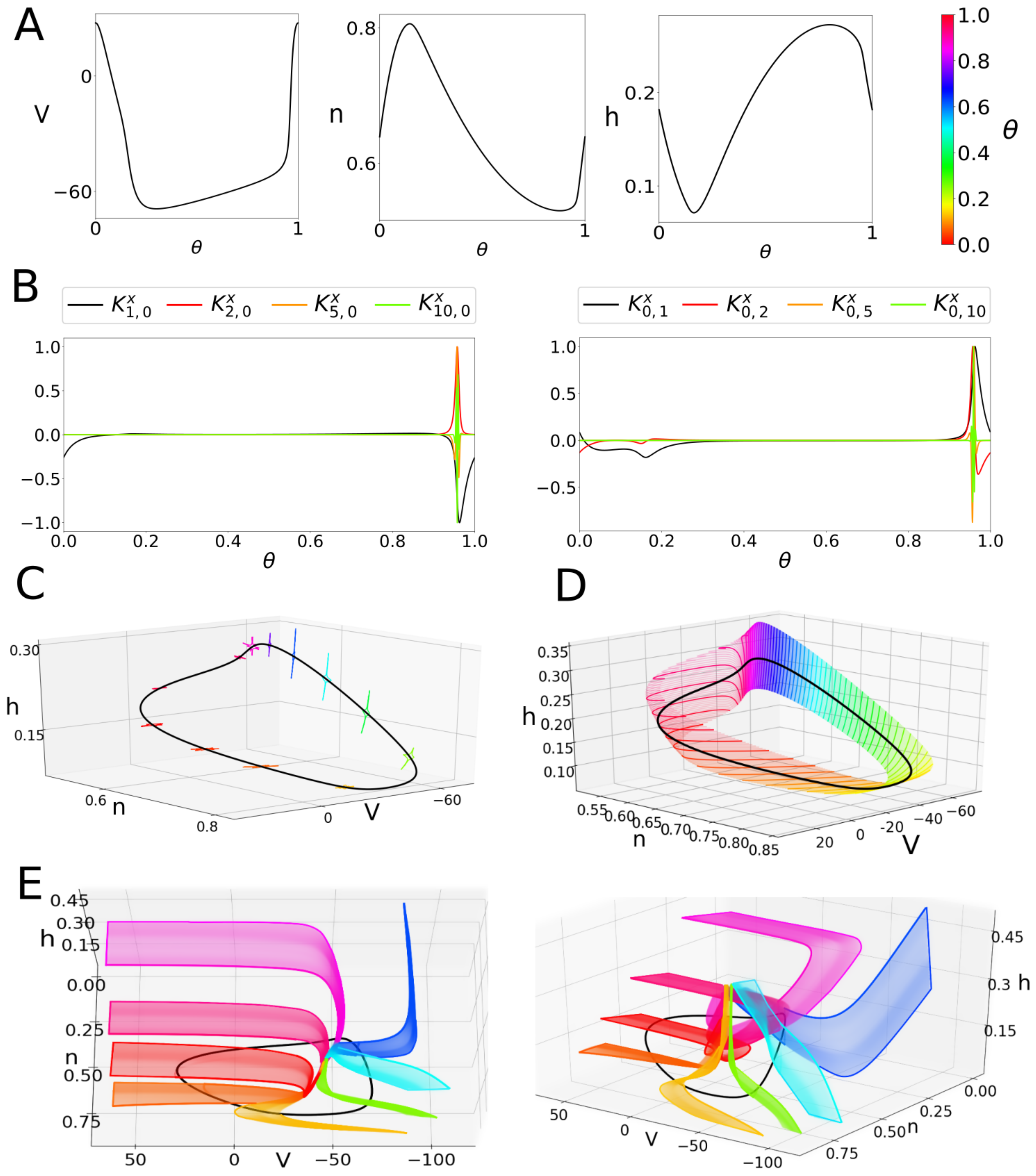}}
	\caption{For the $HH$ equations \eqref{eq:hhEDOs} we show: 
	(A) Coordinates $V, n,$ and $n$ of the parameterization $\gamma(\theta)$ in \eqref{eq:mathDef_2} of the periodic orbit. 
	(B) First coordinate of the coefficients $\bar{K}_{i,0}$ and $\bar{K}_{0,i}$ for $i=1,2,5,10$ of the local approximation $\bar K$ normalized by the maximum. 
	(C) Local isochrons $\mathcal{I}^{loc}_{\theta_i}$ for $\theta_i=i/2048$ and $i=0,    64,   128,   256,   512,   768,  1024,  1280,
        1536,  1792,  1920,  1984$. 
        (D) Slow attracting 2D manifold $\mathcal{S}$ of $\Gamma$, corresponding the isostable $\mathcal{A}^1_0$, and 1D leaves $\mathcal{S}_{\theta_i}$ (solid curves) for $\theta_i=i32/2048$ and $i=0,\ldots,63$. (E) Two perspectives of the globalized isochrons $\mathcal{I}_{\theta_i}$ corresponding to the phases $\theta_i$ with $i=0,   128,   256,   512,  1024,  1280,  1792,  1920$. In panels CDE the black curve corresponds to the limit cycle.
	%$\theta = 0.0, 0.0625, 0.125, 0.25, 0.5, 0.625, 0.875, 0.9375$
	} \label{fig:hhPanel}
\end{figure}

\begin{figure}
	\centering
	{\includegraphics[width=160mm]{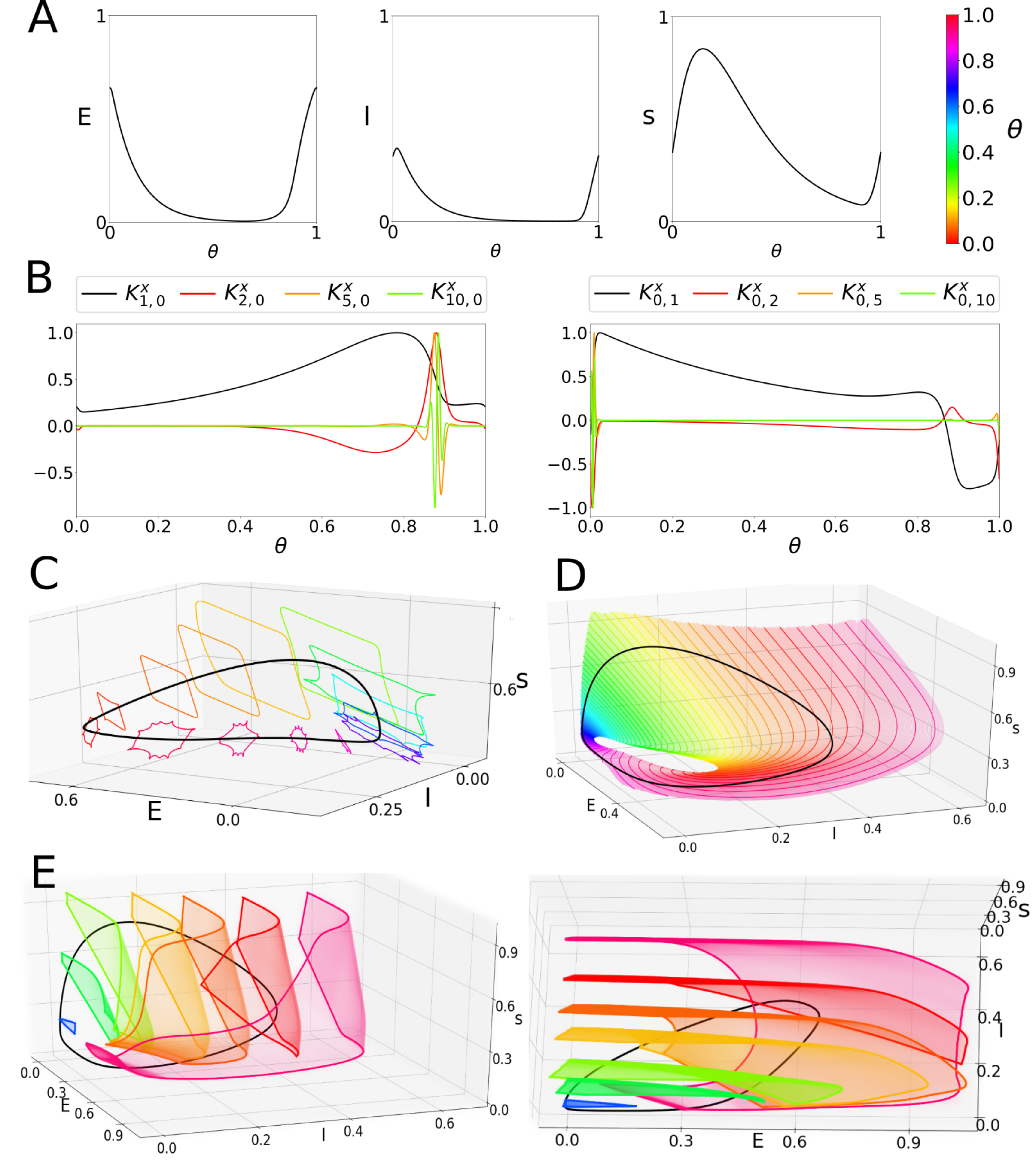}}
	\caption{For the $WC_{Syn}$ equations \eqref{eq:wcEDOs3D} we show: 
	(A) Coordinates of $E, I,$ and $s$ of the parameterization $\gamma(\theta)$ in \eqref{eq:mathDef_2} of the periodic orbit.
	(B) First coordinate of the coefficients $\bar{K}_{i,0}$ and $\bar{K}_{0,i}$ for $i=1,2,5,10$ of the local approximation $\bar K$ normalized by the maximum.  
	(C) Local isochrons $\mathcal{I}^{loc}_{\theta_i}$ for $\theta_i=i/2048$ and $i=0, 64, 128, 192, 256, 384, 512, 768, 1280, 1792, 1920, 1984 $. (D) Slow 2D attracting manifold $\mathcal{S}$ of $\Gamma$, corresponding the isostable $\mathcal{A}^1_0$, and 1D leaves $\mathcal{S}_{\theta_i}$ (solid curves) for $\theta_i=i32/2048$ and $i=0,\ldots,63$. (E) Two perspectives of the globalized isochrons $\mathcal{I}_{\theta_i}$ corresponding to the phases $\theta_i$ with $i=0,   128,   256,   512,   768,  1280,  1920$. In panels CDE the black curve corresponds to the limit cycle.
	%$\theta = 0.0, 0.0625, 0.125, 0.25, 0.375, 0.625, 0.9375$.
	} \label{fig:wcPanel}
\end{figure}

\begin{figure}
	\centering
	{\includegraphics[width=160mm]{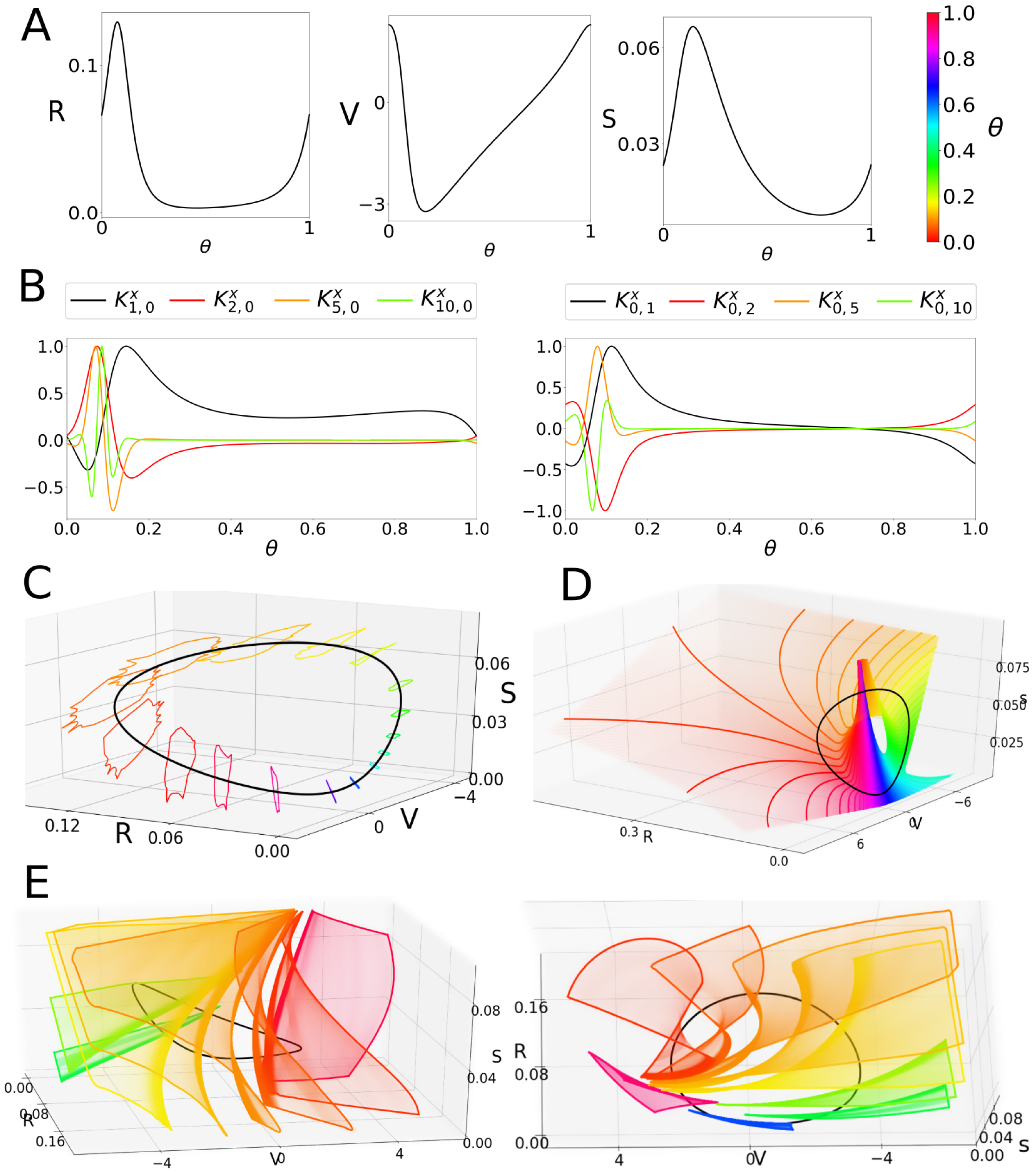}}
	\caption{For the $QIF$ model \eqref{eq:qfEDOs} we show: 
	(A) Coordinates of $R, V,$ and $s$ of the parameterization $\gamma(\theta)$ in \eqref{eq:mathDef_2} of the periodic orbit.
	(B) First coordinate of the coefficients $\bar{K}_{i,0}$ and $\bar{K}_{0,i}$ for $i=1,2,5,10$ of the local approximation $\bar K$ normalized by the maximum. 
	(C) Local isochrons $\mathcal{I}^{loc}_{\theta_i}$ for $\theta_i=i/2048$ and $i=0, 64, 128, 192, 256, 320, 384, 512, 640, 768, 896, 1024, 1280, 1536, 1856, 1984$. (D) Slow attracting 2D manifold $\mathcal{S}$ of $\Gamma$, corresponding the isostable $\mathcal{A}^1_0$, and 1D leaves $\mathcal{S}_{\theta_i}$ (solid curves) for $\theta_i=i32/2048$ and $i=0,\ldots,63$. (E) Two perspectives of the globalized isochrons $\mathcal{I}_{\theta_i}$ corresponding to the phases $\theta_i$ with $i=128,    64,   192,   256,   320,   512,   768,  1984$. In panels CDE the black curve corresponds to the limit cycle.
	%$\theta = 0.0625, 0.03125, 0.09375, 0.125, 0.15625, 0.25, 0.375, 0.96875$.
	} \label{fig:qfPanel}
\end{figure}

For each model, we have computed the approximate parameterization $\bar{K}$ \eqref{eq:truncatedSeries3D}, the slow manifold $\mathcal{S}$ \eqref{eq:slowManifold} (which is the most relevant isostable) and the isochrons $\mathcal{I}_{\theta}$ \eqref{eq:isochronsDef}. In Figs.~\ref{fig:rtPanel} to \ref{fig:qfPanel}, we illustrate some aspects of the elements computed, but we stress that our computations have more data than the ones shown. More precisely, each figure corresponds to a different model and has 5 panels that contain the same information for each model. Panel A shows the three components of the parameterization of the limit cycle $\gamma(\theta)$ defined in \eqref{eq:mathDef_2}. Panel B shows the first component of the periodic functions $K_{i,0}$ and $K_{0,i}$ for 
$i=1,2,5,10$, corresponding to different coefficients of the Taylor expansion of $\bar K$ up to order 10. As the magnitude of the different functions $K_{i,j}$ is very different and depends on the constants  $b_{1}$ and $b_{2}$ (see Section~\ref{sec:section32}), for illustration purposes we have normalized these functions so that the maximum is one (see Table \ref{table:ksPrms} for more details). Panel C shows the domain of accuracy $\Omega_{loc}$ of the local approximation $\bar K$ for some values of $\theta$ or, equivalently, the local isochrons $\iloc$. We have computed $N=2048$ isochrons, that is $\mathcal{I}^{loc}_{\theta_i}$ for $\theta_i=i/2048$ for $i=0,\ldots,2047$. We have not chosen equidistant values of $\theta$ to plot the isochrons. Alternatively, the chosen values of $\theta_k$ have been adapted to each model, so that the isochrons are uniformly distributed in the phase space along the limit cycle to allow the reader to distinguish the isochrons in the plots.  Panel D shows the globalized slow manifold $\mathcal{S}$ \eqref{eq:slowManifold} and 64 leaves $\mathcal{S}^{\theta}$ \eqref{eq:sltheta}, for equidistant values of $\theta$, more precisely for $\theta_i=i 32/2048$, for $i=0,..63$. Finally, panel E shows the globalization of some of the isochrons $\iiso$ computed in panel C from two different perspectives.  The isochron computation has been restricted to a domain $\Omega_c$, based on values that are biophysically plausible. Namely, for the $RT$ model, the gating variables $h$ and $r$ are allowed to vary between 0 and 1. The same criterion is applied to the gating variables $n$ and $h$ in $HH$ model, and $V$ was restricted to $V<60$. For the $WC_{Syn}$ model, the three variables $E, I$ and $s$ are restricted between 0 and 1. Finally, for the $QIF$ model, $R$ and $S$ variables can not be negative, and $-6 < V < 6$. 

\begin{figure}
	\centering
	{\includegraphics[width=160mm]{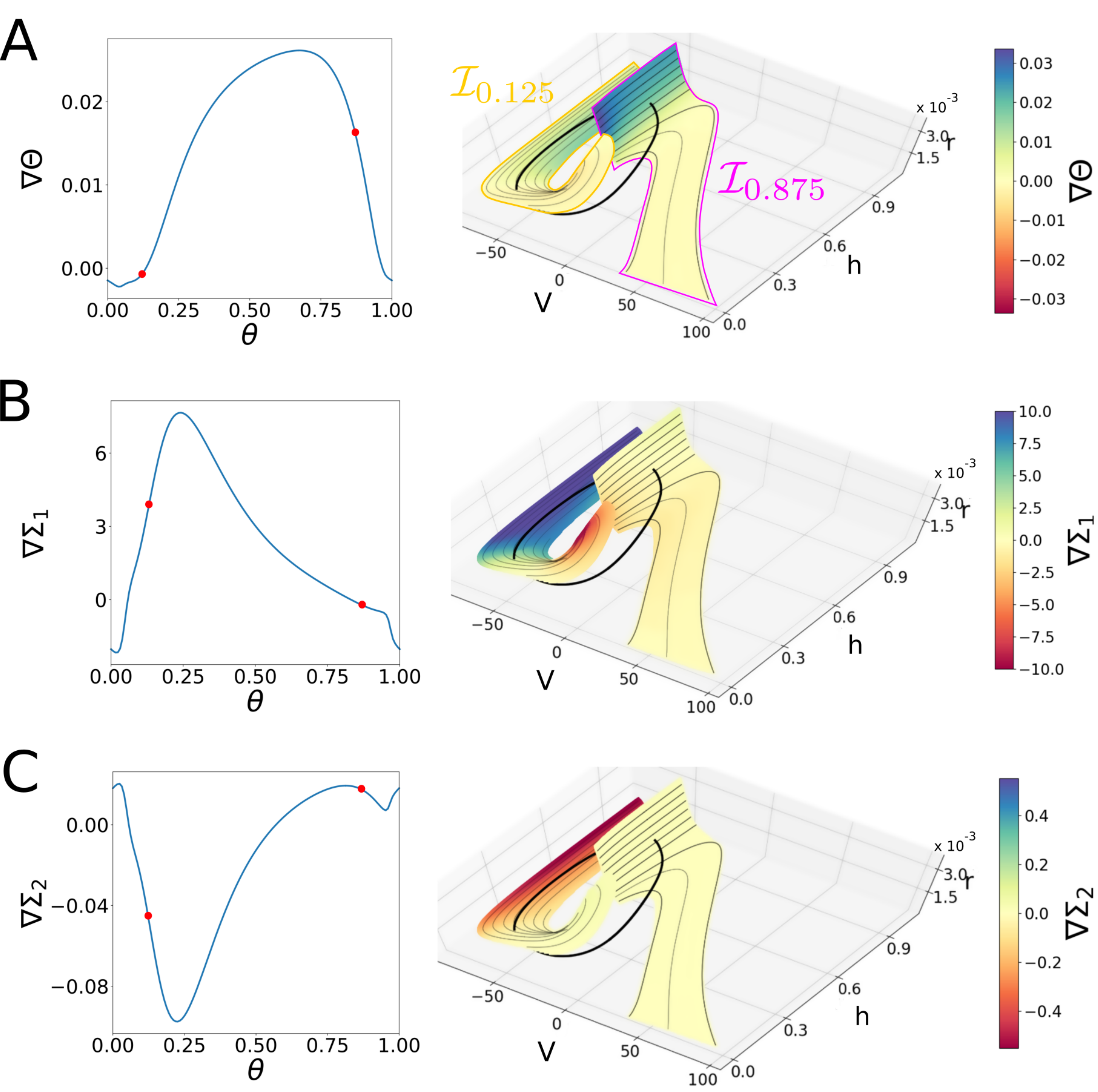}}
	\caption{For the $RT$ model \eqref{eq:rtEDOs} we show the $V$ component of the iPRF $\nabla \Theta$ (A) and the iARFs $\nabla \Sigma_i, i=1,2$ (B,C) (that abusing of notation we denote by the same letters), restricted to the limit cycle, i.e. iPRC and iARCs (left) and evaluated on the isochrons $\mathcal{I}_\theta$, for $\theta = 0.125, 0.875$ (right). The border of each isochron is highlighted with the same colour used to plot the isochron in Panel E in Fig.~\ref{fig:rtPanel}. The red dots on the curves of the left column correspond to the phases of the isochrons on the right column.} \label{fig:rtPanelPrf}
\end{figure}

\begin{figure}
	\centering
	{\includegraphics[width=160mm]{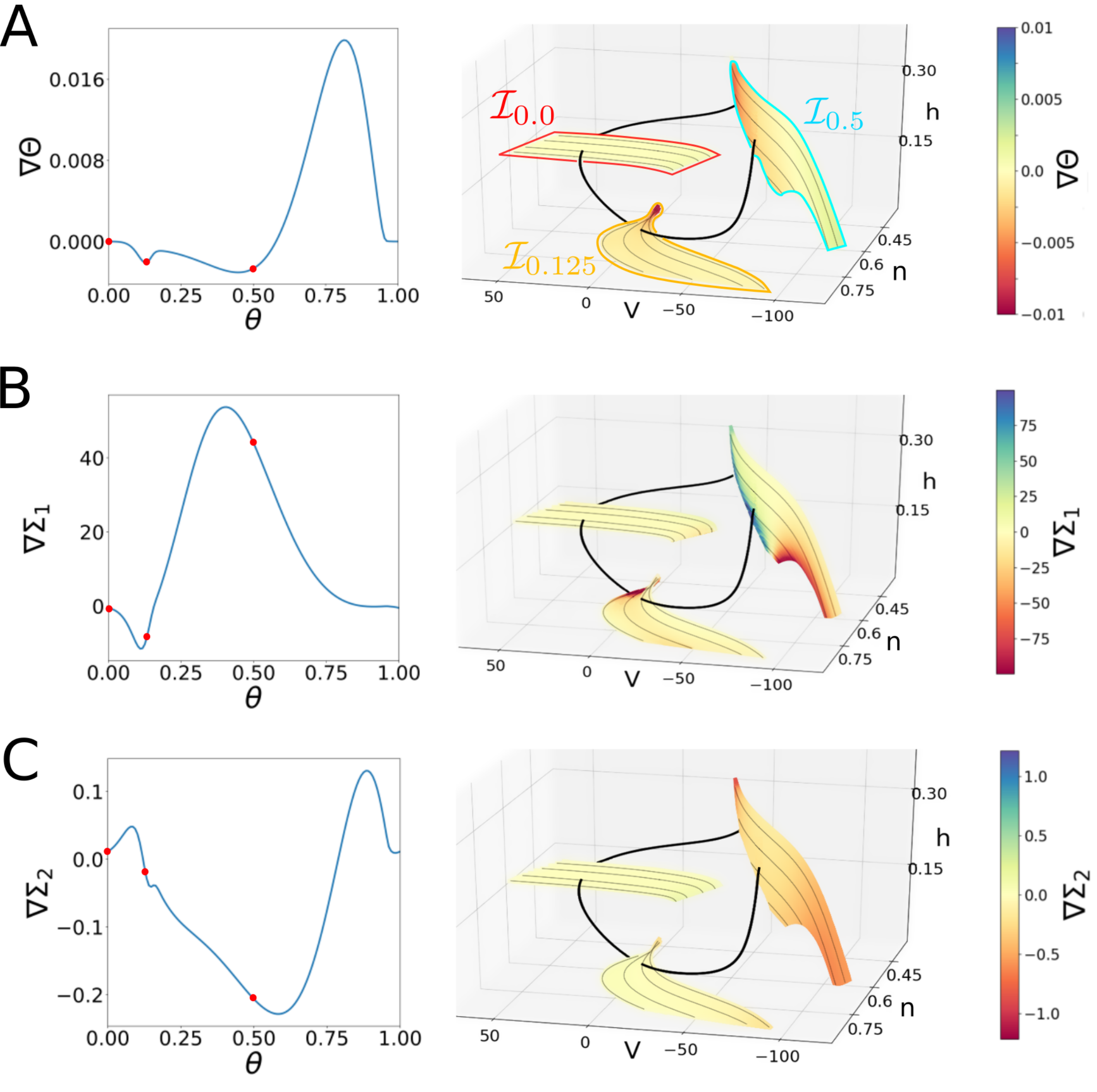}}
	\caption{
	For the $HH$ model \eqref{eq:hhEDOs} we show the $V$ component of the iPRF $\nabla \Theta$ (A) and the iARFs $\nabla \Sigma_i, i=1,2$ (B,C) (that abusing of notation we denote by the same letters), restricted to the limit cycle, i.e. iPRC and iARCs (left) and evaluated on the isochrons $\mathcal{I}_\theta$, for $\theta = 0.0, 0.125, 0.5$ (right). The border of each isochron is highlighted with the same colour used to plot the isochron in Panel E in Fig.~\ref{fig:hhPanel}. The red dots on the curves of the left column correspond to the phases of the isochrons on the right column. } \label{fig:hhPanelPrf}
\end{figure}

\begin{figure}
	\centering
	{\includegraphics[width=160mm]{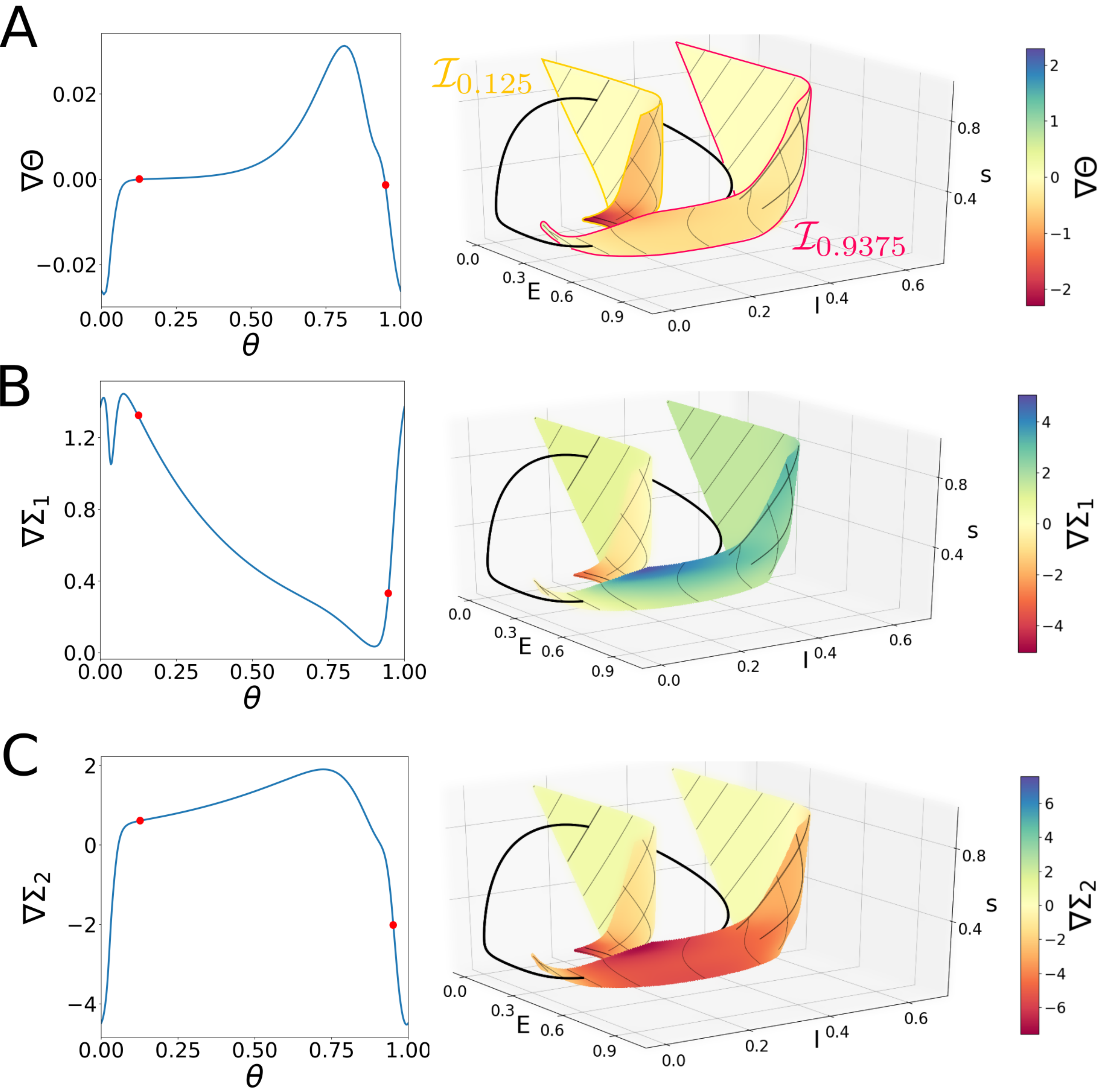}}
	\caption{For the $WC_{Syn}$ model \eqref{eq:wcEDOs3D} we show the $E$ component of the iPRF $\nabla \Theta$ (A) and the iARFs $\nabla \Sigma_i, i=1,2$ (B,C) (that abusing of notation we denote by the same letters), restricted to the limit cycle, i.e. iPRC and iARCs (left) and evaluated on the isochrons $\mathcal{I}_\theta$, for $\theta = 0.125, 0.9375$ (right). The border of each isochron is highlighted with the same colour used to plot the isochron in Panel E in Fig.~\ref{fig:wcPanel}. The red dots on the curves of the left column correspond to the phases of the isochrons on the right column.} \label{fig:wcPanelPrf}
\end{figure}

\begin{figure}
	\centering
	{\includegraphics[width=160mm]{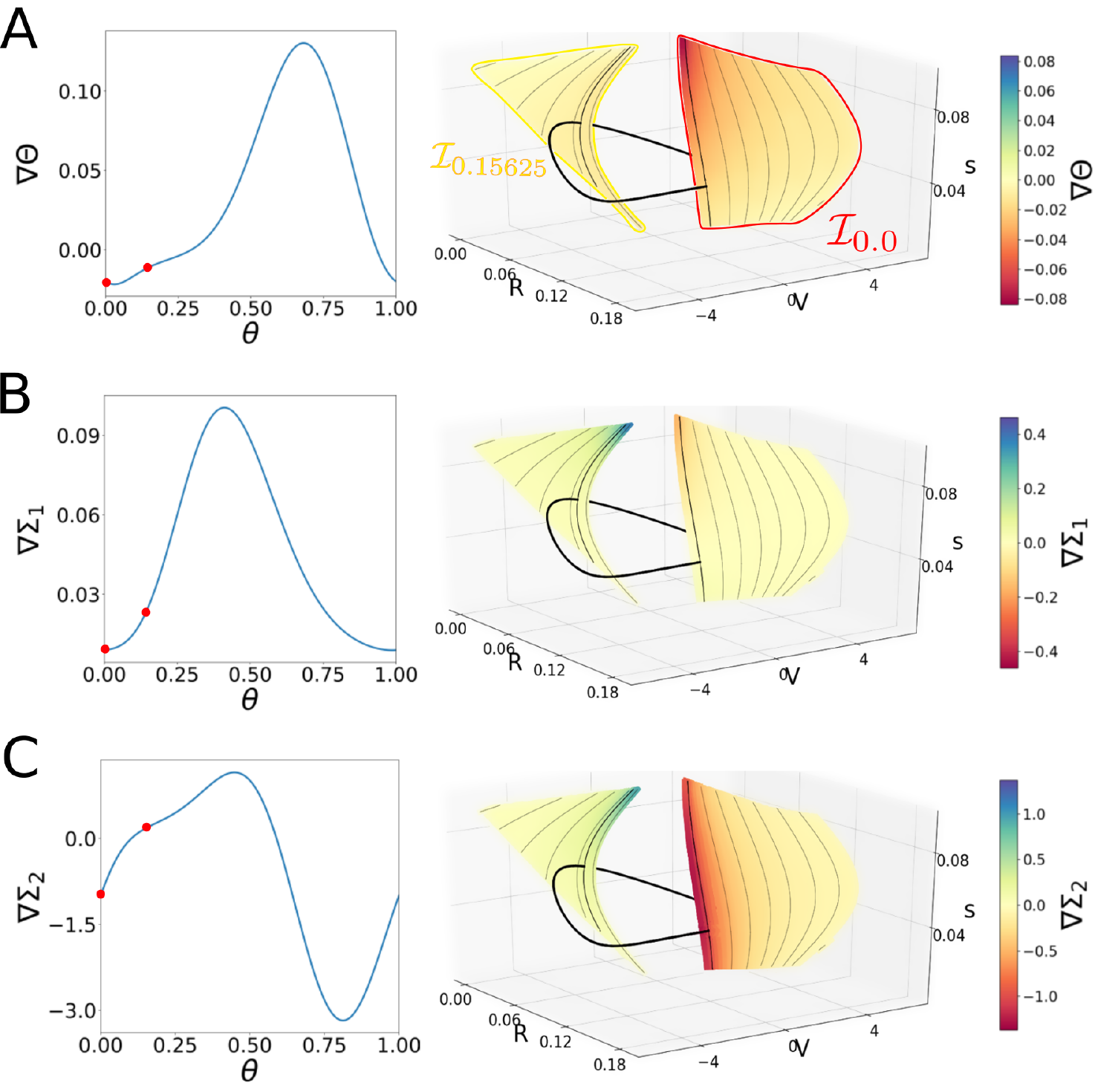}}
	\caption{For the $QIF$ model \eqref{eq:qfEDOs} we show the $V$ component of the iPRF $\nabla \Theta$ (A) and the iARFs $\nabla \Sigma_i, i=1,2$ (B,C) (that abusing of notation we denote by the same letters), restricted to the limit cycle, i.e. iPRC and iARCs (left) and evaluated on the isochrons $\mathcal{I}_\theta$, for $\theta = 0.0, 0.15625$ (right). The border of each isochron is highlighted with the same colour used to plot the isochron in Panel E in Fig.~\ref{fig:qfPanel}. The red dots on the curves of the left column correspond to the phases of the isochrons on the right column. } \label{fig:qfPanelPrf}
\end{figure}

Moreover, we have also computed for each model the iPRF $\nabla \Theta$ and the iARFs $\nabla \Sigma_i, i=1,2$ (Figs.~\ref{fig:rtPanelPrf} to \ref{fig:qfPanelPrf}). In each figure on the right, we show the first component of these vector-valued functions that, abusing language, we denote by the same letters (see Remark~\ref{rem:sobreprf}). On the left, we plot the restriction of these functions to the limit cycle, the iPRC and iARCs (see Remark~\ref{rm:prcsRM}).  Since the computation of the functions $\nabla \Theta(x)$, $\nabla \Sigma_1(x)$ and $\nabla \Sigma_2(x)$ is done in parallel with the computation of points on the isochrons (see Section~\ref{sec:isoComp3D}), we have evaluated these functions for those computed points on the isochrons. For illustration purposes, we plot these functions for two or three isochrons indicated in each caption. Moreover, in the figures we just plot a part of the isochron that permits a proper visualization of the functions $\nabla \Theta$ and $\nabla \Sigma_i$, $i=1,2$, even though our computations are done for the whole globalized isochrons.

%In the four Figures, panel A shows the time series during a period. Notice how mechanisms across which the zero phase (which is or can be thought as a neuronal or population spike in all four cases) is reach is different. For the RT model, the spike is reached by the opening of both channels $h$ and $r$. So both channels contribute to the spike rise. By contrast the other neuron model, the $HH$, the opening of $h$ channels leads to a sudden rise of the membrane voltage $V$ which lead to the opening of $n$ channels which cause a sadden fall. In this case, a channel leads to the opening and the other to the clausure. 

%

\begin{table}[H]
	\renewcommand{\arraystretch}{1.5}
	\begin{center}
		\begin{tabular}{|c|c|c|c|c|c|c|c|c|c|c|}
			\hline
			Model  & $b_{1}$    & $K^{max}_{1,0}$ & $K^{max}_{2,0}$ & $K^{max}_{5,0}$ & $K^{max}_{10,0}$ & $b_{2}$ & $K^{max}_{0,1}$ & $K^{max}_{0,2}$ & $K^{max}_{0,5}$ & $K^{max}_{0,10}$ \\
			\hline
			$RT$  &  0.5   & 1.2 & 0.05 & $7\cdot10^{-6}$ & $1.5\cdot10^{-11}$ & 0.5 & 51.2 & 59.5 & $1.2\cdot10^{3}$ & $1.2\cdot10^{5}$ \\
			\hline
			$HH$   & 2     & 7.8 & 1.2 & $1.3\cdot10^{-2}$ & $1.2\cdot10^{-5}$ & 2 & 23.7 & 9.1 & 3.3 & 3.2 \\
			\hline
			$WC_{Syn}$    & 1     & 3.4 & 7.8 & $3\cdot10^{3}$ & $4\cdot10^{8}$ & 1 & 0.54 & 3.64 & $1.2\cdot10^{3}$ & $1\cdot10^{8}$ \\
			\hline
			$QIF$  & 0.2     & 4.58 & 5.9 & 8.6 & 18.1 & 1 & 2.3 & 2.4 & 3.05 & 4.38 \\
			\hline
		\end{tabular}
	\end{center}
	\renewcommand{\arraystretch}{1}
	\caption{Numerical values for the constants $b_1$, $b_2$ and $K^{max}_{i,j}=\max_{\theta \in \mathbb{T}} |\bar{K}^X_{i,j}(\theta)|$, where we denote by $\bar{K}^X_{i,j}(\theta)$ the first component of the function $\bar{K}_{i,j}(\theta)$.}
	\label{table:ksPrms}
\end{table}

As we observed in Remark \ref{rm:remark32}, the choice of constants $b_{1}$, $b_{2} \in \mathbb{R}$ affects the magnitude of the monomials $K_{\alpha, m-\alpha}$ and therefore the domain of accuracy of $\bar K$. Indeed, depending on the choice of $b_{1}$, $b_{2}$, successive orders of $K_{\alpha, m-\alpha}$ can either blow up or vanish, since the monomial $K_{\alpha, m-\alpha}$ will be multiplied by a factor $b_{1}^\alpha b_{2}^{m-\alpha}$. In Table~\ref{table:ksPrms} we provide the values of the constants $b_{1}$ and $b_{2}$ that we have used for each model as well as the maximum values of some representative monomials. The criteria for the choice has been to keep the monomials at order one (as much as possible), thus maximizing the domain of accuracy of the local approximation $\bar K$. 

\newpage

For all models we have computed the Taylor expansion up to the same degree $L=10$ and the same error tolerances $E_{loc}=10^{-8}$ and $E_{tail}=10^{-10}$. For the $QIF$ model the monomials can be kept of order 1 up to degree 10, but for the $HH$ and $RT$ models some monomials decay to zero at degree 5. Thus, beyond degree 5, for the $HH$ and $RT$ models, adding more monomials in the Taylor expansion does not allow us to increase the size (in the Euclidean metrics) of the domain $\Omega_{loc}$ and, in consequence, the size of the local isochrons $\mathcal{I}^{loc}_\theta$. Indeed, when we compare the sizes of the local isochrons for the different models, they are smaller for the RT and HH models than the others (compare panel C in for the $RT$ and $HH$ model with panel C for the $QIF$ model). Moreover, the isochrons for the $HH$ model look elongated (panel C in Fig.~\ref{fig:hhPanel}). Indeed, in the direction of $\sigma_1$, the coefficients $\bar{K}_{ij}$ decay to  0, while in the direction of $\sigma_2$ they remain of order 1 (see Table~\ref{table:ksPrms}).

Besides the purely numerical considerations regarding the coefficients $\bar{K}_{ij}$, we want to highlight that the linear terms $\bar{K}_{10}$ and $\bar{K}_{01}$ (black curves in panel B of Figs~\ref{fig:rtPanel}-\ref{fig:qfPanel}) provide information about the attraction to the limit cycle. More precisely, when the attraction to the limit cycle is not homogeneous along the cycle, the functions $\bar{K}_{10}$ and $\bar{K}_{01}$ show dramatic changes along a period, as it happens for instance for the $RT$ and $HH$ models (see Figs.~\ref{fig:rtPanel}B and \ref{fig:hhPanel}B). Then, this has also consequences in the size and shape of the local isochrons. Take, for instance, the values of the coefficients $\bar{K}_{ij}$ for the first component $\bar{K}^X_{ij}$ shown in panels B.  The phases $\theta$ at which $\bar{K}^X_{ij}$ are close to zero (approximately $0.4-0.8$ in Fig.~\ref{fig:rtPanel}B and $0.2-0.8$ in Fig.~\ref{fig:hhPanel}B), the corresponding local isochron $\iloc$ does not extend along the first component $V$ and it is almost parallel to planes with constant $V$ (see the isochrons in blue colours Fig.~\ref{fig:rtPanel}C and in blue/green colours in Fig.~\ref{fig:hhPanel}C). For the same reason, if we look at the 3 components of the vector-valued functions $\bar{K}_{ij}$, for those values of $\theta$ for which the 3 components are close to zero, the local isochrons will be small. 
Thus, local isochrons $\iloc$ show great size and shape differences for different values of $\theta$, whenever the range of the functions $\bar{K}_{ij}$ is wide. However, whenever $\bar{K}_{ij}$ are more uniform along a period the size of the local isochrons is also more homogeneous for the different phases (see for instance Fig.\ref{fig:wcPanel}C and their respective values of $\bar{K}_{ij}$ in Fig.\ref{fig:wcPanel}B). 

The slow manifold $\mathcal{S}$ in panel D illustrates the different geometries underlying the approximation to the limit cycle $\Gamma$. For instance, when the system is slow in one variable as in the $RT$ model, the slow manifold has a cylinder-like shape in the direction of the slow variable $r$ (see Fig.~\ref{fig:rtPanel}D and Fig.~\ref{fig:hhPanel}D for negative values of the voltage $V$). The $QIF$ model has also a slow-fast dynamics (see the characteristic exponents $\lambda_1$ and $\lambda_2$ in Table~\ref{table:modelPrms}), but since the slow manifold is not in the direction of any coordinate axis, its cylindrical shape should became visible after a linear change of variables (Fig.~\ref{fig:qfPanel}D). Even if for completeness we also show the slow manifold for the $WC_{syn}$ model in Fig.~\ref{fig:wcPanel}, we stress that this manifold has not a special significance in this case since the characteristic exponents $\lambda_1$ and $\lambda_2$ are of the same size (see Table~\ref{table:modelPrms}).

Furthermore, the slow-fast dynamics is also reflected in the isochron distribution. In general, the isochrons are not homogeneously distributed in phase space for equidistant values of $\theta$. Indeed, they accumulate in the regions where the dynamics is slow, and they appear separated when the dynamics is fast. Notice that in Figs.~\ref{fig:rtPanel}C, \ref{fig:hhPanel}C and \ref{fig:qfPanel}C, one can see that isochrons for phases between 0 and 0.2 get distributed in approximately half of the limit cycle (warm coloured isochrons), while to cover the rest of the limit cycle we need the isochrons of phases between 0.2 and 1 (cool coloured isochrons). Consequently, if we consider equidistant values of $\theta$, we would observe that the isochrons accumulate in a part of the limit cycle. This can be appreciated in Panel D of Figs.~\ref{fig:rtPanel}, \ref{fig:hhPanel} and \ref{fig:qfPanel}, where the leaves $\mathcal{S}_{\theta}$, corresponding to the projection of the isochrons $\iloc$ onto the slow manifold, are equidistant in $\theta$. This phenomenon is different in Fig.~\ref{fig:wcPanel}C, where the isochrons are better distributed for equidistant values of $\theta$ (see also Fig.~\ref{fig:wcPanel}D). This is a consequence of the non slow-fast nature of the model $WC_{syn}$.

%\textcolor{red}{les isocrones bastant planes en els models de neurona tipus Slow fast, comparar amb el 2D de llavehuguet}

%Moreover, the shape of the isochrons is affected by ?? \textcolor{red}{les isocrones per als colors propers al roig, creixen en plans on el V es almost constant, mentre que en phases properes al blau ($\theta = 0.5$) es el contrari.}  This is refected in the iPRCs. Typically for phases near zero, the iPRCs is almost zero, whereas the maximum is achieved for phases in between 0.5 and 0.75. See Figure~\ref{fig:}

%\textcolor{red}{distribucio no homogenia de les isocrones i la relacio amb les PRCs}
%Notice that in Fig3E several  isochrons are parallel to the voltage directions Fig 7A right

The geometry of the isochrons is linked to the iPRF and iARFs shown in Figs.~\ref{fig:rtPanelPrf}-\ref{fig:qfPanelPrf}. For instance, when the isochrons get very close to each other, $\nabla \Theta$ increases dramatically because small perturbations can cause the trajectories to jump to an isochron with a very different phase. Indeed, in a neighbourhood of the unstable fixed point (a phaseless set for these models) the isochrons pack (see Figs.~\ref{fig:rtPanel}-\ref{fig:qfPanel}E) and the amplitude variable $\sigma$ tends to infinity, thus a perturbation acting on points therein generates large shifts both in phase and in amplitude (results not shown). Alongside, for systems with isochrons that are parallel to the coordinate plane containing the $V$-axis (as it happens in the $HH$ model, see Fig.~\ref{fig:hhPanel}E), the changes in phase due to perturbations in the direction of $V$ are going to be zero. In Figs.\ref{fig:rtPanelPrf}-\ref{fig:qfPanelPrf}A (right) we plot $\nabla \Theta$ along the isochron. Clearly, we see that $\nabla \Theta$ is close to 0 in Fig.~\ref{fig:hhPanelPrf}A (right) which corresponds to Fig.~\ref{fig:hhPanel}E. As we mentioned, one can also appreciate the increment of $\nabla \Theta$ on the isochrons when they accumulate near the phaseless set but, as explained before, we just show a part of the isochron before $\nabla \Theta$ blows up. Notice though that the effect of perturbations acting on points lying on the same isochron can be very different depending whether the point is close or far from the limit cycle.

In Figs.~\ref{fig:rtPanelPrf}-\ref{fig:qfPanelPrf} the magnitude of $\nabla \Sigma$ needs to be put in context. Indeed, since the amplitude variables $\sigma_1$ and $\sigma_2$ are scaled by the factors $b_1$ and $b_2$ one needs to control how the size of $\sigma$ translates to Euclidean distance to the limit cycle. In Section~\ref{sec:section5} we discuss this issue in more detail for a concrete example.  

\newpage

\section{Study of Perturbations using the Phase-Amplitude Variables}\label{sec:section5}

%study of weakly perturbed nonlinear oscillators

The phase reduction is a powerful tool to reduce high dimensional dynamics to a single equation and understand the dynamics that emerges when nonlinear oscillators are weakly perturbed. It has been extensively used to study weakly coupled oscillators and synchronization properties \cite{hoppensteadt2012}, specially in the neuroscience context \cite{ErmentroutTerman2010}, but it has several limitations. Indeed, the phase reduction assumes that the dynamics tkes place on the limit cycle or very close to it. 
However, if the attraction to the limit cycle is weak or the time between perturbations is short compared to the strength of the perturbation, this strong assumption is not valid anymore and the phase reduction becomes unreliable. 

%However, factors like the rate of convergence to the oscillator, strong forcing or high stimulation frequency may invalidate the above assumption and raise the question of how is the phase variation away from an attractor.

In this Section we show that the phase-amplitude approach based on the parameterization method can account for an accurate description of the phase variation away from the limit cycle $\Gamma$. The phase-amplitude description does not reduce the dimension, but allows us to identify the transversal directions that are relevant for the dynamics and suggests possible dynamical reductions. %More precisely, we will see those associated to Floquet multipliers that are close to 0, can be neglected. slow manifold.

%In some cases, only a fet amplitude coordinates are necessary. Indeed, those associated to Floquet multipliers that are close to 0, can be neglected. We will show how this phase-amplitude coordinates can capture behaviors that are not captured by the phase reduction.

In general, we consider perturbations of our original system $X$ in \eqref{eq:mathDef_1} such that the perturbed system is of the form
\begin{align}\label{eq:sispert}
\dot{x} = X_A(x, t) = X(x) + A p(t),
\end{align}
where $p(t)$ is a vector-valued function on $\mathbb{R}^d$ depending on time.
The dynamics of the perturbed system in phase-amplitude variables $(\theta,\sigma)=(\theta,\sigma_1,\ldots,\sigma_{d-1}) \in \mathbb{T} \times \mathbb{R}^{d-1}$ is given by 
\begin{equation}\label{eq:systa}
\begin{array}{rcl}
\dot{\theta} &= &\frac{1}{T} + A \nabla \Theta (K(\theta,\sigma))\cdot p(t), \\
%\dot{\sigma} &=& \Lambda \sigma + A \nabla \Sigma (K(\theta,\sigma))\cdot p(t), 
\dot{\sigma}_i &=& \lambda_i \sigma_i + A \nabla \Sigma_i (K(\theta,\sigma))\cdot p(t), \quad i=1,\ldots,d-1,
%\dot{\sigma}_2 &=& \lambda_2 \sigma_2 + A \nabla \Sigma_2 (K(\theta,\sigma))\cdot p(t), 
\end{array}
\end{equation}
where $\nabla \Theta (K(\theta,\sigma))$ and $\nabla \Sigma_i (K(\theta,\sigma))$, for $i=1,...,d-1$, are the iPRF and iARFs, respectively, defined in \eqref{eq:arcs_2}.

\begin{remark}
We want to stress that the distance to the limit cycle is not a limitation for our method. Indeed, $\nabla \Theta$ and $\nabla \Sigma$ are known for any point in $\Omega$. Recall that they are known locally as Taylor expansions around the limit cycle up to high order (order $L$) and can be globalized using equations \eqref{eq:prcs_4b} (see Section~\ref{sec:phaseAmpFun}).
Thus, the system above is exact.
\end{remark}

\begin{figure}
	\centering
	{\includegraphics[width=160mm]{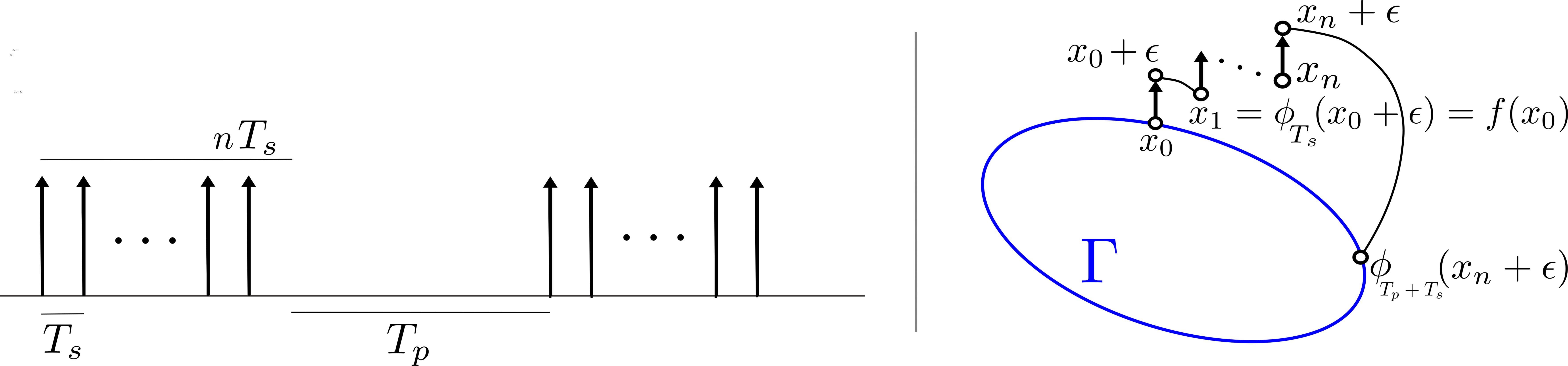}}
	\caption{The stimulus $p(t)$ in \eqref{eq:stimulus} consists of a train of $n$ pulses of size $\epsilon$ separated by a time $T_s$, which is repeated periodically and separated by a time $T_p$ (left). Scheme of the effects of the perturbation on the trajectory. Notice that the dynamics does not occur on the limit cycle (right).}\label{fig:perturbation3D}
\end{figure}

For our study, we will use a perturbation consisting of a train of $n$ pulses of size $\varepsilon$ separated by a time $T_s$. Each train of pulses is repeated periodically with period $T_p$ (see Fig.\ref{fig:perturbation3D} left), thus, generating a periodic perturbation $p(t)$ of period $T_{total}=nT_s+T_p$. Mathematically, 
\begin{equation}\label{eq:stimulus}
p(t)=\varepsilon \vec{v} \sum_{k=0}^{n} \delta(t-kT_s), \quad t \in [0, T_{total}),  
\end{equation}
where $\vec{v}$ is a $d$-dimensional vector representing the direction of the pulsatile perturbation. 

Notice that although this perturbation consists of small kicks (size $\varepsilon$), it manages to displace the trajectories away from the limit cycle $\Gamma$ if the interpulse interval $T_s$ is short and $n$ is large (see Fig.\ref{fig:perturbation3D} right). Indeed, when $T_s \rightarrow 0$, the perturbation becomes close to a delta perturbation of amplitude $n\cdot\epsilon$. 
%In other words, this perturbation is a way of studying large pulsatile perturbations by means of the iPRF and iARFs.

Since the perturbation is periodic, the dynamics can be described by the stroboscopic map, given by the flow of \eqref{eq:sispert} at time $T_{total}$ starting at $t=0$, which has 
the following expression 
\begin{equation}\label{eq:gemmaMap3D}
F(x):=\phi_{\substack{T_p}} \circ \underbrace{f \circ \cdots \circ f}_{n} (x),
\end{equation}
where $f(x):=\phi_{\substack{T_s}}(x+\varepsilon \vec v)$ and $\phi_t$ is the flow of the unperturbed vector field $X$. Notice that, knowing the position of the trajectory at a given kick $n$, the map $f$ provides the position of the trajectory at the time immediately preceding the next kick $n+1$.

Since the perturbation is pulsatile, the map $f$ can be obtained analytically in terms of the phase-amplitude variables and has the following expression:
%For pulses of weak amplitude, one can use the  and iARFs to write the map  
%
\begin{equation}\label{eq:fullMap3Dnoaprox}
\begin{aligned}
\bar \theta & = \Theta(K(\theta, \sigma) + \varepsilon \vec{v})  + \frac{T_s}{T}, \quad \quad  \text{(mod 1)} \quad k \in \mathbb{N} \\
\bar \sigma_i & = \left ( \Sigma_i(K(\theta, \sigma)+\varepsilon \vec{v}) \right )e^{\lambda_i T_s}, \quad \quad i = 1,\ldots, d-1. 
\end{aligned}
\end{equation}

Using that the pulses are of size $|\varepsilon| \ll 1$, the map above can be approximated using the iPRF and iARFs and has the following approximated expression $\bar f$:
%For pulses of weak amplitude, one can use the iPRFs and iARFs to write the map  
%
\begin{equation}\label{eq:fullMap3D}
\begin{aligned}
\bar \theta & = \theta + \varepsilon \nabla\Theta(K(\theta, \sigma))\cdot \vec{v} + \frac{T_s}{T}, \quad \quad  \text{(mod 1)} \quad k \in \mathbb{N} \\
\bar \sigma_i & = \big(\sigma_i + \varepsilon \nabla\Sigma_i(K(\theta, \sigma))\cdot \vec{v}\big)e^{\lambda_i T_s}, \quad \quad i = 1,\ldots, d-1. 
\end{aligned}
\end{equation}

Moreover, recall that the map $\phi_{T_p}$ is described exactly in phase-amplitude variables by the expression 
\eqref{eq:aboveEq}. Thus, we have a semi-analytical expression for the map $F$ in phase-amplitude variables.

\begin{remark}
Expression \eqref{eq:fullMap3D} for the map $f$ is not exact. Indeed, it uses that if $\varepsilon$ is small then the iPRF $\nabla \Theta$ and the iARFs $\nabla \Sigma_i$ provide a good approximation in first order of the PRF and ARFs, respectively (see Eq.~\eqref{eq:arcs_2}). For higher order approximations of the PRC see \cite{suvak2010quadratic, takeshita2010higher}. 
%For an exact expression of the map $f$ we need to integrate numerically system \eqref{eq:systa}, which in this case (although the iPRF and iARF also appear in them and might lead to a confusion) it is exact.
\end{remark}

As we already mentioned, the system with one phase and $d-1$ amplitude variables is not reduced, as the number of variables is the same as in the original system. However, we will use it to show that our computation in the new coordinates can capture the dynamics with high accuracy. 

%In high-dimensional systems it is often the case that some characteristic exponnents $\lambda_i$ are close to zero, and therefore, the corresponding coordinate $\sigma_i$ can be assumed to decay rapidly and assumed 0 and its corresponding equation removed in system \eqref{}, resulting in a reduction in the dimension.

\subsection{Dynamical Reductions}

In some cases, the particularities of the system or the perturbation allow for a reduction of the full phase-amplitude system. The \emph{phase reduction} is the most extended reduction. It assumes that perturbations are weak enough so that the trajectories are not displaced far away from the limit cycle. The study of the dynamics is then reduced to control the phase on the limit cycle. Thus, the transversal directions $\sigma_i$ are assumed to be 0 and the map $\bar f|_{\Gamma}$ is just a stroboscopic map of a circle to itself of the form,
\begin{equation}\label{eq:phaseMap}
\bar \theta = \theta + \varepsilon \nabla\Theta(K(\theta, 0))\cdot \vec{v} + \frac{T_s}{T}, \quad \quad  \text{(mod 1)}, 
\end{equation}
where $\nabla\Theta(K(\theta, 0))$ is the classical iPRC (see Remark~\ref{rm:prcsRM}).

Nevertheless, when studying large perturbations that displace the trajectories far from the limit cycle there exist more suitable reductions that can account for the dynamics. Since in high dimensional systems ($d > 2$), there exist more than one characteristic exponents, it is often the case that some of them are very small, as it happens close to a bifurcation. Recall that the amplitude variables $\sigma_i$ decay to zero at a rate that depends on its associated Floquet exponent $\lambda_i$. Therefore, those associated to Floquet exponents with large moduli can be assumed to be zero and its equation removed from \eqref{eq:fullMap3D}, thus reducing the dimension of the system. In the particular case that all variables $\sigma_i$ are assumed to be 0, except the one corresponding to the smallest (in modulus) Floquet exponent that we assume to be $\lambda_{d-1}$ without loss of generality, we say that the dynamics is reduced to the slow manifold $\mathcal{S}$ defined in \eqref{eq:slowManifold}. In this case, the map $\bar f|_{\mathcal{S}}$ writes as:
\begin{equation}\label{eq:variables_prc21}
\begin{aligned}
\bar \theta & = \theta + \varepsilon \nabla\Theta(K(\theta, 0, \ldots, 0,\sigma_{d-1}))\cdot \vec{v} + \frac{T_s}{T}, \quad \quad  \text{(mod 1)}  \\
\bar \sigma_{d-1} & = \big(\sigma_{d-1} + \varepsilon \nabla\Sigma_{d-1}(K(\theta, 0, \ldots, 0, \sigma_{d-1}))\cdot \vec{v}\big)e^{\lambda_{d-1} T_s}.
\end{aligned}
\end{equation} 
The dynamics, then, reduces to two equations and, furthermore, the number of monomials in $\bar K$ drastically reduces. We refer to this reduction as the \emph{slow manifold reduction}. 

\subsection{Examples}

In this Section we illustrate our methodology and compare the different dynamical reductions discussed in the previous Section. To do so, we choose the $RT$ model in \eqref{eq:rtEDOs}, as in \cite{wilson2018greater}. Recall that it is a $3$-dimensional model with an attracting limit cycle that has two associated Floquet exponents: $\lambda_2=-0.022$ and $\lambda_1=-0.368$, which is approximately 16 times greater (in modulus) than $\lambda_2$ (see Table~\ref{table:modelPrms}). Therefore, it is reasonable to explore also the slow manifold reduction for this system.

We apply the perturbation \eqref{eq:gemmaMap3D} (see Fig.~\ref{fig:perturbation3D}) with $n=100$ pulses of amplitude $\varepsilon=-0.1$, separated by a time interval $T_s=0.001$. The interval between input trains will be set to $T_p=8.394 \approx T$, where $T$ is the period of the unperturbed limit cycle. The perturbation will be in the voltage direction, that is, $\vec{v}=(1,0,0)$. Notice that the sign of $\varepsilon$ is negative, therefore the perturbation is inhibitory. We apply inhibitory perturbations because the changes in $\sigma_2$ are greater for inhibitory perturbations than excitatory, therefore they displace the trajectory further away from the limit cycle as observed in \cite{wilson2018greater} (results not shown). 
%By using this perturbation we expect that the dynamics occurs far from the limit cycle. 

\begin{figure}
	\centering
	{\includegraphics[width=160mm]{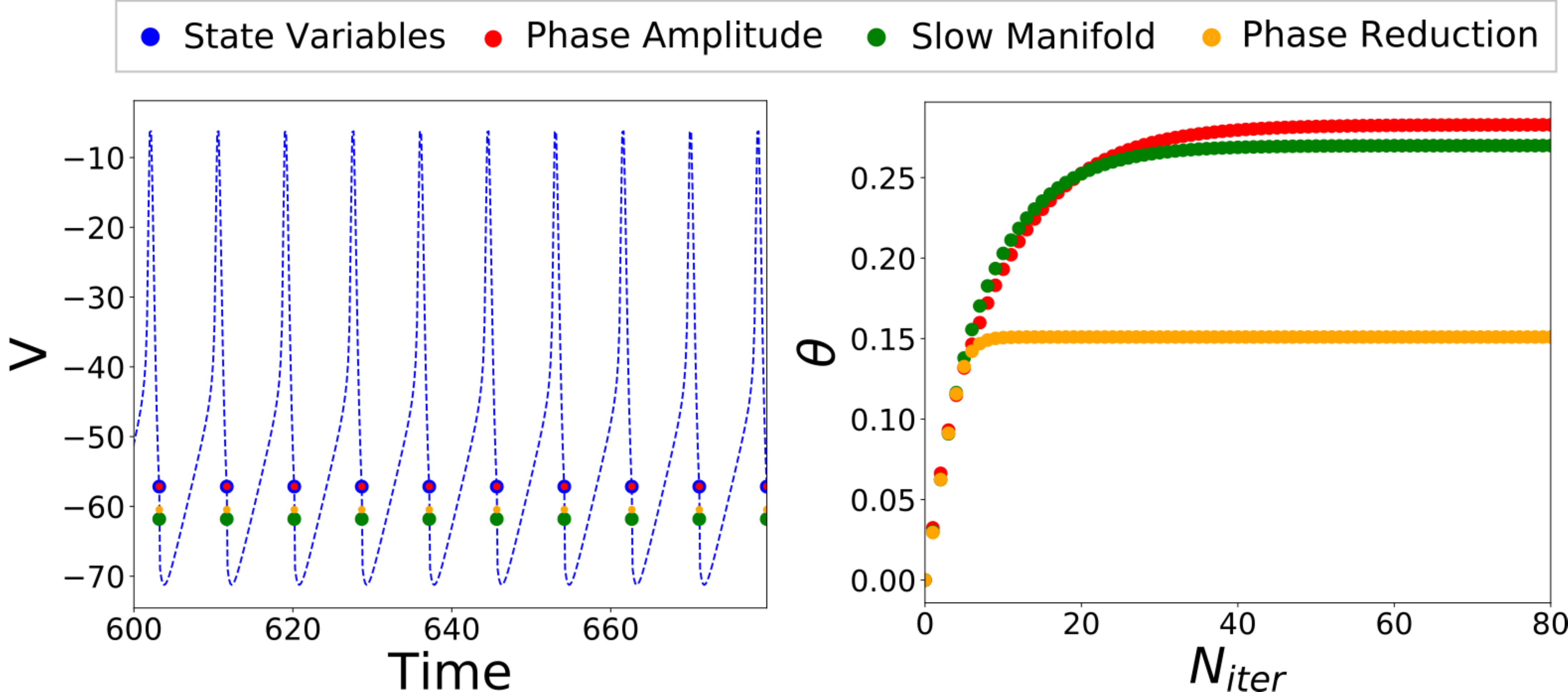}}
	\caption{(Left) Voltage coordinate of the trajectory starting at $x = K(0,0,0)$ for the perturbed system \eqref{eq:sispert} (dashed blue curve) and voltage position for the last ten iterates of the four different stroboscopic maps. (Right) Phase coordinate for the first 80 iterates of the stroboscopic maps given by the full phase-amplitude map $\bar{F}$ (red), the slow manifold reduction map $\bar{F}|_{\mathcal{S}}$ (green) and the phase reduction map $\bar{F}|_{\Gamma}$ (yellow). }\label{fig:panellPertRT}
\end{figure}

We compute the iterates using four different stroboscopic maps. The first one is the 3-dimensional map $F$ in \eqref{eq:gemmaMap3D} obtained by direct integration of the vector field \eqref{eq:sispert} in the original coordinates $(V,h,n) \in \mathbb{R}^3$. We refer to this map as the \emph{state variables map}. This map provides the exact description of the dynamics. The second one is an approximation of the map $F$ given in the phase-amplitude coordinates $(\theta,\sigma) \in \mathbb{T} \times \mathbb{R}^2$. We use expression of $\bar f$ in \eqref{eq:fullMap3D} to approximate $f$ and Eq.~\eqref{eq:aboveEq} for $\phi_{T_p}$. We denote it by $\bar{F}$ and refer to it as the \emph{Phase-Amplitude map}. Notice that this map is still 3-dimensional. The other two maps are dimensional reductions of the {Phase-Amplitude map}. One is obtained by setting the variable $\sigma_1=0$ (see Eq.~\eqref{eq:variables_prc21}), and it is 2-dimensional. We refer to this one as the \emph{slow manifold reduction} and denote it by $\bar{F}|_{\mathcal{S}}$. The other one is the classical 1-dimensional \emph{phase map} $\bar{F}|_{\Gamma}$ and is obtained by setting all amplitude variables to zero (see Eq.~\eqref{eq:phaseMap}) and keeping only the phase variable.

Next we discuss the performance of the above mentioned maps. For all of them we compute $N=80$ iterates of the map starting on the point with 0-phase on the limit cycle. Notice that this point has $\sigma_1=\sigma_2=0$ amplitudes, therefore, it is well-defined for all maps and reductions. We observe that, for all maps, the iterates tend to a fixed point, but it is different for each map. Fig.~\ref{fig:panellPertRT} shows the phase coordinate for 80 iterates of the maps $\bar F$, $\bar F|_{\mathcal{S}}$ and $\bar F|_{\Gamma}$ (right) and the evolution of the voltage variable along the corresponding trajectory with the values of the voltage at the last 10 iterates of the maps $F$, $\bar F$, $\bar F|_{\mathcal{S}}$ and $\bar F|_{\Gamma}$ (left).  In Table \ref{table:puntFix} we present the coordinates of the fixed point obtained for the different maps. We provide the fixed points in state variables (obtained evaluating the parameterization $K$) and phase-amplitude variables, when applicable. Notice that only the phase-amplitude variables provide an accurate description of the dynamics of the system. Nevertheless, the slow manifold reduction can capture the value of the phase variable. 

\begin{figure}
\centering
\includegraphics[width=150mm]{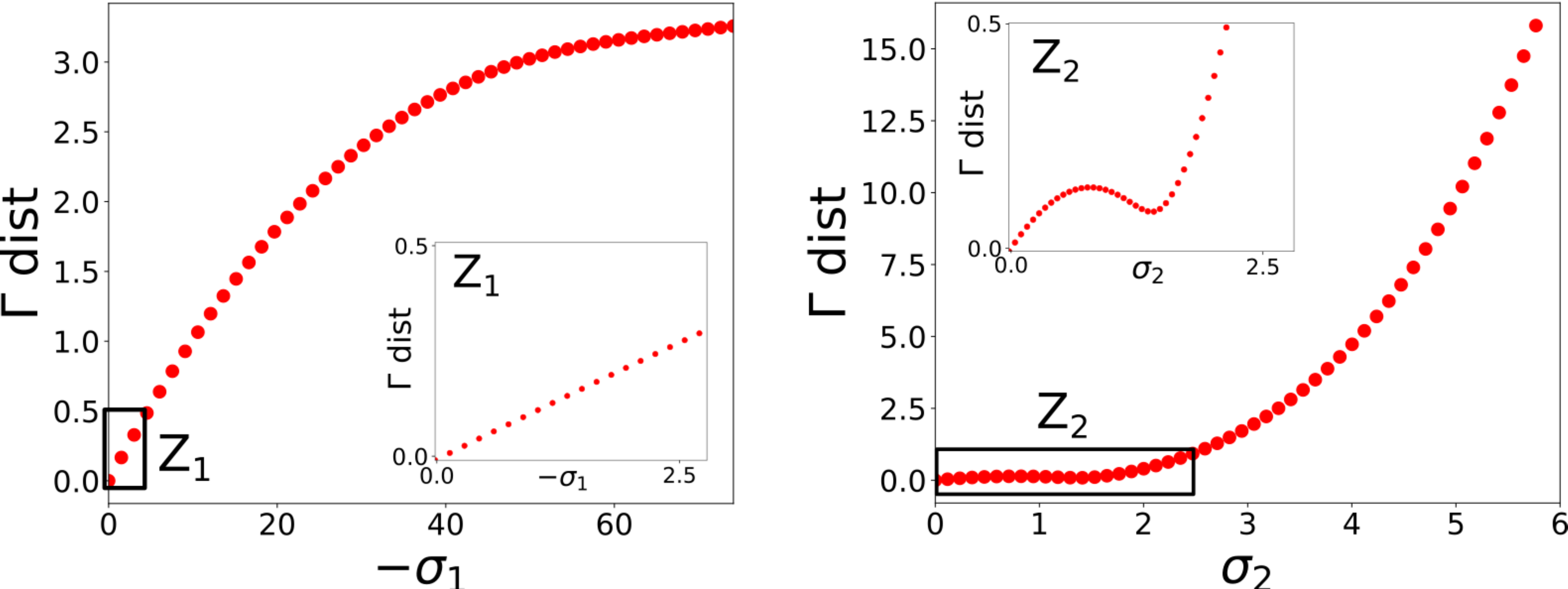}
\caption{For the \emph{RT} model \eqref{eq:rtEDOs} we show the Euclidean distance to the limit cycle as a function of the value of the amplitude coordinate $\sigma_1$ (left) and $\sigma_2$ (right). More precisely, for $\theta=0.283$, we show
$\Vert K(\theta,\sigma_1,0)-K(\theta,0,0)\Vert$ (left) and $\Vert K(\theta,0,\sigma_2)-K(\theta,0,0)\Vert$ (right) }\label{fig:proves}
\end{figure}

By looking at the $(\theta,\sigma)$ coordinates of the fixed points it might seem that the assumption that some amplitude variable is zero is very strong. However, we want to recall that the amplitude variables $\sigma_1$ and $\sigma_2$ are scaled by the constants $b_1$ and $b_2$, thus, their absolute values do not provide enough information on the distance to the limit cycle. Indeed, in the computation performed for the $RT$ model, we can see that small variations in $\sigma_1$ translate to small variations in the distance to the limit cycle (see Fig.~\ref{fig:proves} left). A different behaviour occurs for $\sigma_2$, where small variations translate to large variations in the distance to the limit cycle (see Fig.~\ref{fig:proves} right). Thus, the assumption that $\sigma_1=0$ in the slow manifold reduction is less strong that the  assumption that $\sigma_2=0$ in the phase reduction (see Table~\ref{table:puntFix} and Fig~\ref{fig:proves}).

%Figure~\ref{fig:proves} shows the Euclidean distance to the limit cycle as a function of the values of $\sigma_1$ (left) and $\sigma_2$ (right).

\begin{figure}
	\centering
	{\includegraphics[width=160mm]{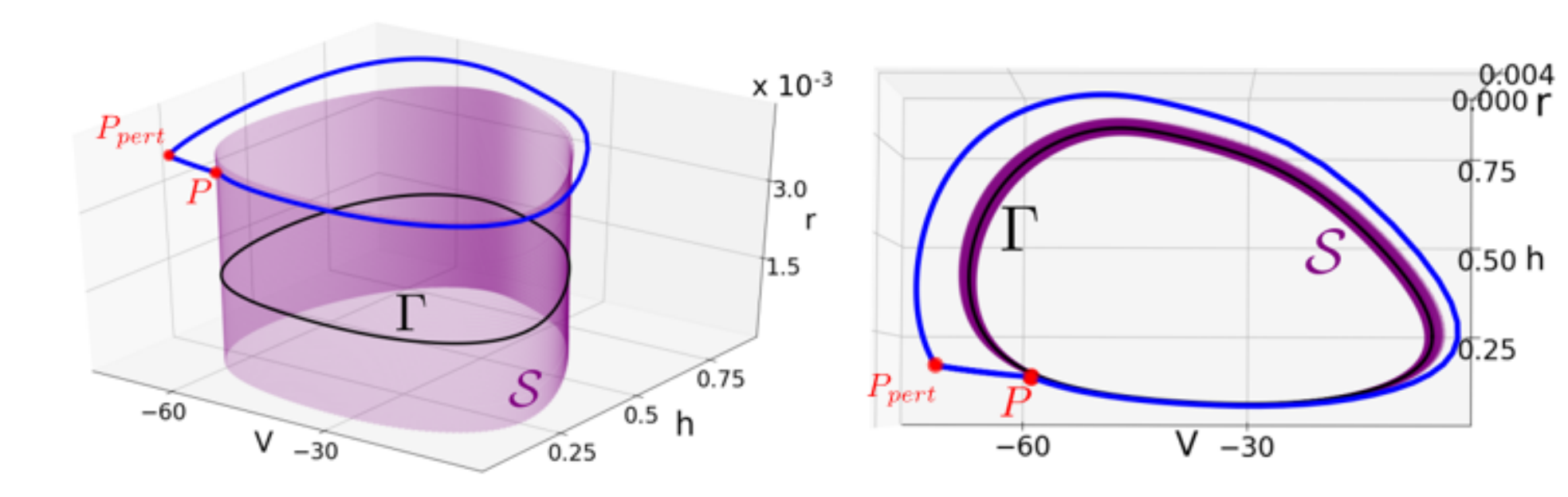}}
	\caption{Trajectory of the perturbed system \eqref{eq:sispert} (blue curve) corresponding to the fixed point $P$ of the stroboscopic map \eqref{eq:gemmaMap3D} (i.e $F(P)=P$), and the points $P$ and $P_{pert}=f^n(P)$ together with the slow manifold $\mathcal{S}$ (purple) and the limit cycle $\Gamma$ (black curve) from two different perspectives. Notice that the fixed point $P$ lies on the slow manifold and the perturbation displaces trajectories away from the slow manifold ($P_{pert}$). Once the perturbation vanishes the trajectory approaches again the slow manifold $\mathcal{S}$. }\label{fig:panellSlowMan}
\end{figure}

Results in Table \ref{table:puntFix} can be interpreted by means of Fig.~\ref{fig:panellSlowMan}, which shows the limit cycle $\Gamma$, the slow manifold $\mathcal{S}$ and the trajectory for a time $T_{total}$ of the perturbed system starting at the fixed point $P$ of the map $F$. We also show the point $P_{pert}$, which provides the position of the trajectory after the perturbation has been applied, that is $P_{pert}=f^n(P)$.  Observe that the trajectory $\phi_{t}(P)$ for $t \in [0,T_{total})$ is not on the limit cycle $\Gamma$, which explains the poor description of the dynamics when using the phase map \eqref{eq:phaseMap}. Notice that the perturbation displaces the trajectory away from the slow manifold (see position of $P_{pert}$), but the trajectory relaxes back to the slow manifold $\mathcal{S}$ when the perturbation is removed (see position of point $P$). 
In this example, the phase-amplitude map is the most suitable to describe properly the dynamics. We would like to highlight the good performance of the phase-amplitude map \eqref{eq:fullMap3D} (see Table~3). However, if we consider a larger value of $T_s$ so that the trajectories relax back to the slow manifold but not to the limit cycle after every kick, we expect that the slow manifold map will provide as well an accurate description of the dynamics. We stress that for the phase map to provide an accurate description of the dynamics, this time should be much longer. Another possibility to improve the phase reduction could be to use higher order Taylor expansions of the PRC besides the infinitesimal PRC. We leave this explorations for future work.
%For this reason, the performance of the slow manifold map is not accurate enough, although the point $P$ is very close to the slow manifold. 
%Indeed, the interpulse interval $T_s$ is still too short to allow the trajectory to relax back to the slow manifold and therefore, the iterates computed assuming that the trajectory is on the slow manifold $\mathcal{S}$ using the map $\bar{f}|_{\mathcal{S}}$ have an error which affects predictions. 

\begin{table}[H]
	\renewcommand{\arraystretch}{1.5}
	\begin{center}
		\begin{tabular}{|c|c|c|c|c|c|c|}
			\hline
			Map     % & Symbol & 
			& $V$ & $h$ & $r$ & $\theta$ & $\sigma_1$ & $\sigma_2$ \\
			\hline
			State Variables $F$        %& $P = F(P)$ 
			& $-57.16$ & $0.135$ & $0.00383$ & ---- & ---- & ----\\
			\hline
			Phase-Amplitude $\bar{F}$ %&  $P_{\Omega} = \bar{F}(P_{\Omega})$   
			& $-57.13$ & $0.135$ & $0.00377$ & $0.283$ & $-2.37$ & $4.98$ \\
			\hline
			Slow Manifold Reduction $\bar{F}|_\mathcal{S}$ %&  $P_{\mathcal{S}} = \bar{F}|_\mathcal{S}(P_{\mathcal{S}})$  
			& $-61.81$ & $0.197$ & $0.00314$ & $0.269$ & $0$ & $3.439$\\
			\hline
			Phase Reduction   $\bar{F}|_\Gamma$ % &  $P_{\Gamma} = \bar{F}|_\Gamma(P_{\Gamma})$  
			& $-60.458$ & $0.175$ & $0.0017$ & $0.15$ & $0$ & $0$ \\
			\hline					
		\end{tabular}
	\end{center}
	\renewcommand{\arraystretch}{1}
	\caption{Values of the $V$, $h$ and $r$ coordinates of the fixed point $P$ of the stroboscopic map $F$ in the state variables, the fixed point $P_{\Omega}$ of the full phase-amplitude map $\bar{F}$, the fixed point $P_{\mathcal{S}}$ of the slow manifold reduction $ \bar{F}|_\mathcal{S}$ and the fixed point $P_{\Gamma}$ of the phase reduction map $ \bar{F}|_\Gamma$. We also show the values for the coordinates $\theta$, $\sigma_1$ and $\sigma_2$ of the corresponding fixed point.}
	\label{table:puntFix}
\end{table}

\section{Conclusions and discussion}\label{sec:section6}

In this paper we extend the applications of the parameterization method to study the phase-amplitude description for $d$-dimensional non-linear oscillators with $d \geq 2$. More precisely, we consider systems having a hyperbolic attracting limit cycle $\Gamma$ and we provide a computational method to obtain a parameterization $K$ of its attracting invariant manifold (which in fact is its whole basin of attraction) in terms of the phase-amplitude variables. The use of this method has several advantages. On one hand, it automatically provides a geometrical portrait of the oscillator through the computation of the isochrons and the isostables, including the slow submanifold. On the other hand, it provides the iPRF and iARFs, which allow us to accurately track the effects of a given perturbation $p(t)$ beyond the phase reduction. 

 The theoretical and numerical methods presented herein extend previous results for planar systems \cite{guillamon2009computational, huguet2013computation} to higher dimensional systems ($d \geq 2$), which are based on the classical parameterization method. In this paper, we use efficient numerical tools to obtain a semi-analytical approximation of the parameterization $K$ using Fourier-Taylor expansions around the limit cycle $\Gamma$. Whereas the lower order terms in the Taylor expansion (degree 0 and 1) are straightforward to obtain using variational equations (see for instance \cite{brown2004phase}), the higher order terms (degree larger than 1) involve solving a homological equation \eqref{eq:homologicalEqs} for each term in the Taylor expansion. There are two main challenges to solve these equations. On one hand, to obtain the homological equations we need to compute the Taylor expansion of the composition of the vector field $X$ with the terms of the Taylor expansion of $K$ of lower order previously computed. Some works \cite{wilson2020} propose a symbolic manipulator to compute these derivatives, but for high orders and high dimensional systems this strategy is computationally very expensive. Alternatively, we propose to use automatic differentiation techniques as it is typical in the parameterization method \cite{haro2016}. These techniques avoid the use of symbolic differentiation, while provide the exact recurrent formulas for the terms of the Taylor expansion, which are faster and without numerical errors besides roundoff errors. On the other hand, one needs to find a periodic solution of each homological equation. In \cite{wilson2020} the authors design a Newton method to find a periodic orbit for each equation. In this paper, we adapt the methodology in \cite{castelli2015parameterization} (see also \cite{huguet2013computation}), based on Floquet theory, to reduce the problem to solve a diagonal linear system for the Fourier coefficients, which has an exact explicit solution. Of course, we can easily switch between Fourier and real space using the FFT algorithms. 
 The algorithms are efficient in the sense that if we discretize the functions using $N$ points, solving the homological equation requires $\mathcal{O}(N)$ storage and $\mathcal{O}(N \log N)$ operations in Fourier discretization.
 Moreover, this method takes advantage of the form of the equations to provide a general solution for all orders where only the computation of the terms $B_{\alpha,m-\alpha}$ in \eqref{eq:homologicalEqs} is required at each step (see Remark~\ref{rem:general}).  Thus, our method allows us to obtain high order analytical approximations of the map $K$ at a low computational cost. Indeed, as we indicate in Section~\ref{sec:section42}, for Taylor expansions of order $L=10$ and $N=2048$  Fourier coefficients, the numerical computation of $\bar{K}$ only takes a few seconds (around 10s) on a regular laptop.

% The computed $\bar K$ provides analytically the local isochrons.

Our approach also proposes a numerical strategy to globalize the local isochrons given by the approximate parameterization $\bar K$. There exist several approaches in the literature for computing global (un)stable manifolds of a vector field, see \cite{Krauskopfetal05, simo1990analytical}. In \cite{guillamon2009computational} the authors adapted the method in \cite{simo1990analytical} which uses backwards integration to compute $1$-dimensional global isochrons. The procedure proposes an efficient way of choosing the points in the local domain to obtain a homogeneous coverage of the isochron when globalizing it. The higher dimensional case $d > 2$, addressed in this paper, is more difficult because there are several expanding directions and using rough backwards integration produces points that are concentrated along the direction of the largest multiplier. In this paper, we solve this drawback and we provide the details for the 3-dimensional case. Our procedure is based on \cite{simo1990analytical}, which computes first the invariant slow submanifold  $\mathcal{S}$ \eqref{eq:slowManifold}, corresponding to the direction of the smallest (in modulus) multiplier. It is known that the trajectories asymptotically approach the limit cycle along this manifold. Furthermore, we use this manifold to globalize the isochron along the isostables (see Fig.~\ref{fig:globalMethod3d}C). Finally, we use the parameterization $K$ to compute the iPRF and the iARFs analytically for points in the local domain and we globalize these functions to the whole basin of attraction by using a modified adjoint equation \eqref{eq:prcs_4b}. We illustrate the methodology for different models in neuroscience. 

Finally, in Section \ref{sec:section5} we use the computed parameterization $K$ to extend the study of phase dynamics for perturbed systems beyond the limit cycle. We show first that the full phase-amplitude description captures the perturbed dynamics accurately. Moreover, we have explored the scope of applicability of two useful reductions, namely the slow manifold reduction and the phase reduction. While the limitations of the phase reduction are already well-known (see for instance \cite{castejon2013phase, AshComNic2016, wilson2018greater}) we have been able to identify the limitations of the slow manifold reduction for a particular example. More precisely, we observe that even for systems with two Floquet exponents with a large difference in size and one of them close to zero, 
%that is, for which it is easy to identify a slow manifold, this property is not enough to assume that 
the slow manifold reduction might not be accurate enough. Indeed, if the interpulse interval is too short, the slow manifold reduction might not be sufficient and the method would require to consider the full phase-amplitude system. When this period increases, we expect that the performance of the slow manifold is guaranteed.

We acknowledge that we assumed that the Floquet multipliers of the limit cycle are all different and real. We highlight that we can easily tackle the extension to the case with complex conjugate multipliers, just adapting the techniques in \cite{castelli2015parameterization}. Of special interest is the case where the Floquet exponents are resonant (see Remark~\ref{rem:nres0} and Remark~\ref{rem:nonres}). In this case, it is not possible to conjugate the dynamics to a linear vector field and one needs to consider higher order terms in the dynamics of the amplitude variables $\sigma$ (see \cite{cabre2003parameterization} for a thorough discussion and \cite{James16} for an application to the case of the stable manifold of an equilibrium point). A detailed study of this case is an interesting topic for future research.  

The study of the oscillatory dynamics beyond the phase reduction has become a topic of growing interest in the last decade (see \cite{AshComNic2016,ermentroutetal19}). We briefly discuss how our method based on the parameterization method, compares with other approaches in the literature. In \cite{wilson2018greater} the authors provide a method to compute the iPRF $\nabla \Theta$ and iARFs $\nabla \Sigma_i$, $i=1,\ldots,d-1$ in a neighbourhood of the limit cycle by means of a linear order expansion of these functions in the amplitude variable, using the generalized adjoint equations. In a recent paper \cite{wilson2020}, the author extends the previous results to obtain a higher order expansion  of the iPRF and iARFs in the amplitude variables in a neighbourhood of the limit cycle. Moreover, he also writes the classical homological equations for computing the parameterization $K$ \cite{cabre2003parameterization}. We want to stress that the strategy to solve these equations is different from ours, as discussed above. 

Several groups use the spectral properties of the Koopman operator \cite{budivsic2012applied, mezic2005spectral} to provide an alternative  approach to the study of the phase-amplitude dynamics of a limit cycle. In particular, the Koopman operator has been used to compute global isochrons and isostables of a given limit cycle by forward integration \cite{mauroy2012use}, as well as PRFs and ARFs. This methodology has been used to compute the amplitude coordinate corresponding to the dominant Floquet multiplier to study  perturbed oscillatory dynamics \cite{mauroy2018global, shirasaka2017phase}. We emphasize that our method provides all transversal directions in the same computation using the same procedure.
%Thus, it can only capture the full dynamics for planar systems \cite{mauroy2018global}. 

%Thus, the methodology in this paper complements the approaches based on the Koopman operator specially for non-planar cases.

%Although our methodology is not the same, as it relies in the obtention of $K(\theta, \sigma)$, both approaches consider the computation of high order terms in $\sigma$ to increase the accuracy of the phase. This is why   

 %However, similarly to these methodologies, we also need to compute high order derivatives of the vector field at the limit cycle. At this point, we would like to stress how fundamental becomes the automatic differentiation for our purposes. The usage of automatic differentiation, far from being a challenging technique, is a wide used tool which permits to obtain high order derivatives up to n-order automatically. Thus, it considerably facilitates the obtention of these high order terms, as it avoids the needing of writing higher order derivatives by hand (or symbolically) and compute them. Once the high order terms of the derivatives are obtained, they easily expand the domain of $K(\theta, \sigma)$ through the simplified (via Floquet normal form) homological equation. Therefore, thanks to the automatic differentiation we are able to provide a simple algorithm to obtain $K(\theta, \sigma)$, and use it go beyond the phase-reduction by simply  evaluating a polynomial at $\sigma$ and a Fourier-series at $\theta$.

We acknowledge here that we only have applied our methods to neural oscillators, even if they are valid for any nonlinear oscillator. Neural oscillatory activity is widely observed at different levels of organization \cite{buzsaki2006rhythms}, so the phase reduction of neural oscillators has been profusely studied from the single neuron to the network level \cite{rinzel1998analysis, brown2004phase, izhikevich2007, AshComNic2016}. Moreover, the PRCs have been an important tool in neuroscience both from the theoretical and experimental perspective \cite{smeal2010phase, gutkin2005phase, achuthan2009phase, canavier2015phase, schultheiss2011phase}. However, other areas of biology, life sciences or control theory offer more examples of oscillatory activity \cite{strogatzbook, monga2019phase} to which this methodology can be applied.

Over the last decade, there has been an increasing need of understanding phase dynamics from biological data \cite{tass2007phase}. Indeed, the interest in therapies focusing on stimulating cellular tissue at particular phases of certain pathological rhythms is increasing since it has been reported to suppress them \cite{holt2016phasic, rosenblum2004delayed, azodi2015phase}. This, on its turn, has stimulated the study and computation of phase dynamics and PRCs in the stochastic domain \cite{schwabedal2013phase, thomas2014asymptotic, rosenblum2019numerical}. Since our methodology provides accurate descriptions and efficient computational techniques,  we consider its possible extensions in the stochastic domain as a very promising and interesting area of research.

In conclusion, in this paper we show how the parameterization method provides a clear and solid theoretical framework which
can be used in an efficient computable way to obtain both a geometrical and dynamical accurate description of the oscillatory dynamics. 
We hope that the application of the different techniques herein can be useful to gain insight and alternative perspectives to oscillatory dynamics.

\section*{Acknowledgments}
This work has been partially funded by the Spanish MINECO-FEDER Grant PGC2018-098676-B-100 (AEI/FEDER/UE), the Catalan Grant 2017SGR1049 (GH, AP, TS). GH acknowledges the RyC project RYC-2014-15866. TS is supported by the
Catalan Institution for research and advanced studies via an ICREA academia price 2018. AP
acknowledges the FPI Grant from project MINECO-FEDER-UE MTM2012-31714. AP acknowledges the long-term strategic development financing of the Institute of Computer Science (RVO:67985807) of the Czech Academy of Sciences.

\section*{Data Availability Statement} 

Data sharing is not applicable to this article as no new data were created or analyzed in this study.

\appendix
\setcounter{equation}{0}
\renewcommand{\theequation}{\thesection\arabic{equation}}
\section{Computation of the Floquet normal form}\label{sec:floquetComputation}

In this Appendix we provide the method to compute the matrices $\mathcal{Q}$ and $R$ that appear in formula \eqref{eq:floquetDifInv}.
By Floquet theory \cite{floquet1883equations}, the monodromy matrix satisfies
\begin{equation}
M = \Phi(T) = \mathcal{Q}(T)e^{TR} = \mathcal{Q}(0)e^{TR} = \Phi(0)e^{TR} = e^{TR},
\end{equation}
where we have used $\Phi(0) = Id$ and $\mathcal{Q}(T)$ is a T-periodic matrix.

Then, if there exists a matrix $C$ such that 
\begin{equation}
\Phi(T) = CDC^{-1} \quad \quad \quad \text{where} \quad D = diag(1, \mu_1, \mu_2),
\end{equation}
one can find the matrix $R$ in \eqref{eq:floqNf} as
\begin{equation}
R = \frac{1}{T} C  \begin{pmatrix}
0 & &   \\
& \ln(\mu_1) & \\ 
& & \ln(\mu_2)
\end{pmatrix} C^{-1} = CJ C^{-1},
\end{equation}
where $J$ is given in \eqref{eq:tereJ}, with $\lambda_0 = 0$, $\lambda_i = \ln(\mu_i)/T$ $i=1,2$. Therefore, the matrix $\mathcal{Q}$ in \eqref{eq:floqNf} is given by
\begin{equation}
\mathcal{Q}(t) = \Phi(t) e^{-tR} = \Phi(t) C \begin{pmatrix}
1 & &   \\
& e^{\frac{-t}{T}\ln(\mu_1)} & \\ 
& & e^{\frac{-t}{T}\ln(\mu_2)}
\end{pmatrix} C^{-1} = \Phi(t) C \begin{pmatrix}
1 & &   \\
& e^{-t \lambda_1} & \\ 
& & e^{-t \lambda_2}
\end{pmatrix} C^{-1}.
\end{equation}
We recall that $\mu_1$ and $\mu_2$ are the Floquet multipliers and $\lambda_1$ and $\lambda_2$ the Floquet exponents.

\section{Automatic Differentiation}\label{sec:autoDifAp}
\setcounter{equation}{0}

In this Section we aim to illustrate the Automatic Differentiation techniques. We recall that we use this technique in Section \ref{sec:sectionNum} to obtain the terms $B_{\alpha, m-\alpha}$ in \eqref{eq:homologicalEqs}, corresponding to the terms of degree $m$ of the Taylor expansion of $X(K(\theta, \sigma))$ around $\sigma = 0$ (see Eq. \ref{eq:Bexpand}). A possible approach to obtain the coefficients of the Taylor expansion is the to compute the appropriate derivatives. 
Writing and computing explicitly (using symbolic differentiation) the derivatives up to a given order is a costly task. Nevertheless, if the vector field $X$ is analytic, it can be written as a combination of algebraic operations (sum, product, etc.) and elementary transcendental functions (sin, cos, exp, log, power, etc.), and we can overcome this drawback by using automatic differentiation techniques. Automatic differentiation avoids the need of calculating by hand the $n$-th order derivatives by obtaining 
the terms in \eqref{eq:Bexpand} using simple recurrent relations. These relations are obtained by applying systematically the chain rule to each of the operations that compose the function. We stress that this procedure, which is faster than symbolic differentiation, provides the terms without numerical error. This appendix aims to show the basics of this technique for a Taylor expansion in two variables. For more details we refer the reader to \cite{griewank2008evaluating, haro2016, JorbaZ05}. 

Start by considering a function $f$ expressed in power series in the variables $\sigma_1, \sigma_2 \in \mathbb{R}$,
\begin{equation}\label{eq:equacioPerF}
f(\theta, \sigma) = \sum_{m=0}^{\infty} \sum_{\alpha=0}^{m} f_{\alpha, m-\alpha}(\theta) \sigma^{\alpha}_1 \sigma^{m-\alpha}_2,
\end{equation}
where $f_{\alpha, m-\alpha}(\theta)$ are periodic coefficients in $\theta$. 

Given an analytic function $\varphi$, our goal is to find the Taylor expansion of the composition $\varphi(f)$, that is,
\begin{equation}\label{eq:equacioPerPhi}
\varphi(f)(\theta, \sigma) = \sum_{m=0}^{\infty} \sum_{\alpha=0}^{m} [\varphi(f)]_{\alpha, m-\alpha}(\theta) \sigma^{\alpha}_1 \sigma^{m-\alpha}_2.
\end{equation}

Even if the coefficients of these Taylor expansions depend periodically on the variable $\theta$, this feature does not play any role in the procedure we are going to expose.
For this reason, from now on we do not write explicitly this dependence. Therefore, we will write $f(\sigma)$ instead of $f(\theta,\sigma)$ and  $f_{\alpha, m-\alpha}$ instead of $f_{\alpha, m-\alpha}(\theta)$.

We introduce the radial derivative, which is defined as
\begin{equation}
Rf(\sigma) = \nabla f(\sigma) \cdot \sigma = \sum_{i=1}^{d=2} \frac{\partial f}{\partial \sigma_i}(\sigma) \sigma_i,
\end{equation}
and has two useful properties. The first one, is that when applied to the Taylor expansion \eqref{eq:equacioPerF} it satisfies the following relationship
\begin{equation}\label{eq:eulerIdentity}
Rf(\sigma) = \sum_{m=0}^{\infty} m \sum_{\alpha=0}^{m} f_{\alpha, m-\alpha} \sigma^{\alpha}_1 \sigma^{m-\alpha}_2,
\end{equation}
which is known as Euler's identity. The second one comes up when applying the chain rule to the composition $\varphi(f(\sigma))$ in \eqref{eq:equacioPerPhi},
\begin{equation}\label{eq:chainRule}
R \varphi(f(\sigma)) = \varphi'(f(\sigma)) R \varphi(f(\sigma)).
\end{equation}

Next, we will show that for the case of $\varphi(x)$ being an elementary function, the combination of \eqref{eq:eulerIdentity} and \eqref{eq:chainRule}, allows us to compute the $n$-th order coefficient of \eqref{eq:equacioPerPhi} starting by the already known $\varphi(f_0)$ as initial seed. We illustrate the methodology for the example $\varphi(x) = \exp(x)$. Writing the Taylor expansion of the function $\exp(f)(\sigma)$ as
\begin{equation}\label{eq:equacioPerExp}
\exp(f)(\sigma) =  \sum_{m=0}^{\infty} \sum_{\alpha=0}^{m} e_{\alpha, m-\alpha} \sigma^{\alpha}_1 \sigma^{m-\alpha}_2,
\end{equation}
and since $\varphi'(x)= \varphi(x)$, substituting \eqref{eq:equacioPerExp} in \eqref{eq:chainRule} and using \eqref{eq:eulerIdentity} yields 
\begin{equation}\label{eq:bigEquation}
\sum_{m=0}^{\infty} m \sum_{\alpha=0}^{m} e_{\alpha, m-\alpha} \sigma^{\alpha}_1 \sigma^{m-\alpha}_2 = \left(\sum_{n=0}^{\infty} \sum_{\beta=0}^{n} e_{\beta, n-\beta} \sigma^{\beta}_1 \sigma^{n-\beta}_2 \right)\left(\sum_{l=0}^{\infty} l \sum_{\mu=0}^{l} f_{\mu, l-\mu} \sigma^{\mu}_1 \sigma^{l-\mu}_2 \right).
\end{equation}
Then, collecting the terms of the same order, we end up with the following recurrence equation
\begin{equation}\label{eq:recEquation}
e_{\alpha, m-\alpha} = \frac{1}{m} \sum_{k=0}^{\alpha-1} \sum_{l=0}^{m-\alpha} (m - k - l) e_{k,l} f_{\alpha-k, m-\alpha-l}, \quad \quad \quad \text{for} \quad m \geq 1, \alpha=0,\ldots,m.
\end{equation}
Notice that the computation of $e_{\alpha, m-\alpha}$ only depends on terms $e_{k,l}$, where $0 \leq k+l \leq m-1$.

Recurrence equations for a large variety of transcendental functions can be found in Table 2.1 in \cite{haro2016}.

\section{Power Expansions of the iPRFs $\nabla \Theta$ and iARFs $\nabla \Sigma$}\label{sec:expPotenciesGradient}
\setcounter{equation}{0}

In some cases, one might be interested in the explicit expression of the functions $\nabla \Theta$ and $\nabla \Sigma$ in power expansions in $\sigma$. Next, since formula \eqref{eq:usefulSpace} uses implicitly the fact that the functions $\nabla \Theta$ and $\nabla \Sigma_i$, $i=1,2$ are known at any order $L$, we show that \eqref{eq:usefulSpace} provides also the terms in $\sigma$. To do so, notice that using \eqref{eq:mjFourierTaylor} we have (recall that $K_0(\theta):=K_{00}(\theta)$, see Eq.~\eqref{eq:eqMzero}): 
\[ 
\begin{array}{rcl}
DK(\theta, \sigma) & = &
\sum_{m=0}^{\infty} \sum_{\alpha=0}^{m}
\big [ \begin{array}{c|c|c} 
K'_{\alpha, m-\alpha}(\theta) &
(\alpha+1) K_{\alpha+1, m-\alpha}(\theta)  &
(m- \alpha+1) K_{\alpha, m-\alpha+1}(\theta)  \\
\end{array}
\big ]
\sigma^{\alpha}_1 \sigma^{m-\alpha}_2 \\
\\
& = &
\big [ \begin{array}{c|c|c} 
K'_{0}(\theta) &
K_{10}(\theta)  &
K_{01} \\
\end{array} \big]
+
\big [ \begin{array}{c|c|c} 
K'_{10}(\theta) &
2K_{20}(\theta)  &
K_{11} (\theta)\\
\end{array} \big ]
\sigma_1 
\\
\\
&& +
\big [ \begin{array}{c|c|c} 
K'_{01}(\theta) &
K_{11}(\theta)  &
2 K_{02} (\theta)\\
\end{array} \big ]
\sigma_2 
+ \mathcal{O}_2(\sigma_1, \sigma_2),
\end{array}
\]
and recall that $K_{\alpha,m-\alpha}(\theta)$ are 3-dimensional column vectors. Let us introduce the following Taylor expansion
\[
\begin{bmatrix} \nabla\Theta(x) \\ \nabla\Sigma_1(x) \\ \nabla\Sigma_{2}(x) \end{bmatrix} =
\sum_{m=0}^{\infty} \sum_{\alpha=0}^{m}
\begin{bmatrix}  
\nabla \Theta^{(\alpha, m-\alpha)} (\theta)    \\ 
\nabla \Sigma_1^{(\alpha, m-\alpha)} (\theta)    \\ 
\nabla \Sigma_{2}^{(\alpha, m-\alpha)} (\theta)   
\end{bmatrix}
\sigma^{\alpha}_1 \sigma^{m-\alpha}_2,
\]
where we recall that $\nabla \Theta^{(\alpha, m-\alpha)}$ and $\Sigma_i^{(\alpha, m-\alpha)}$, $i=1,2$, are 3-dimensional row vectors.

Thus, substituting the above Taylor expansions in \eqref{eq:prodidentity} for $d=3$ and collecting the $0$-th order terms we have that 
\begin{equation}\label{eq:0inv}
\big [
\begin{array}{c|c|c} 
K'_0 (\theta) & K_{10} (\theta) & K_{01} (\theta) \\
\end{array}
\big ]
\begin{bmatrix} 
\nabla\Theta^{(0)} (\theta)  \\ 
\nabla\Sigma_1^{(0)}(\theta) \\
\nabla\Sigma_{2}^{(0)}(\theta)  
\end{bmatrix} = 
Id_{3x3}.
\end{equation}
From the expression above we have
\[ \nabla\Theta^{(0)} (\theta) = \frac{1}{det(K'_0 (\theta),K_{10} (\theta), K_{01} (\theta))} K_{10}(\theta) \times K_{01} (\theta).\]
That is, the iPRC $\nabla \Theta^{(0)}(\theta)$ is orthogonal to the attracting linear eigenspace of the limit cycle spanned by $K_{10}(\theta)$ and $K_{01} (\theta)$. Moreover,
\[ \nabla\Sigma_1^{(0)} (\theta) = \frac{1}{det(K'_0 (\theta),K_{10} (\theta), K_{01} (\theta))} K_{01}(\theta) \times K'_{0} (\theta),\]
and
\[ \nabla\Sigma_2^{(0)} (\theta) = \frac{1}{det(K'_0 (\theta),K_{10} (\theta), K_{01} (\theta))} K'_{0}(\theta) \times K_{10} (\theta).\]
Thus, using expressions \eqref{eq:eqMzeroSols} for $K_{10}(\theta)$ and $K_{01}(\theta)$ and \eqref{eq:eqMzero} for $K'_{0}(\theta)$, we have an explicit expression for the iPRC and iARCs (see remark \ref{rm:prcsRM}). 

Consider now the terms of order 1 in $\sigma_1$, we have
\[
\big [
\begin{array}{c|c|c} 
K'_0 (\theta) & K_{10} (\theta) & K_{01} (\theta) \\
\end{array}
\big ]
\begin{bmatrix} 
\nabla\Theta^{(1,0)} (\theta)  \\ 
\nabla\Sigma_1^{(1,0)}(\theta) \\
\nabla\Sigma_{2}^{(1,0)}(\theta)  
\end{bmatrix} 
+
\big [
\begin{array}{c|c|c} 
K'_{10}(\theta) &
2K_{20}(\theta)  &
K_{11} (\theta)\\
\end{array}
\big ]
\begin{bmatrix} 
\nabla\Theta^{(0)} (\theta)  \\ 
\nabla\Sigma_1^{(0)}(\theta) \\
\nabla\Sigma_{2}^{(0)}(\theta)  
\end{bmatrix} 
= 0,
\]
and, therefore, using \eqref{eq:0inv} we have
\[
\begin{bmatrix} 
\nabla\Theta^{(1,0)} (\theta)  \\ 
\nabla\Sigma_1^{(1,0)}(\theta) \\
\nabla\Sigma_{2}^{(1,0)}(\theta) 
\end{bmatrix} 
=
-
\begin{bmatrix} 
\nabla\Theta^{(0)} (\theta)  \\ 
\nabla\Sigma_1^{(0)}(\theta) \\
\nabla\Sigma_{2}^{(0)}(\theta)  
\end{bmatrix} 
\big [
\begin{array}{c|c|c} 
K'_{10}(\theta) &
2K_{20}(\theta)  &
K_{11} (\theta)\\
\end{array}
\big ]
\begin{bmatrix} 
\nabla\Theta^{(0)} (\theta)  \\ 
\nabla\Sigma_1^{(0)}(\theta) \\
\nabla\Sigma_{2}^{(0)}(\theta)  
\end{bmatrix} .
\]
Equivalently, collecting terms of order 1 in $\sigma_2$, we have
\[
\begin{bmatrix} 
\nabla\Theta^{(0,1)} (\theta)  \\ 
\nabla\Sigma_1^{(0,1)}(\theta) \\
\nabla\Sigma_{2}^{(0,1)}(\theta)  
\end{bmatrix} 
=
-\begin{bmatrix} 
\nabla\Theta^{(0)} (\theta)  \\ 
\nabla\Sigma_1^{(0)}(\theta) \\
\nabla\Sigma_{2}^{(0)}(\theta)  
\end{bmatrix} 
\big [
\begin{array}{c|c|c} 
K'_{01}(\theta) &
K_{11}(\theta)  &
2 K_{02} (\theta)\\
\end{array}
\big ]
\begin{bmatrix} 
\nabla\Theta^{(0)} (\theta)  \\ 
\nabla\Sigma_1^{(0)}(\theta) \\
\nabla\Sigma_{2}^{(0)}(\theta)  
\end{bmatrix} 
.
\]

Moreover, using the following remark we can compute the functions $\nabla \Theta$ and $\nabla \Sigma_i$, $i=1,2$ at any order $L$.
\begin{remark} \label{rem:allorders}
	If 
	\[ \left ( \sum_{i=0}^{\infty} A_i(\theta) \sigma^i \right ) \left ( \sum_{i=0}^{\infty} B_i(\theta) \sigma^i \right )= Id\]
	where $A_i(\theta)$ and $B_i(\theta)$ are real-valued $d \times d$ matrices, and $\sigma \in \mathbb{R}^d$, we have that,
	\[
	\begin{array}{rcl}
	A_0 B_0 = Id &\Rightarrow & B_0  =  A_0^{-1},\\
	A_0 B_1 + A_1 B_0 =0 & \Rightarrow & B_1  =  - A_0^{-1} A_1 B_0=- B_0 A_1 B_0, \\
	\sum_{j=0}^i A_{i-j} B_j =0 & \Rightarrow & B_i  =  - A_0^{-1} \left ( \sum_{j=0}^{i-1} A_{i-j} B_j \right)=  - B_0 \left ( \sum_{j=0}^{i-1} A_{i-j} B_j \right)\quad \textrm{for } i \geq 2.
	\end{array}
	\]
\end{remark}

\begin{remark}
	The functions $\nabla \Theta^{(0)}$ and $\nabla \Sigma^{(0)}$ correspond to the iPRC and iARFs.
	%functions $Z$ and $\mathcal{I}_i$, respectively, in \cite{wilson2018greater}. The first one, is indeed the iPRC. Moreover, $\nabla \Theta^{(1,0)}$ and $\nabla \Theta^{(0,1)}$ correspond to the functions $C_1$, $C_2$, respectively, while $\nabla \Sigma_i^{(1,0)}$ and $\nabla \Sigma_i^{(0,1)}$ correspond to the functions $D_i^1$ and $D_i^2$ for $i=1,2$, respectively (see formula (47) in \cite{wilson2018greater}).
	%\marginpar{Moreover if we denote by $\bm{v}^{(j)}$ the vector $\bm{v}^{(j)}=(0,\ldots,v_j=1,\ldots,0)$, we have that $\nabla \Theta^{\bm{v}^{(j)}}$ correspond to the functions $C_j$ and $\Sigma_i^{\bm{v}^{(j)}}$ to the functions $D_i^j$.} 
	Notice that, using Remark~\ref{rem:allorders}, our method provides the explicit expansions of $\nabla \Theta$ and $\nabla \Sigma$ in $\sigma$ at any order $L$. Although the method is presented only for the case $d=3$, it extends straightforwardly to any dimension $d$.
\end{remark}

\section{Model Parameters}\label{sec:appendix}
\setcounter{equation}{0}
\subsection*{RT model}

	For the single neuron model \eqref{eq:rtEDOs} we have used the following functions
	\begin{equation*}
	\begin{aligned}
	I_L = g_L(V - V_L), \quad & \quad
	I_{Na} = g_{Na}m^3_{\infty}h(V-V_{Na}),\\
	I_K = g_K(.75(1-h))^4(V-V_K), \quad & \quad
	I_T = g_Tp^2_{\infty}r(V-V_{T}),\\
	h_{\infty}(V) = \frac{1}{1 + \exp((v+41)/4)}, \quad & \quad
	r_{\infty}(V) = \frac{1}{1 + \exp((v+84)/4)},\\
	m_{\infty}(V) = \frac{1}{1 + \exp(-(v+37)/7)}, \quad & \quad
	p_{\infty}(V) = \frac{1}{1 + \exp(-(v+60)/6.2)},\\
	\tau_r(V) = 28 + \exp(-(v+25)/10.5), \quad & \quad \tau_h(V) = 1/(a_h(V) + b_h(V)),\\
	\quad a_h(V) = 0.128 \exp(-(v+46)/18), \quad & \quad
	b_h(V) = \frac{4}{1 + \exp(-(v+23)/5)},
	\end{aligned}
	\end{equation*}
	and the following set of parameters
	\begin{equation*}
	\begin{split}
	\mathcal{P}_{RT} = \{C_m = 1, g_L = 0.05, V_L = -70, g_{Na} = 3, V_{Na} = 50,\\ g_K = 5, V_K = -90, g_T = 5, V_T = 0\}.
	\end{split}
	\end{equation*}
	
\subsection*{HH model}

	For the $HH$ model \eqref{eq:hhEDOs} we have used the following functions
	\begin{equation*}
	I_L = g_L(V - V_L), \quad \quad I_{Na} = g_{Na}m^3_{\infty}h(V-V_{Na}), \quad \quad
	I_K = g_K n^4(V-V_K), 
	\end{equation*}
	\begin{equation*}
	n_{\infty}(V) = \frac{1}{1 + \exp(-(V+53)/15)}, \quad \quad
	h_{\infty}(V) = \frac{1}{1 + \exp((V+62)/7)},
	\end{equation*}
	\begin{equation*}
	m_{\infty}(V) = \frac{1}{1 + \exp(-(V+40)/9)},
	\end{equation*}
	\begin{equation*}
	\tau_h(V) = 7.4 \exp(-((67+V)/20)^2) + 1.2, \quad \quad
	\tau_n(V) = 4.7 \exp(-((79+V)/50)^2) + 1.1,
	\end{equation*}
	and the following set of parameters
	\begin{equation*}
	\mathcal{P}_{HH} = \{ C_m = 1, g_L = 0.1, V_L = -75.6, g_{Na} = 30, V_{Na} = 55, g_K = 9, V_K = -77 \}.
	\end{equation*}
	
\subsection*{WC-Syn Model}

	For the extension of the Wilson-Cowan equations \eqref{eq:wcEDOs3D} we have used the following functions
	\begin{equation*}
	\delta_{E}(x) = \frac{1}{1 + \exp(-a_E(v-\theta_E))}, \quad
	\delta_{I}(x) = \frac{1}{1 + \exp(-a_I(v-\theta_I))},
	\end{equation*}
	and the following set of parameters
	\begin{equation*}
	\begin{split}
	\mathcal{P}_{WC} = \{P = 4.5, \tau_e = 3, a = 8, b = 16, a_E = 3, \theta_E = 4, Q = 0, \\d = 3, \tau_i = 3, c = 7, a_I = 2, \theta_I = 3, \tau_r = 1, \tau_d = 6\}.
	\end{split}
	\end{equation*}	
	
\subsection*{QIF Model}

	For the $QIF$ model \eqref{eq:qfEDOs} we have used the following parameters
	\begin{equation*}
	\mathcal{P}_{QF} = \{ \tau_m = 10, \Delta = 0.3, J = 21, \Theta = 4, \tau_d = 5 \}.
	\end{equation*}

\bibliographystyle{alpha}
%NOLINEAL16 you have to edit the file end.tex, changing bibjimenez to your bibfile name should read as:
\bibliography{bibIsos}

\newcommand{\etalchar}[1]{$^{#1}$}
\begin{thebibliography}{vdBMJR16}

\bibitem[AAG15]{azodi2015phase}
Ramin Azodi-Avval and Alireza Gharabaghi.
\newblock Phase-dependent modulation as a novel approach for therapeutic brain
  stimulation.
\newblock {\em Frontiers in computational neuroscience}, 9:26, 2015.

\bibitem[AC09]{achuthan2009phase}
Srisairam Achuthan and Carmen~C Canavier.
\newblock Phase-resetting curves determine synchronization, phase locking, and
  clustering in networks of neural oscillators.
\newblock {\em Journal of Neuroscience}, 29(16):5218--5233, 2009.

\bibitem[ACN16]{AshComNic2016}
Peter Ashwin, Stephen Coombes, and Rachel Nicks.
\newblock Mathematical frameworks for oscillatory network dynamics in
  neuroscience.
\newblock {\em The Journal of Mathematical Neuroscience}, 6(1):2, Jan 2016.

\bibitem[BMH04]{brown2004phase}
Eric Brown, Jeff Moehlis, and Philip Holmes.
\newblock On the phase reduction and response dynamics of neural oscillator
  populations.
\newblock {\em Neural computation}, 16(4):673--715, 2004.

\bibitem[BMM12]{budivsic2012applied}
Marko Budi{\v{s}}i{\'c}, Ryan Mohr, and Igor Mezi{\'c}.
\newblock Applied koopmanism.
\newblock {\em Chaos: An Interdisciplinary Journal of Nonlinear Science},
  22(4):047510, 2012.

\bibitem[Buz06]{buzsaki2006rhythms}
Gyorgy Buzsaki.
\newblock {\em Rhythms of the Brain}.
\newblock Oxford University Press, 2006.

\bibitem[BY78]{brigham1978fast}
E~Oran Brigham and CK~Yuen.
\newblock The fast fourier transform.
\newblock {\em IEEE Transactions on Systems, Man, and Cybernetics},
  8(2):146--146, 1978.

\bibitem[Can15]{canavier2015phase}
Carmen~C Canavier.
\newblock Phase-resetting as a tool of information transmission.
\newblock {\em Current opinion in neurobiology}, 31:206--213, 2015.

\bibitem[CB19]{coombes2019next}
Stephen Coombes and Aine Byrne.
\newblock Next generation neural mass models.
\newblock In {\em Nonlinear Dynamics in Computational Neuroscience}, pages
  1--16. Springer, 2019.

\bibitem[CFdlL03a]{cabre2003parameterization}
Xavier Cabr{\'e}, Ernest Fontich, and Rafael de~la Llave.
\newblock The parameterization method for invariant manifolds {I}: manifolds
  associated to non-resonant subspaces.
\newblock {\em Indiana University mathematics journal}, pages 283--328, 2003.

\bibitem[CFdlL03b]{cabre2003parameterization2}
Xavier Cabr{\'e}, Ernest Fontich, and Rafael de~la Llave.
\newblock The parameterization method for invariant manifolds {II}: regularity
  with respect to parameters.
\newblock {\em Indiana University mathematics journal}, pages 329--360, 2003.

\bibitem[CFDLL05]{cabre2005parameterization}
Xavier Cabr{\'e}, Ernest Fontich, and Rafael De~La~Llave.
\newblock The parameterization method for invariant manifolds {III}: overview
  and applications.
\newblock {\em Journal of Differential Equations}, 218(2):444--515, 2005.

\bibitem[CG20]{castejon2020phase}
Oriol Castej{\'o}n and Antoni Guillamon.
\newblock Phase-amplitude dynamics in terms of extended response functions:
  Invariant curves and arnold tongues.
\newblock {\em Communications in Nonlinear Science and Numerical Simulation},
  81:105008, 2020.

\bibitem[CGH13]{castejon2013phase}
Oriol Castej{\'o}n, Antoni Guillamon, and Gemma Huguet.
\newblock Phase-amplitude response functions for transient-state stimuli.
\newblock {\em The Journal of Mathematical Neuroscience}, 3(1):13, 2013.

\bibitem[CLMJ15]{castelli2015parameterization}
Roberto Castelli, Jean-Philippe Lessard, and Jason~D Mireles~James.
\newblock Parameterization of invariant manifolds for periodic orbits {I}:
  Efficient numerics via the floquet normal form.
\newblock {\em SIAM Journal on Applied Dynamical Systems}, 14(1):132--167,
  2015.

\bibitem[DDMG16]{Detrixheetal16}
Miles Detrixhe, Marion Doubeck, Jeff Moehlis, and Frederic Gibou.
\newblock A fast eulerian approach for computation of global isochrons in high
  dimensions.
\newblock {\em SIAM Journal on Applied Dynamical Systems}, 15:1501--1527, 2016.

\bibitem[DRM17]{devalle2017firing}
Federico Devalle, Alex Roxin, and Ernest Montbri{\'o}.
\newblock Firing rate equations require a spike synchrony mechanism to
  correctly describe fast oscillations in inhibitory networks.
\newblock {\em PLoS computational biology}, 13(12):e1005881, 2017.

\bibitem[EK91]{ErmentroutKopell91}
G.~B. Ermentrout and N.~Kopell.
\newblock Multiple pulse interactions and averaging in systems of coupled
  neural oscillators.
\newblock {\em J. Math. Biol.}, 29(3):195--217, 1991.

\bibitem[EPW19]{ermentroutetal19}
Bard Ermentrout, Youngmin Park, and Dan Wilson.
\newblock Recent advances in coupled oscillator theory.
\newblock {\em Philosophical Transactions of the Royal Society A: Mathematical,
  Physical and Engineering Sciences}, 377(2160):20190092, 2019.

\bibitem[ET10]{ErmentroutTerman2010}
B.G. Ermentrout and D.H. Terman.
\newblock {\em {Mathematical foundations of neuroscience.}}
\newblock New York : Springer, 2010.

\bibitem[Flo83]{floquet1883equations}
Gaston Floquet.
\newblock Sur les equations differentielles lineaires.
\newblock {\em Ann. ENS [2]}, 12(1883):47--88, 1883.

\bibitem[GER05]{gutkin2005phase}
Boris~S Gutkin, G~Bard Ermentrout, and Alex~D Reyes.
\newblock Phase-response curves give the responses of neurons to transient
  inputs.
\newblock {\em Journal of neurophysiology}, 94(2):1623--1635, 2005.

\bibitem[GH09]{guillamon2009computational}
Antoni Guillamon and Gemma Huguet.
\newblock A computational and geometric approach to phase resetting curves and
  surfaces.
\newblock {\em SIAM Journal on Applied Dynamical Systems}, 8(3):1005--1042,
  2009.

\bibitem[GH13]{guckenheimer2013nonlinear}
John Guckenheimer and Philip Holmes.
\newblock {\em Nonlinear oscillations, dynamical systems, and bifurcations of
  vector fields}, volume~42.
\newblock Springer Science \& Business Media, 2013.

\bibitem[Guc75]{guckenheimer1975}
John Guckenheimer.
\newblock Isochrons and phaseless sets.
\newblock {\em Journal of Mathematical Biology}, 1(3):259--273, 1975.

\bibitem[GW08]{griewank2008evaluating}
Andreas Griewank and Andrea Walther.
\newblock {\em Evaluating derivatives: principles and techniques of algorithmic
  differentiation}, volume 105.
\newblock Siam, 2008.

\bibitem[HCF{\etalchar{+}}16]{haro2016}
{\`A}lex Haro, Marta Canadell, Jordi-Lluis Figueras, Alejandro Luque, and
  Josep-Maria Mondelo.
\newblock {\em The Parameterization Method for Invariant Manifolds}.
\newblock Springer, 2016.

\bibitem[HdlL13]{huguet2013computation}
Gemma Huguet and Rafael de~la Llave.
\newblock Computation of limit cycles and their isochrons: fast algorithms and
  their convergence.
\newblock {\em SIAM Journal on Applied Dynamical Systems}, 12(4):1763--1802,
  2013.

\bibitem[HI12]{hoppensteadt2012}
Frank~C Hoppensteadt and Eugene~M Izhikevich.
\newblock {\em Weakly connected neural networks}, volume 126.
\newblock Springer Science \& Business Media, 2012.

\bibitem[HP70]{hirsch1970stable}
Morris~W Hirsch and Charles~C Pugh.
\newblock Stable manifolds and hyperbolic sets.
\newblock In {\em Global Analysis (Proc. Sympos. Pure Math., Vol. XIV,
  Berkeley, Calif., 1968)}, pages 133--163, 1970.

\bibitem[HPS77]{HirschPS77}
M.W. Hirsch, C.C. Pugh, and M.~Shub.
\newblock {\em Invariant manifolds}, volume 583 of {\em Lecture Notes in Math.}
\newblock Springer-Verlag, Berlin, 1977.

\bibitem[HWS{\etalchar{+}}16]{holt2016phasic}
Abbey~B Holt, Dan Wilson, Max Shinn, Jeff Moehlis, and Theoden~I Netoff.
\newblock Phasic burst stimulation: a closed-loop approach to tuning deep brain
  stimulation parameters for parkinson’s disease.
\newblock {\em PLoS computational biology}, 12(7):e1005011, 2016.

\bibitem[Izh07]{izhikevich2007}
Eugene~M Izhikevich.
\newblock {\em Dynamical systems in neuroscience}.
\newblock MIT press, 2007.

\bibitem[JZ05]{JorbaZ05}
\`Angel Jorba and Maorong Zou.
\newblock A software package for the numerical integration of {ODE}s by means
  of high-order {T}aylor methods.
\newblock {\em Experiment. Math.}, 14(1):99--117, 2005.

\bibitem[KOD{\etalchar{+}}05]{Krauskopfetal05}
B.~Krauskopf, H.~M. Osinga, E.~J. Doedel, M.~E Henderson, J.~Guckenheimer,
  A.~Vladimirsky, M.~Dellnitz, and O.~Junge.
\newblock A survey of methods for computing (un)stable manifolds of vector
  fields.
\newblock {\em International Journal of Bifurcation and Chaos},
  15(03):763--791, 2005.

\bibitem[Kur03]{kuramoto2003chemical}
Yoshiki Kuramoto.
\newblock {\em Chemical oscillations, waves, and turbulence}.
\newblock Courier Corporation, 2003.

\bibitem[Mez05]{mezic2005spectral}
Igor Mezi{\'c}.
\newblock Spectral properties of dynamical systems, model reduction and
  decompositions.
\newblock {\em Nonlinear Dynamics}, 41(1-3):309--325, 2005.

\bibitem[MM12]{mauroy2012use}
Alexandre Mauroy and Igor Mezi{\'c}.
\newblock On the use of fourier averages to compute the global isochrons of
  (quasi) periodic dynamics.
\newblock {\em Chaos: An Interdisciplinary Journal of Nonlinear Science},
  22(3):033112, 2012.

\bibitem[MM16]{mauroy2016global}
Alexandre Mauroy and Igor Mezi{\'c}.
\newblock Global stability analysis using the eigenfunctions of the koopman
  operator.
\newblock {\em IEEE Transactions on Automatic Control}, 61(11):3356--3369,
  2016.

\bibitem[MM18]{mauroy2018global}
Alexandre Mauroy and Igor Mezi{\'c}.
\newblock Global computation of phase-amplitude reduction for limit-cycle
  dynamics.
\newblock {\em Chaos: An Interdisciplinary Journal of Nonlinear Science},
  28(7):073108, 2018.

\bibitem[MM19]{monga2019optimal}
Bharat Monga and Jeff Moehlis.
\newblock Optimal phase control of biological oscillators using augmented phase
  reduction.
\newblock {\em Biological cybernetics}, 113(1-2):161--178, 2019.

\bibitem[MMM13]{mauroy2013isostables}
Alexandre Mauroy, Igor Mezi{\'c}, and Jeff Moehlis.
\newblock Isostables, isochrons, and koopman spectrum for the action--angle
  representation of stable fixed point dynamics.
\newblock {\em Physica D: Nonlinear Phenomena}, 261:19--30, 2013.

\bibitem[MPR15]{montbrio2015macroscopic}
Ernest Montbri{\'o}, Diego Paz{\'o}, and Alex Roxin.
\newblock Macroscopic description for networks of spiking neurons.
\newblock {\em Physical Review X}, 5(2):021028, 2015.

\bibitem[MWMM19]{monga2019phase}
Bharat Monga, Dan Wilson, Tim Matchen, and Jeff Moehlis.
\newblock Phase reduction and phase-based optimal control for biological
  systems: a tutorial.
\newblock {\em Biological cybernetics}, 113(1-2):11--46, 2019.

\bibitem[OM10]{osinga2010continuation}
Hinke~M Osinga and Jeff Moehlis.
\newblock Continuation-based computation of global isochrons.
\newblock {\em SIAM Journal on Applied Dynamical Systems}, 9(4):1201--1228,
  2010.

\bibitem[OPC04]{OprisanPC04}
S.A. Oprisan, A.A. Prinz, and C.C. Canavier.
\newblock Phase resetting and phase locking in hybrid circuits of one model and
  one biological neuron.
\newblock {\em Biophysical Journal}, 87(4):2283 -- 2298, 2004.

\bibitem[PCHS17]{PerezCerveraHS17}
A.~P\'erez-Cervera, G.~Huguet, and T.M. Seara.
\newblock Computation of invariant curves in the analysis of periodically
  forced neural oscillators.
\newblock In {\em Nonlinear Systems, Vol. 2: Nonlinear Phenomena in Biology,
  Optics and Condensed Matter}. Springer, 2017.

\bibitem[PCMSH19]{PerezCervera2019}
Alberto P{\'e}rez-Cervera, Tere M-Seara, and Gemma Huguet.
\newblock A geometric approach to phase response curves and its numerical
  computation through the parameterization method.
\newblock {\em Journal of Nonlinear Science}, Jul 2019.

\bibitem[PCSH20]{PerezSH20}
Alberto Pérez-Cervera, Tere~M. Seara, and Gemma Huguet.
\newblock Phase-locked states in oscillating neural networks and their role in
  neural communication.
\newblock {\em Communications in Nonlinear Science and Numerical Simulation},
  80:104992, 2020.

\bibitem[PE16]{park2016weakly}
Youngmin Park and Bard Ermentrout.
\newblock Weakly coupled oscillators in a slowly varying world.
\newblock {\em Journal of computational neuroscience}, 40(3):269--281, 2016.

\bibitem[PN15]{pyragas2015phase}
Kestutis Pyragas and Viktor Novi{\v{c}}enko.
\newblock Phase reduction of a limit cycle oscillator perturbed by a strong
  amplitude-modulated high-frequency force.
\newblock {\em Physical Review E}, 92(1):012910, 2015.

\bibitem[PRK01]{PIK01}
A.~Pikovsky, M.~G. Rosenblum, and J.~Kurths.
\newblock {\em Synchronization, A Universal Concept in Nonlinear Sciences}.
\newblock Cambridge University Press, 2001.

\bibitem[RE98]{rinzel1998analysis}
John Rinzel and G~Bard Ermentrout.
\newblock Analysis of neural excitability and oscillations.
\newblock {\em Methods in neuronal modeling}, 2:251--292, 1998.

\bibitem[RP04]{rosenblum2004delayed}
Michael Rosenblum and Arkady Pikovsky.
\newblock Delayed feedback control of collective synchrony: An approach to
  suppression of pathological brain rhythms.
\newblock {\em Physical Review E}, 70(4):041904, 2004.

\bibitem[RP19]{rosenblum2019numerical}
Michael Rosenblum and Arkady Pikovsky.
\newblock Numerical phase reduction beyond the first order approximation.
\newblock {\em Chaos: An Interdisciplinary Journal of Nonlinear Science},
  29(1):011105, 2019.

\bibitem[RT04]{rubin2004high}
Jonathan~E Rubin and David Terman.
\newblock High frequency stimulation of the subthalamic nucleus eliminates
  pathological thalamic rhythmicity in a computational model.
\newblock {\em Journal of computational neuroscience}, 16(3):211--235, 2004.

\bibitem[SD10]{suvak2010quadratic}
{\"O}nder Suvak and Alper Demir.
\newblock Quadratic approximations for the isochrons of oscillators: a general
  theory, advanced numerical methods, and accurate phase computations.
\newblock {\em IEEE Transactions on Computer-Aided Design of Integrated
  Circuits and Systems}, 29(8):1215--1228, 2010.

\bibitem[SEW10]{smeal2010phase}
Roy~M Smeal, G~Bard Ermentrout, and John~A White.
\newblock Phase-response curves and synchronized neural networks.
\newblock {\em Philosophical Transactions of the Royal Society of London B:
  Biological Sciences}, 365(1551):2407--2422, 2010.

\bibitem[Sim90]{simo1990analytical}
Carles Sim{\'o}.
\newblock On the analytical and numerical approximation of invariant manifolds.
\newblock In {\em Les M{\'e}thodes Modernes de la M{\'e}canique C{\'e}leste.
  Modern methods in celestial mechanics}, pages 285--329, 1990.

\bibitem[SKN17]{shirasaka2017phase}
Sho Shirasaka, Wataru Kurebayashi, and Hiroya Nakao.
\newblock Phase-amplitude reduction of transient dynamics far from attractors
  for limit-cycling systems.
\newblock {\em Chaos: An Interdisciplinary Journal of Nonlinear Science},
  27(2):023119, 2017.

\bibitem[SP13]{schwabedal2013phase}
Justus~TC Schwabedal and Arkady Pikovsky.
\newblock Phase description of stochastic oscillations.
\newblock {\em Physical review letters}, 110(20):204102, 2013.

\bibitem[SPB11]{schultheiss2011phase}
Nathan~W Schultheiss, Astrid~A Prinz, and Robert~J Butera.
\newblock {\em Phase response curves in neuroscience: theory, experiment, and
  analysis}.
\newblock Springer Science \& Business Media, 2011.

\bibitem[Str94]{strogatzbook}
Steven~H. Strogatz.
\newblock {\em {Nonlinear Dynamics and Chaos: With Applications to Physics,
  Biology, Chemistry and Engineering}}.
\newblock Westview Press, 1994.

\bibitem[Tas07]{tass2007phase}
Peter~A Tass.
\newblock {\em Phase resetting in medicine and biology: stochastic modelling
  and data analysis}.
\newblock Springer Science \& Business Media, 2007.

\bibitem[TF10]{takeshita2010higher}
Daisuke Takeshita and Renato Feres.
\newblock Higher order approximation of isochrons.
\newblock {\em Nonlinearity}, 23(6):1303, 2010.

\bibitem[TL14]{thomas2014asymptotic}
Peter~J Thomas and Benjamin Lindner.
\newblock Asymptotic phase for stochastic oscillators.
\newblock {\em Physical review letters}, 113(25):254101, 2014.

\bibitem[TS10]{tiesinga2010}
Paul~HE Tiesinga and Terrence~J Sejnowski.
\newblock Mechanisms for phase shifting in cortical networks and their role in
  communication through coherence.
\newblock {\em Frontiers in human neuroscience}, 4:196, 2010.

\bibitem[vdBMJR16]{James16}
Jan~Bouwe van~den Berg, Jason~D. Mireles~James, and Christian Reinhardt.
\newblock Computing (un)stable manifolds with validated error bounds:
  non-resonant and resonant spectra.
\newblock {\em J. Nonlinear Sci.}, 26(4):1055--1095, 2016.

\bibitem[WC72]{wilson1972excitatory}
Hugh~R Wilson and Jack~D Cowan.
\newblock Excitatory and inhibitory interactions in localized populations of
  model neurons.
\newblock {\em Biophysical journal}, 12(1):1--24, 1972.

\bibitem[WE18]{wilson2018greater}
Dan Wilson and Bard Ermentrout.
\newblock Greater accuracy and broadened applicability of phase reduction using
  isostable coordinates.
\newblock {\em Journal of mathematical biology}, 76(1-2):37--66, 2018.

\bibitem[WE19]{WilsonE19}
Dan Wilson and Bard Ermentrout.
\newblock Augmented phase reduction of (not so) weakly perturbed coupled
  oscillators.
\newblock {\em SIAM Review}, 61(2):277--315, 2019.

\bibitem[Wil20]{wilson2020}
Dan Wilson.
\newblock Phase-amplitude reduction far beyond the weakly perturbed paradigm.
\newblock {\em Physical Review E}, 101, 2020.

\bibitem[Win67]{winfree1967biological}
Arthur~T Winfree.
\newblock Biological rhythms and the behavior of populations of coupled
  oscillators.
\newblock {\em Journal of theoretical biology}, 16(1):15--42, 1967.

\bibitem[Win74]{winfree1974patterns}
AT~Winfree.
\newblock Patterns of phase compromise in biological cycles.
\newblock {\em Journal of Mathematical Biology}, 1(1):73--93, 1974.

\bibitem[Win01]{winfree2001geometry}
Arthur~T Winfree.
\newblock {\em The geometry of biological time}, volume~12.
\newblock Springer Science \& Business Media, 2001.

\bibitem[WLTC13]{wedgwood2013phase}
Kyle~CA Wedgwood, Kevin~K Lin, Ruediger Thul, and Stephen Coombes.
\newblock Phase-amplitude descriptions of neural oscillator models.
\newblock {\em The Journal of Mathematical Neuroscience}, 3(1):2, 2013.

\bibitem[WM16]{moehliswilsonpre2016}
Dan Wilson and Jeff Moehlis.
\newblock Isostable reduction of periodic orbits.
\newblock {\em Physical Review E}, 94(5):052213, 2016.

\end{thebibliography}
	
\end{document}